\documentclass[11pt]{article}
\usepackage{graphicx}
\usepackage{subfigure}
\usepackage{caption}
\usepackage[hidelinks]{hyperref}
\usepackage[margin=1in]{geometry}

\usepackage{algpseudocode}
\usepackage{amsmath}
\usepackage{graphicx}
\usepackage{enumerate}
\usepackage[authoryear]{natbib}
\usepackage{url}
\usepackage{enumitem}

\usepackage{amsmath,amssymb,amsthm,bm,mathtools}
\usepackage{algorithm}
\usepackage{dsfont,multirow,hyperref,setspace,enumerate}
\hypersetup{colorlinks,linkcolor=blue,citecolor=blue,urlcolor={black}}
\usepackage{lscape}
\usepackage{afterpage}
\usepackage{hyperref}

\usepackage{microtype}
\usepackage{wrapfig}
\allowdisplaybreaks
\usepackage{algorithm}
\usepackage{algpseudocode}

\usepackage{textcomp}
\usepackage{stfloats}
\usepackage{verbatim}
\usepackage{array}

\usepackage{tabularx,booktabs}
\newcolumntype{Y}{>{\centering\arraybackslash}X}

\theoremstyle{definition}
\newtheorem{thm}{Theorem}
\newtheorem{lem}{Lemma}
\newtheorem{prop}{Proposition}[section]

\newtheorem{cor}{Corollary}

\theoremstyle{definition}
\newtheorem{assumption}{Assumption}
\newtheorem{defn}{Definition}
\newtheorem{example}{Example}
\newtheorem{rmk}{Remark}
\newtheorem{condition}{Condition}
\newtheorem{conjecture}{Conjecture}

\newtheorem{innercustomgeneric}{\customgenericname}
\providecommand{\customgenericname}{}
\newcommand{\newcustomtheorem}[2]{%
  \newenvironment{#1}[1]
  {%
   \renewcommand\customgenericname{#2}%
   \renewcommand\theinnercustomgeneric{##1}%
   \innercustomgeneric
  }
  {\endinnercustomgeneric}
}

\newcustomtheorem{customexample}{Example}

\usepackage{appendix}
\usepackage{fullpage} 
\usepackage{microtype}
\usepackage{wrapfig}
\mathtoolsset{showonlyrefs}

\newcommand{\of}[1]{\left(#1\right)}
\newcommand{\off}[1]{\left[#1\right]}
\newcommand{\offf}[1]{\left\{#1\right\}}
\newcommand{\aabs}[1]{\left|#1\right|}
\newcommand{\ang}[1]{\left\langle#1\right\rangle}

\newcommand{\Mat}{\text{Mat}}

\usepackage{booktabs}

\allowdisplaybreaks

\newcommand*{\KeepStyleUnderBrace}[1]{
  \mathop{%
    \mathchoice
    {\underbrace{\displaystyle#1}}%
    {\underbrace{\textstyle#1}}%
    {\underbrace{\scriptstyle#1}}%
    {\underbrace{\scriptscriptstyle#1}}%
  }\limits
}

\usepackage{stmaryrd}

\def\mepsilon{\boldsymbol{\epsilon}}
\def\mtheta{\boldsymbol{\theta}}
\def\mTheta{\boldsymbol{\Theta}}
\def\bbR{\mathbb{R}}
\def\bbP{\mathbb{P}}
\def\bbE{\mathbb{E}}
\def\ind{\mathds{1}}
\def\mone{\bf 1}
\def\mat{\text{Mat}}
\def\pad{\text{Pad}}

\def\ma{\bm{a}}
\def\mb{\bm{b}}

\def\me{\bm{e}}

\def\mp{\bm{p}}

\def\ms{\bm{s}}

\def\mv{\bm{v}}
\def\mw{\bm{w}}
\def\mx{\bm{x}}
\def\my{\bm{y}}

\def\mA{\bm{A}}
\def\mB{\bm{B}}
\def\mC{\bm{C}}
\def\mD{\bm{D}}
\def\mE{\bm{E}}

\def\mI{\bm{I}}

\def\mI{\bm{I}}
\def\mM{\bm{M}}

\def\mP{\bm{P}}
\def\mQ{\bm{Q}}
\def\mR{\bm{R}}
\def\mS{\bm{S}}
\def\mT{\bm{T}}
\def\mU{\bm{U}}
\def\mV{\bm{V}}
\def\mW{\bm{W}}
\def\mX{\bm{X}}
\def\mY{\bm{Y}}
\def\mZ{\bm{Z}}

\def\mone{\mathbf{1}}
\def\mx{\bm x}
\def\mZ{\bm Z}
\def\mA{\bm A}
\def\mB{\bm B}
\def\mS{\bm S}
\def\mC{\bm C}
\def\mI{\bm I}
\def\mT{\bm T}

\def\mX{\bm X}
\def\mY{\bm Y}
\def\mM{\bm M}

\def\tZ{\mathcal{Z}}
\def\tA{\mathcal{A}}

\def\tE{\mathcal{E}}

\def\tH{\mathcal{H}}

\def\tN{\mathcal{N}}
\def\tO{\mathcal{O}}
\def\tP{\mathcal{P}}
\def\tQ{\mathcal{Q}}

\def\tS{\mathcal{S}}
\def\tT{\mathcal{T}}

\def\tX{\mathcal{X}}
\def\tY{\mathcal{Y}}
\def\tZ{\mathcal{Z}}

\def\entry#1{\llbracket #1 \rrbracket}

\def\entry#1{\llbracket #1 \rrbracket}

\def\bbR{\mathbb{R}}

\newcommand{\onorm}[1]{\left\lVert#1\right\rVert}

\newcommand{\onormSize}[2]{#1\lVert#2#1\rVert}

\newcommand{\FnormSize}[2]{#1\lVert#2#1\rVert_F}

\DeclareMathOperator*{\argmin}{arg\,min}
\DeclareMathOperator*{\argmax}{arg\,max}

\usepackage[parfill]{parskip}
\usepackage{algorithm}
\usepackage{algpseudocode}
\usepackage{amsmath}

\DeclarePairedDelimiter{\ceil}{\lceil}{\rceil}
\DeclarePairedDelimiter{\floor}{\lfloor}{\rfloor}

\algnewcommand\algorithmicinput{\textbf{Input:}}
\algnewcommand\algorithmicoutput{\textbf{Output:}}
\algnewcommand\INPUT{\item[\algorithmicinput]}
\algnewcommand\OUTPUT{\item[\algorithmicoutput]}

\newcommand\Algphase[1]{%
\vspace*{-.7\baselineskip}\Statex\hspace*{\dimexpr-\algorithmicindent-2pt\relax}\rule{\textwidth}{0.4pt}%
\Statex\hspace*{-\algorithmicindent}\textbf{#1}%
\vspace*{-.7\baselineskip}\Statex\hspace*{\dimexpr-\algorithmicindent-2pt\relax}\rule{\textwidth}{0.4pt}%
}
\usepackage{caption}

\title{Multiway Spherical Clustering via Degree-Corrected\\ Tensor Block Models\footnote{This paper was presented in part at 25th International Conference on Artificial Intelligence and Statistics (AISTATS).}}
\date{}
\author{%
Jiaxin Hu \\
University of Wisconsin -- Madison\\
\texttt{jhu267@wisc.edu} \\
\and
Miaoyan Wang \\
University of Wisconsin -- Madison\\
\texttt{miaoyan.wang@wisc.edu}\\
}

\begin{document}

\maketitle

\begin{abstract}
We consider the problem of multiway clustering in the presence of unknown degree heterogeneity. Such data problems arise commonly in applications such as recommendation system, neuroimaging, community detection, and hypergraph partitions in social networks. The allowance of degree heterogeneity provides great flexibility in clustering models, but the extra complexity poses significant challenges in both statistics and computation. Here, we develop a degree-corrected tensor block model with estimation accuracy guarantees. We present the phase transition of clustering performance based on the notion of angle separability, and we characterize three signal-to-noise regimes corresponding to different statistical-computational behaviors. In particular, we demonstrate that an intrinsic statistical-to-computational gap emerges only for tensors of order three or greater. Further, we develop an efficient polynomial-time algorithm that provably achieves exact clustering under mild signal conditions. The efficacy of our procedure is demonstrated through two data applications, one on human brain connectome project, and another on Peru Legislation network dataset. 
\end{abstract}


\section{Introduction}

Multiway arrays have been widely collected in various fields including social networks \citep{anandkumar2014tensor}, neuroscience \citep{wang2017bayesian}, and computer science \citep{koniusz2016sparse}. Tensors effectively represent the multiway data and serve as the foundation in higher-order data analysis. One data example is from multi-tissue multi-individual gene expression study~\citep{wang2019three,hore2016tensor}, where the data tensor consists of expression measurements indexed by (gene, individual, tissue) triplets. Another example is \emph{hypergraph} network~\citep{ghoshdastidar2017uniform,ghoshdastidar2017consistency,ahn2019community,ke2019community} in social science. A $K$-uniform hypergraph can be naturally represented as an order-$K$ tensor, where each entry indicates the presence of $K$-way hyperedge among nodes (a.k.a.\ entities). In both examples, identifying the similarity among tensor entities is important for scientific discovery. 

\begin{figure}[t]
    \centering
    \includegraphics[width = .9\textwidth]{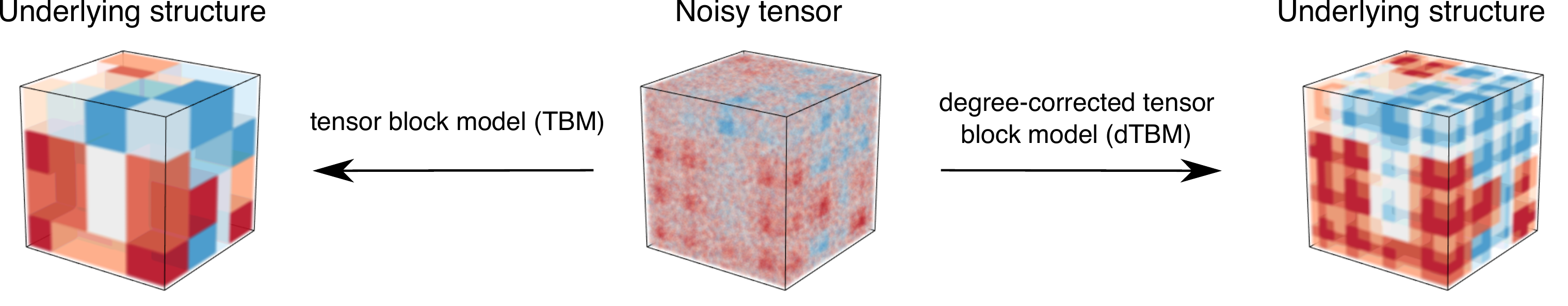}
    \caption{Examples for order-3 tensor block model (TBM) with and without degree correction. Both TBM and dTBM have four communities on each mode, while dTBM allows a richer structure with degree heterogeneity.
    }
    \label{fig:intro}
\end{figure}

We study the problem of multiway clustering based on a data tensor. The goal of multiway clustering is to identify a checkerboard structure from a noisy data tensor. Figure~\ref{fig:intro} illustrates the noisy tensor and the underlying checkerboard structures discovered by multiway clustering methods. In the hypergraph example, the multiway clustering aims to identify the underlying block partition of nodes based on their higher-order connectivities; therefore, we also refer to the clustering as \emph{higher-order clustering}. The most common model for higher-order clustering is called \emph{tensor block model} (TBM)~\citep{wang2019multiway}, which extends the usual matrix stochastic block model~\citep{abbe2018community} to tensors. The matrix analysis tools, however, are sub-optimal for higher-order clustering. Developing tensor tools for solving block models has received increased interest recently~\citep{ wang2019multiway,chi2020provable,han2022exact}. 

The classical tensor block model suffers from drawbacks to model real world data in spite of the popularity. The key underlying assumption of block model is that all nodes in the same community are exchangeable; i.e., the nodes have no individual-specific parameters apart from the community-specific parameters. However, the exchangeability assumption is often non-realistic. Each node may contribute to the data variation by its own multiplicative effect. We call the unequal node-specific effects the \emph{degree heterogeneity}. Such degree heterogeneity appears commonly in social networks. Ignoring the degree heterogeneity may seriously mislead the clustering results. For example, the regular block model fails to model the member affiliation in the Karate Club network~\citep{bickel2009nonparametric} without addressing degree heterogeneity. 

The \emph{degree-corrected tensor block model} (dTBM) has been proposed recently to account for the degree heterogeneity~\citep{ke2019community}. The dTBM combines a higher-order checkerboard structure with degree parameter $\mtheta=(\mtheta(1),\ldots,\mtheta(p))^T$ to allow heterogeneity among $p$ nodes.  Figure~\ref{fig:intro} compares the underlying structures of TBM and dTBM with the same number of communities. The dTBM allows varying values within the same community, thereby allowing a richer structure. To solve dTBM, we project clustering objects to a unit sphere and perform iterative clustering based on angle similarity. We refer to the algorithm as the \textit{spherical} clustering; detailed procedures are in Section~\ref{sec:alg}. The spherical clustering avoids the estimation of nuisance degree heterogeneity. The usage of angle similarity brings new challenges to the theoretical results, and we develop new polar-coordinate based techniques in the proofs. 

\begin{figure}[t]
    \centering
    \includegraphics[width = \textwidth]{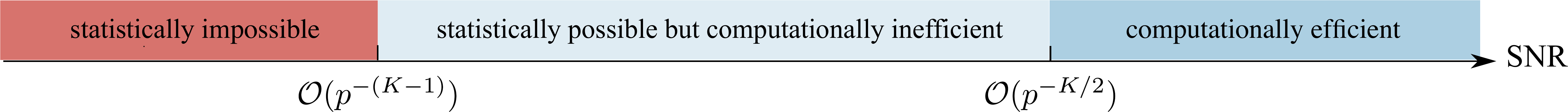}
    \caption{SNR thresholds for statistical and computational limits in order-$K$ dTBM with dimension $(p,...,p)$ and $K \geq 2$. The SNR gap between statistical possibility and computational efficiency  exists only for tensors with $K \geq 3$. }
    \label{fig:phase_axis}
\end{figure}

{\bf Our contributions.} The primary goal of this paper is to provide both statistical and computational guarantees for dTBM. Our main contributions are summarized below.
\begin{itemize}[leftmargin=*]

 \item We develop a general dTBM and establish the identifiability for the uniqueness of clustering using the notion of angle separability.
 
\item  We present the phase transition of clustering performance with respect to three different statistical and computational behaviors.  We characterize, for the first time, the critical signal-to-noise (SNR) thresholds in dTBMs, revealing the intrinsic distinctions among (vector) one-dimensional clustering, (matrix) biclustering, and (tensor) higher-order clustering. Specific SNR thresholds and algorithm behaviors are depicted in  Figure~\ref{fig:phase_axis}. 
        
 \item We provide an angle-based algorithm that achieves exact clustering \emph{in polynomial time} under mild conditions. Simulation and data studies demonstrate that our algorithm outperforms existing higher-order clustering algorithms. 
\end{itemize}
The last two contributions, to our best knowledge, are new to the literature of dTBMs. 

{\bf Related work.} 
Our work is closely related to but also distinct from several lines of existing research. Table~\ref{tab:comp} summarizes the most relevant models. 

\begin{itemize}[wide]
    \item \textit{Block model for clustering.} The block model such as stochastic block model (SBM) and degree-corrected SBM has been widely used for matrix clustering problems. The theoretical properties and algorithm performance for matrix block models have been well-studied~\citep{gao2018community}; see the review paper~\citep{abbe2018community} and the references therein. However, The tensor counterparts are relatively less understood. 
    
    \item  \textit{Tensor block model.} The (non-degree) tensor block model (TBM) is a higher-order extension of SBM, and its statistical-computational properties are investigated in recent literatures~\citep{wang2019multiway, han2022exact, ghoshdastidar2017consistency}. Some works~\citep{ahn2018hypergraph} study the TBM with sparse observations, while, others~\citep{wang2019multiway, han2022exact} and our work focus on the dense regime. Extending results from non-degree to degree-corrected model is highly challenging. Our dTBM parameter space is equipped with angle-based similarity and nuisance degree parameters. The extra complexity makes the Cartesian coordinates based analysis~\citep{han2022exact} non-applicable to our setting. Towards this goal, we have developed a new polar coordinates based analysis to control the model complexity. We have also developed a new angle-based iteration algorithm to achieve optimal clustering rates \emph{without the need of estimating nuisance degree parameters}.

    \item \textit{Degree-corrected block model.} The hypergraph degree-corrected block model (hDCBM) and its variant have been proposed in the literature~\citep{ke2019community, yuan2022testing}. For this popular model, however, the optimal statistical-computational rates remain an open problem. Our main contribution is to provide a sharp statistical and computational critical phase transition in dTBM literature. In addition, our algorithm results in a faster \emph{exponential} error rate, in contrast to the \emph{polynomial} rate in~\cite{ke2019community}. The original hDCBM~\citep{ke2019community} is designed for binary observations only, and we extend the model to both continuous and binary observations. We believe our results are novel and helpful to the community. See Figure~\ref{fig:phase_axis} for overview of our results. 
    
    \item \textit{Global-to-local algorithm strategy.} Our methods generalize the recent global-to-local strategy for matrix learning~\citep{gao2018community,chi2019nonconvex,yun2016optimal} to tensors~\citep{han2022exact,ahn2018hypergraph,kim2018stochastic}. Despite the conceptual similarity, we address several fundamental challenges associated with this non-convex, non-continuous problem. We show the insufficiency of the conventional tensor HOSVD~\citep{de2000multilinear}, and we develop a weighted higher-order initialization that relaxes the singular-value gap separation condition. Furthermore, our local iteration leverages the angle-based clustering in order to avoid explicit estimation of degree heterogeneity. Our bounds reveal the interesting interplay between the computational and statistical errors. We show that our final estimate \emph{provably} achieves the exact clustering within only polynomial-time complexity. 
    
\end{itemize}

\begin{table}[t]
\resizebox{\textwidth}{!}{%
    \begin{tabular}{c|cccccc}
\hline
    & \cite{gao2018community}& \cite{ahn2018hypergraph} &\cite{han2022exact}& \cite{ghoshdastidar2017consistency} &\cite{ke2019community} & \textbf{Ours}\\
    \hline
      Allow tensors of arbitrary order & $\times$ & $\surd$&  $\surd$ & $\surd$ &$\surd$ & $\surd$  \\
        Allow degree heterogeneity &  $\surd$ & $\times$ & $\times$ & $\surd$ &$\surd$ & $\surd$ \\
        Singular-value gap-free clustering &  $\surd$ & $\surd$ & $\surd$ & $\times$ & $\times$ &$\surd$ \\
        Misclustering rate (for order $K^*$)& - & ${p^{-(K-1)} \alpha^{-1}}^{**}$ & $\exp(-p^{K/2})$ & $p^{-1}$& $p^{-2}$ &  $\exp(-p^{K/2})$\\
        Consider sparse observation & $\times$ & $\surd$ & $\times$ & $\times$ & $\times$ &$\times$\\
        \hline
    \end{tabular}
    }
    \caption{ Comparison between previous methods with our method. $^*$We list the result for order-K tensors with $K \geq 3$ and general number of communities $r = \tO(1)$. $^{**}$The parameter $\alpha = f(p) > 0$ denotes the sparsity level which is some function of dimension $p$. }\label{tab:comp}
\end{table}

{\bf Notation.} We use lower-case letters (e.g., $a,b$) for scalars, lower-case boldface letters (e.g., $\ma,\mtheta$) for vectors, upper-case boldface letters (e.g., $\mX,\mY$) for matrices, and calligraphy letters (e.g., $\tX,\tY$) for tensors of order three or greater. We use $\mone_p$ to denote a vector of length $p$ with all entries to be 1. We use $|\cdot|$ for the cardinality of a set and $\ind\{\cdot\}$ for the indicator function. For an integer $p\in\mathbb{N}_{+}$, we use the shorthand $[p]= \offf{1,2,...,p}$. For a length-$p$ vector $\ma$, we use $a(i)\in\mathbb{R}$ to denote the $i$-th entry of $\ma$, and use $\ma_{I}$ to denote the sub-vector by restricting the indices in the set $I\subset [p]$.  We use  $\onorm{\ma}=\sqrt{\sum_{i}a^2(i)}$ to denote the $\ell_2$-norm, $\onorm{\ma}_1=\sum_i |a_i|$ to denote the $\ell_1$ norm of $\ma$. For two vector $\ma, \mb$ of the same dimension, we denote the angle between $\ma, \mb$ by 
\begin{equation}
    \cos \of{\ma, \mb} = \frac{\ang{ \ma, \mb}}{ \onorm{\ma} \onorm{\mb} },
\end{equation}
where $\ang{\ma,\mb}$ is the inner product of two vectors and $\cos \of{\ma, \mb} \in [-1,1]$. We make the convention that $\cos \of{\ma, \mb} = \cos \of{\ma^T, \mb^T}$. 

Let $\tY  \in \bbR^{p_1 \times \cdots \times p_K}$ be an order-$K$ $(p_1,...,p_K)$-dimensional tensor. We use $\tY(i_1,\ldots,i_K)$ to denote the $(i_1,\ldots,i_K)$-th entry of $\tY$. The multilinear multiplication of a tensor $\tS\in \bbR^{r_1\times \cdots \times r_K}$ by matrices $\mM_k \in\mathbb{R}^{p_k\times r_k}$ results in an order-$K$ $(p_1,\ldots,p_K)$-dimensional tensor $\tX$, denoted
\[
\tX=\tS \times_1 \mM_1 \times \cdots \times_K \mM_K,
\]
where the entries of $\tX$ are defined by
\begin{align}
    &\tX(i_1,\ldots, i_K)  =\sum_{(j_1,\ldots,j_K)}\tS(j_1,\ldots,j_K)\mM_1(i_1,j_1)\cdots \mM_K(i_K,j_K).
\end{align} 
\normalsize
For a matrix $\mY$, we use $\mY_{i:}$ (respectively, $\mY_{:i}$) to denote the $i$-th row (respectively, $i$-th column) of the matrix. Similarly, for an order-3 tensor, we use $\tY_{::i}$ to denote the $i$-th matrix slide of the tensor. We use $\text{Ave}(\cdot)$ to denote the operation of taking averages across elements and $\text{Mat}_k(\cdot)$ to denote the unfolding operation that reshapes the tensor along mode $k$ into a matrix.   {For a symmetric tensor $\tX \in\mathbb{R}^{p\times \cdots \times p}$, we omit the subscript and use $\text{Mat}(\tX)\in\mathbb{R}^{p\times p^{K-1}}$ to denote the unfolding.}  For two sequences $\{a_p\}, \{b_p\}$, we denote $a_p\lesssim b_p$ or $a_p=\tO(b_p)$ if $\lim_{p\to\infty}a_p /b_p\leq c${, $a_p \gtrsim b_p$ or $a_p = \Omega(b_p)$ if $\lim_{p\to\infty}a_p /b_p\geq c$, for some constant $c > 0$,}  $a_p=o(b_p)$ if $\lim_{p\to\infty}a_p/b_p =0$, and $a_p \asymp b_p$ if both $b_p \lesssim a_p$ and $a_p\lesssim b_p$. Throughout the paper, we use the terms ``community'' and ``clusters'' exchangeably.

{\bf Organization.} The rest of this paper is organized as follows. Section~\ref{sec:model} introduces the degree-corrected tensor block model (dTBM) with three motivating examples and presents the identifiability of dTBM under the angle gap condition. We show the phase transition and the existence of statistical-computational gaps for the higher-order dTBM in Section~\ref{sec:limits}. In Section~\ref{sec:alg}, we provide a polynomial-time two-stage algorithm with misclustering rate guarantees. { Extension to Bernoulli models is also presented. In Section~\ref{sec:tbm}, we compare our work with non-degree tensor block models.} Numerical studies including the simulation, comparison with other methods, and two real dataset analyses are in Sections~\ref{sec:simulation}-\ref{sec:real}. The main technical ideas we develop for addressing main theorems are provided in Section~\ref{sec:mainproof}. Detailed proofs and extra theoretical results are provided in Appendix.

\section{Model formulation and motivations}\label{sec:model}

\subsection{Degree-corrected tensor block model}

Suppose that we have an order-$K$ data tensor $\tY \in \bbR^{p \times \cdots \times p}$. 
Assume that there exist $r \geq 1$  disjoint communities among the $p$ nodes. We represent the community assignment by a function $z \colon [p]\mapsto[r]$, where $z(i) = a$ for $i$-th node that belongs to the $a$-th community. Then, $z^{-1}(a)=\{i\in[p]\colon z(i)=a\}$ denotes the set of nodes that belong to the $a$-th community, and $|z^{-1}(a)|$ denotes the number of nodes in the $a$-th community. Let $\mtheta=(\theta(1),\ldots,\theta(p))^T$ denote the degree heterogeneity for $p$ nodes. We consider the order-$K$ dTBM~\citep{ghoshdastidar2017consistency,ke2019community},
\begin{equation}\label{eq:model_margin}
    \tY(i_1,\ldots,i_K) = \tS(z(i_1),\ldots,z(i_K)) \prod_{k = 1}^K\theta_{i_k} + \tE(i_1,\ldots,i_K), 
\end{equation}
\normalsize
where $\tS \in \bbR^{r \times \cdots \times r}$ is an order-$K$ tensor collecting the block means among communities, and 
$\tE\in\mathbb{R}^{p\times \cdots \times p}$ is a noise tensor consisting of independent zero-mean sub-Gaussian entries with variance bounded by $\sigma^2$. 
The unknown parameters are $z$, $S$, and $\mtheta$. The dTBM can be equivalently written in a compact form of tensor-matrix product:
\begin{equation}\label{eq:model_tensor}
\bbE\tY = \tS \times_1 \mTheta \mM \times_2 \cdots \times_K  \mTheta \mM,
\end{equation}
where $\mTheta = \text{diag}(\theta(1),...,\theta(p)) \in \bbR^{p \times p}$ is a diagonal matrix, $\mM \in \offf{0,1}^{p \times r}$ is the membership matrix associated with community assignment $z$ such that $\mM(i,j)=\ind\{z(i)=j\}$. By definition, each row of $\mM$ has one copy of 1's and 0's elsewhere. Note that the discrete nature of $\mM$ renders our model~\eqref{eq:model_tensor} more challenging than Tucker decomposition. We call a tensor $\tY$ an $r$-block tensor with degree $\mtheta$ if $\tY$ admits dTBM~\eqref{eq:model_tensor}   {and let $\tX = \bbE \tY$ denote the mean tensor.} The goal of clustering is to estimate $z$ from a single noisy tensor $\tY$. We are particularly interested in the high-dimensional regime where $p$ grows whereas $r=\tO(1)$. 

For ease of notation, we have focused on the case with symmetric mean tensor $\bbE \tY$. This assumption simplifies the notation because all modes have the same $(\mTheta, \mM, z)$; the noise tensor $\tE$ and the data tensor $\tY$ are still possibly asymmetric. In general, we allow asymmetric mean tensors with $\{(\mTheta_k, \mM_k, z_k)\}_{k=1}^K$, one for each mode. The extension can be found in Appendix B.

\subsection{Motivating examples}\label{subsec:motiv} Here, we provide four applications to illustrate the practical necessity of dTBM.

\paragraph{Tensor block model} Consider the model~\eqref{eq:model_tensor}. Let $\theta(i)=1$ for all $ i \in [p]$. The model~\eqref{eq:model_tensor} reduces to the tensor block model, which is widely used in previous clustering algorithms~\citep{wang2019multiway,chi2020provable,han2022exact}. The theoretical results in TBM serve as benchmarks for dTBM.  

\paragraph{Community detection in hypergraphs} The hypergraph network is a powerful tool to represent the complex entity relations with higher-order interactions~\citep{ke2019community}. A typical undirected hypergraph is denoted as $H = (V,E)$, where $V = [p]$ is the set of nodes and $E$ is the set of undirected hyperedges. Each hyperedge in $E$ is a subset of $V$, and we call the hyperedge an order-$K$ edge if the corresponding subset involves $K$ nodes. We call $H$ a $K$-uniform hypergraph if $E$ only contains order-$K$ edges. 

It is natural to represent the $K$-uniform hypergraph using a binary order-$K$ adjacency tensor. Let $\tY \in \{0,1\}^{p \times \cdots \times p}$ denote the adjacency tensor, where the entries encode the presence or absence of order-$K$ edges among $p$ nodes. Specifically, for all $(i_1,\ldots,i_K) \in [p]^K$, we have
\begin{equation}
    \tY(i_1,...,i_K) =  \begin{cases}
    1  & \text{if }  (i_1,...,i_K) \in E,\\
    0 & \text{if }  (i_1,...,i_K) \notin E.
    \end{cases}
\end{equation}

Assume that there exist $r$ disjoint communities among $p$ nodes, and the connection probabilities depend on the community assignments and node-specific parameters. 
Then, the equation~\eqref{eq:model_tensor} models $\bbE \tY$ with unknown degree heterogeneity $\mtheta$ and sub-Gaussianity parameter $\sigma^2 = 1/4$.

\paragraph{Multi-layer weighted network} 
Multi-layer weighted network data consists of multiple networks over the same set of nodes. One representative example is the brain connectome data~\citep{zhang2019tensor}. The multi-layer weighted network $\tY$ has dimension of $p \times p \times L$, where $p$ denotes the number of brain regions of interest, and $L$ denotes the number of layers (networks). Each of the $L$ networks describes one aspect of the brain connectivity, such as functional connectivity or structural connectivity. The resulting tensor $\tY$ consists of a mixture of slices with various data types.

Assume that there exist $r$ disjoint communities among $p$ nodes and $r_l$ disjoint communities among the $L$ layers. The multi-layer network community detection is modeled by the general asymmetric dTBM model~\eqref{eq:model_tensor}
\begin{equation}
    \bbE \tY = \tS \times_1 \mTheta \mM \times_2 \mTheta \mM \times_3 \mTheta_l \mM_l, 
\end{equation}
where $(\mtheta \in \bbR^p, \mM \in \{0,1\}^{p \times r})$ and $ (\mtheta_l \in \bbR^L, \mM_l \in \{0,1\}^{L \times r_l})$ are the degree heterogeneity and membership matrices corresponding to the community structure for $p$ nodes and $L$ layers, respectively. 

\paragraph{Gaussian higher-order clustering} Datasets in various fields such as medical image, genetics, and computer science are formulated as Gaussian tensors. One typical example is the multi-tissue gene expression dataset, which records different gene expressions in different individuals and different tissues. The dataset, denoted as $\tY \in \bbR^{p \times n \times t}$, consists of the expression data for $p$ genes of $n$ individuals in $t$ tissues. 

Assume that there exist $r_1, r_2, r_3$ disjoint clusters for $p$ genes, $n$ individuals, and $t$ tissues, respectively. We apply the general asymmetric dTBM model~\eqref{eq:model_tensor} 
\begin{equation}
    \bbE\tY = \tS \times_1 \mTheta_1 \mM_1 \times_2 \mTheta_2 \mM_2 \times_3 \mTheta_3 \mM_3, 
\end{equation}
where $\{(\mtheta_k, \mM_k)\}_{k=1}^3$ represents the degree heterogeneity and membership for genes, individuals, and tissues.

\begin{rmk}[Comparison with non-degree models]
Our dTBM uses fewer block parameters than TBM. In particular, every non-degree $r_1$-block tensor can be represented by a \emph{degree-corrected} $r_2$-block tensor with $r_2\leq r_1$. In particular, there exist tensors with $r_1=p$ but $r_2=1$, so the reduction in model complexity can be dramatic from $p$ to 1. This fact highlights the benefits of introducing degree heterogeneity in higher-order clustering tasks.
\end{rmk}

\subsection{Identifiability under angle gap condition}\label{subsec:identify}
The goal of clustering is to estimate the partition function $z$ from model~\eqref{eq:model_tensor}. For ease of notation, we focus on symmetric tensors; the extension to non-symmetric tensors are similar. We use $\tP$ to denote the following parameter space for $(z,\tS,\mtheta)$,

\vspace{-0.5cm}
\small
\begin{align}
\tP=  &\bigg\{  (z,\tS,\mtheta)\colon  \ \mtheta\in\mathbb{R}^p_{+},\ 
{c_1p\over r}\leq |z^{-1}(a)|\leq {c_2 p\over r},  c_3\leq \onormSize{}{\text{Mat}(\tS)_{a:}}\leq c_4, \onorm{\mtheta_{z^{-1}(a)}}_1=|z^{-1}(a)|, a\in[r]\bigg\}\label{eq:family}
\end{align}
\normalsize
where $c_i>0$'s are universal constants. We briefly describe the rationale of the constraints in~\eqref{eq:family}. 
First, the entrywise positivity constraint on  $\mtheta\in\mathbb{R}^p_{+}$ is imposed to avoid sign ambiguity between entries in $\mtheta_{z^{-1}(a)}$ and $\tS$. This constraint allows the trigonometric $\cos$ to describe the angle similarity in the Assumption~\ref{assmp:min_gap} below and Sub-algorithm~\hyperref[alg:main]{2} in Section~\ref{sec:alg}. Note that the positivity constraint can be achieved without sacrificing model flexibility, by using a slightly larger dimension of $\tS$ in the factorization~\eqref{eq:model_tensor}; see Example~\ref{ex:positive} below. Second, recall that the quantity $|z^{-1}(a)|$ denotes the number of nodes in the $a$-th community. The constants $c_1, c_2$ in the $|z^{-1}(a)|$ bounds assume the roughly balanced size across $r$ communities.  
{Third, the constant $c_3$ requires that all slides in $\tS$ have non-degenerate norm. Particularly, the lower bound $c_3$ excludes the purely zero slide to avoid trivial non-identifiability of model \eqref{eq:model_tensor}; see Example~\ref{example:c3} below. The upper bound $c_4$ is a technical constraint to avoid the slides with diverging norm as dimension grows.} 
Lastly, the $\ell_1$ normalization $\onormSize{}{\mtheta_{z^{-1}(a)}}_1=|z^{-1}(a)|$ is imposed to avoid the scalar ambiguity between $\mtheta_{z^{-1}(a)}$ and $\tS$. This constraint, again, incurs no restriction to model flexibility but makes our presentation cleaner.   {Our constraints in $\tP$ are mild compared with previous literature; see Table~\ref{tab:para} for comparison.}

\begin{example}[Positivity of degree parameters]~\label{ex:positive}
Here we provide an example to show the positivity constraint on $\mtheta$ incurs no loss on the model flexibility. 
Consider an order-3 dTBM with core tensor $\tS = 1$ and degree $\mtheta = (1,1,-1, -1)^T$. We have the mean tensor 
\begin{equation}
    \tX = \tS \times_1 \mTheta \mM \times_2 \mTheta \mM \times_3 \mTheta \mM,
\end{equation}
where $\mTheta = \text{diag}(\mtheta)$ and $\mM = (1,1,1,1)^T$. Note that $\tX \in \bbR^{4 \times 4 \times 4}$ is a 1-block tensor with \emph{mixed-signed} degree $\mtheta$, and the mode-3 slices of $\tX$ are
\begin{equation}
    \tX_{::1} = \tX_{::2} = - \tX_{::3} = -\tX_{::4}  = \begin{bmatrix}
    1 & 1& -1& -1 \\
    1 & 1& -1& -1 \\
    -1& -1&1 & 1\\
     -1& -1&1 & 1
    \end{bmatrix}.
\end{equation}
Now, instead of original decomposition, we encode $\tX$ as a 2-block tensor with \emph{positive-signed} degree. Specifically, we write
\begin{equation}
     \tX = \tS' \times_1 \mTheta' \mM' \times_2 \mTheta' \mM'  \times_3 \mTheta' \mM' , 
\end{equation}
where $\mTheta' = \text{diag}(\mtheta') = \text{diag}(1,1, 1, 1)$, the core tensor $\tS' \in \bbR^{2 \times 2 \times 2}$ has following mode-3 slices, and the membership matrix $\mM'\in\{0,1\}^{4\times 2}$ defines the clustering $z'\colon[4]\to[2]$; i.e.,
\begin{equation}
    \tS'_{::1} = -\tS'_{::2} = \begin{bmatrix}
    1 & -1 \\
    -1 & 1
    \end{bmatrix}, \quad \mM' = \begin{bmatrix}
   1 & 0 \\
   1 &0 \\
    0 & 1 \\
    0&1
    \end{bmatrix}.
\end{equation}
 The triplet $(z',\tS', \mtheta')$ lies in our parameter space~\eqref{eq:family}. In general, we can always reparameterize an $r$-block tensor with mixed-signed degree using a $2r$-block tensor with positive-signed degree. Since we assume $r=\tO(1)$ throughout the paper, the splitting does not affect the error rates of our interest.
 \end{example}

  \begin{example}[Non-identifiability with purely zero core slice]\label{example:c3} 
 Consider an order-2 dTBM with core tensor $\mS = \begin{pmatrix} 0 & 0\\
    1 & -1
    \end{pmatrix}$ degree matrices $ \mTheta_1 =\mTheta_2= \text{diag}(1,1,1,1)$, and mean tensor 
    \begin{equation}
        \tX =  \mTheta_1 \mM \mS   \mM^T \mTheta_2, \quad \text{with } \mM = \begin{bmatrix} 1 & 0 \\
        1 & 0 \\
        0 & 1 \\
        0 & 1 \\
        \end{bmatrix}.
    \end{equation}
Replacing $\mTheta_1$ by $\mTheta'_1 = (3/2, 1/2, 1,1)$ leads to the same mean tensor $\tX$. 
\end{example}

 \begin{table}[th]
    \centering
    \resizebox{.95\textwidth}{!}{
    \begin{tabular}{c|c ccc}
    \hline
    Assumptions in parameter space&\cite{gao2018community}& \cite{han2022exact}& \cite{ke2019community} & {\bf Ours}\\
    \hline
         Balanced community sizes &$\surd$ &$\surd$&$\surd$ & $\surd$  \\
          Bounded core tensors & $\surd$ & $\times$ & $\surd$ &$\surd$\\
         Balanced degrees& $\surd$&-&$\surd$ &$\surd$\\
         Flexible in-group connections &$\times$ &$\surd$&$\surd$& $\surd$\\
         Gaps among cluster centers & In-between cluster difference & Euclidean gap & Eigen gap & Angle gap\\
         \hline
    \end{tabular}
    }
    \caption{Parameter space comparison between previous work with our assumption.}
    \label{tab:para}
\end{table}

We now provide the identifiability conditions for our model before estimation procedures. When $r=1$, the decomposition~\eqref{eq:model_tensor} is always unique (up to cluster label permutation) in $\tP$, because dTBM is equivalent to the rank-1 tensor family under this case. When $r\geq 2$, the Tucker rank of signal tensor $\bbE\tY$ in~\eqref{eq:model_tensor} is bounded by, but not necessarily equal to, the number of blocks $r$~\citep{wang2019multiway}. Therefore, one can not apply the classical identifiability conditions for low-rank tensors to dTBM. Here, we introduce a key separation condition on the core tensor. 

\begin{assumption}[Angle gap] \label{assmp:min_gap}Let $\mS = \text{Mat}(\tS)$. Assume that the minimal gap between normalized rows of $\mS$ is bounded away from zero; i.e.,
\begin{equation}\label{eq:minimal_gap}
    \Delta_{\min} := \min_{a \neq b \in [r]} \onorm{\frac{\mS_{a:}}{\onormSize{}{ \mS_{a:}}} - \frac{\mS_{b:}}{\onormSize{}{ \mS_{b:}}} }>0, \quad \text{for} \quad r \geq 2.
\end{equation}
\end{assumption}
We make the convention $\Delta_{\min} = 1$ for $r = 1$. Equivalently, \eqref{eq:minimal_gap} says that none of the two rows in $\mS$ are parallel; i.e., $\max_{a \neq b\in [r]}\cos \of{\mS_{a:},\  \mS_{b:}}  = 1-\Delta^2_{\min}/2<1$. The quantity $\Delta_{\min}$ characterizes the non-redundancy among clusters measured by angle separation. The denominators involved in definition~\eqref{eq:minimal_gap} are well posed because of the lower bound on $\onorm{\mS_{a:}}$ in~\eqref{eq:family}. 

Our first main result is the following theorem showing the sufficiency and necessity of the angle gap separation condition for the parameter identifiability under dTBM. 

\begin{thm}[Model identifiability]\label{thm:unique} Consider the dTBM with $r\geq 2$ {and $K \geq 2$}. The parameterization~\eqref{eq:model_tensor} is unique in $\tP$ up to cluster label permutations, if and only if Assumption~\ref{assmp:min_gap} holds.
\end{thm}

The identifiability guarantee for the dTBM is stronger than classical Tucker model. In the Tucker model, the factor matrix $\mM$ is identifiable only up to orthogonal rotations. In contrast, our model does not suffer from rotational invariance. As we will show in Section~\ref{sec:alg}, each column of the membership matrix $\mM$ can be precisely recovered under our algorithm. This property benefits the interpretation of dTBM in practice.

\section{Statistical-computational critical values for higher-order tensors}\label{sec:limits}

\subsection{Assumptions} \label{sec:prelim}

We propose the signal-to-noise ratio~(SNR),
\begin{align}\label{eq:gamma}
  \text{SNR}:= \Delta^2_{\min}/\sigma^2 = p^{\gamma}, 
\end{align}
with varying $\gamma \in \bbR$ that quantifies different regimes of interest. We call $\gamma$ the \emph{signal exponent}. Intuitively, a larger SNR, or equivalently a larger $\gamma$, benefits the clustering in the presence of noise.  With quantification~\eqref{eq:gamma}, we consider the following parameter space,
\begin{equation}\label{eq:gammafamily}
    \tP(\gamma) = \tP\cap\{\tS \text{ satisfies SNR condition~\eqref{eq:gamma} with $\gamma$} \}.
\end{equation}
The $1$-block dTBM does not belong to the space $\tP(\gamma)$ when $\gamma < 0$, due to the convention in Assumption~\ref{assmp:min_gap}. Our goal is to characterize the clustering accuracy with respect to $\gamma$ under the space $\tP(\gamma)$.

In our algorithmic development, we often refer to the regime of balanced degree heterogeneity. We call the degree $\mtheta$ \emph{balanced} if
\begin{equation}\label{eq:degree}
{\min_{a\in[r]} \onormSize{}{\mtheta_{z^{-1}(a)}}=\left(1+o(1)\right)\max_{a\in[r]}\onormSize{}{\mtheta_{z^{-1}(a)}}}.
\end{equation}

The following lemma provides the rationale of balanced degree assumption. We show the close relation between angle gaps in the mean tensor $\tX$ and the core tensor $\tS$ under balanced degree heterogeneity.

\begin{lem}[Angle gaps in $\tX$ and $\tS$]\label{lem:angle_gap_x} Consider the dTBM model~\eqref{eq:model_tensor} under the parameter space $\tP$ in \eqref{eq:family} {with $r \geq 2$. Suppose $\mtheta$ is balanced satisfying~\eqref{eq:degree} {and $\min_{i\in[p]}\theta(i)\geq c$ from some constant $c>0$}.} Then, {as $p \rightarrow \infty$}, for all $i,j$ such that $z(i) \neq z(j)$, we have
\begin{equation}
 \cos(\mX_{i:}, \mX_{j:})\asymp  \cos(\mS_{z(i):}, \mS_{z(j):}),
\end{equation}

where $\mX =\mat(\tX)$ and $\mS = \mat (\tS)$.
\end{lem}
In practice, an estimation algorithm has access to a noisy version of $\tX$ but not $\tS$. 
Our goal is to establish the algorithm performance with respect to the signal $\Delta^2_{\min}$ in the core tensor. By Lemma~\ref{lem:angle_gap_x}, the mapping from the core tensor $\mS_{z(i):}$ to the mean tensor $\mX_{z(i):}$ preserves the angle information $\Delta_{\min}^2$ under balanced degree heterogeneity~\eqref{eq:degree}. Therefore, the balanced degree assumption helps to exclude the cases in which the degree heterogeneity distorts the algorithm guarantees. 

Here, we provide an example to illustrate the insufficiency of $\Delta_{\min}^2$ in the absence of balanced degrees.

\begin{example}[Insufficiency of $\Delta_{\min}^2$ in the absence of balanced degrees] Consider an order-2 $(p,p)$-dimensional dTBM with core matrix
\begin{equation}\label{eq:modelexample}
    \mS = \begin{pmatrix} 1 &a\\
    1 & -a
    \end{pmatrix}, 
\end{equation}
and $ \mtheta \text{ such that } \onormSize{}{\mtheta_{z^{-1}(1)}}^2 = p^m \onormSize{}{\mtheta_{z^{-1}(2)}}^2$,
where {$m \in [-1,1]$} is a scalar parameter controlling the skewness of degrees. Let $\Delta_{\mX}^2$ denote the minimal angle gap of the mean tensor, defined by
\begin{equation}\label{eq:delta_x}
    \Delta_{\mX}^2 \coloneqq \min_{i,j \in [p], z(i) \neq z(j)} \onorm{ \frac{\mX_{i:}}{\onormSize{}{\mX_{i:}}}  -  \frac{\mX_{j:}}{\onormSize{}{\mX_{j:}}}  }, \quad 
\end{equation}
where $ \mX = \mat(\tX)$.
Take $ a = p^{-1/4}$ in the model setup~\eqref{eq:modelexample}. We have 
\begin{align}
    \Delta_{\min}^2 = \frac{2 a^2}{1 + a^2} \asymp p^{-1/2}, \quad
    \Delta_{\mX}^2  = \frac{2\onormSize{}{\mtheta_{z^{-1}(2)}}^{2} a^2}{\onormSize{}{\mtheta_{z^{-1}(1)}}^{2} +  \onormSize{}{\mtheta_{z^{-1}(2)}}^{2}a^2} \asymp p^{-1/2-m}.
\end{align}

{Based on the Theorem~\ref{thm:stats} in Section~\ref{sec:limits}, the dTBM is impossible to solve when $\Delta^2_{\mX} \lesssim p^{-1}$ even though $\Delta_{\min}^2 \asymp p^{-1/2}$
}; that is, the dTBM estimation depends on the relative magnitude of $m$ vs. $1/2$. In such a setting, the proposed signal notion $\Delta^2_{\min}$ alone fails to fully characterize dTBM. 
\end{example}

\begin{rmk}[Flexibility in balanced degree assumption] One important note is that our balance assumption~\eqref{eq:degree} does not preclude the mild degree heterogeneity. In fact, within each of the clusters, we allow the highest degree at the order $\tO(p)$, whereas the lowest degree at the order $\Omega(1)$. This range is more relaxed than previous work \citep{gao2018community} that restricts the highest degree in the sub-linear regime $o(p)$ and the lowest degree at the order $\Omega(1)$. 
\end{rmk}

\begin{rmk}[Similar assumptions in literature]
Similar degree regulations are not rare in literature. In higher-order tensor model \citep{ke2019community}, the degree assumption $\max_{a\in[r]} \onormSize{}{\mtheta_{z^{-1}(a)}} \leq $ $ C \min_{a\in[r]}\onormSize{}{\mtheta_{z^{-1}(a)}}$ is made to ensure degree balance across communities. 
In \cite{gao2018community}, the degree distribution is restricted to ${1\over |z^{-1}(a)|}\sum_{i\in z^{-1}(a)}\theta_i=1+o(1)$ for all communities. 
\end{rmk}

Last, let $\hat z$ and $z$ be the estimated and true clustering functions in the family~\eqref{eq:family}. Define the misclustering error by
\[
\ell(\hat z, z)={1\over p}\min_{\pi \in \Pi}\sum_{i\in[p]}\ind\{\hat z(i)\neq \pi \circ z(i)\},
\]
where $\pi: [r] \mapsto [r]$ is a permutation of cluster labels, $\circ$ denotes the composition operation, and $\Pi$ denotes the collection of all possible permutations. The infimum over all permutations accounts for the ambiguity in cluster label permutation. 

In Sections~\ref{sec:statlimit} and \ref{sec:complimit}, we provide the phase transition of $\ell (\hat z, z)$ for general Gaussian dTBMs~\eqref{eq:model_tensor} without symmetric assumptions. For general (asymmetric) Gaussian dTBMs, we assume Gaussian noise $\tE(i_1,\ldots,i_K)\stackrel{\text{i.i.d.}}{\sim} N(0,\sigma^2)$, and we extend the parameter space~\eqref{eq:family} to allow $K$ clustering functions $\{z_k\}_{k\in[K]}$, one for each mode. For notational simplicity, we still use $z$ and $\tP(\gamma)$ for this general (asymmetric) model. All results should be interpreted as the worst-case results across $K$ modes. 

\subsection{Statistical critical value}\label{sec:statlimit}
The statistical critical value means the SNR required for solving dTBMs with \emph{unlimited computational cost.}   Our following result shows the minimax lower bound for exact recovery and the matching upper bound for maximum likelihood estimator (MLE).  We consider the Gaussian MLE, denoted as $(\hat z_{\text{MLE}},\hat \tS_{\text{MLE}},  \hat \mtheta_{\text{MLE}})$, over the estimation space $\tP$, where 
\begin{equation}\label{eq:mle}
    (\hat z_{\text{MLE}}, \hat \tS_{\text{MLE}}, \hat \mtheta_{\text{MLE}}) = \argmin_{ (z,\tS,\mtheta) \in \tP} \onormSize{}{\tY - \tX(z,\tS,\mtheta)}_F^2.
\end{equation}

\begin{thm}[Statistical critical value]\label{thm:stats} Consider general Gaussian dTBMs with parameter space $\tP(\gamma)$ and $K\geq 2$. Then, we have the following statistical phase transition. 

\begin{itemize}[wide]
    \item 
    \textbf{Impossibility.}  {Assume $p \rightarrow \infty$ and $2 \leq r\lesssim p^{1/3}$.} Let $\tP_{\tS}(\gamma) \coloneqq \{ \tS: c_3\leq \onormSize{}{\text{Mat}(\tS)_{a:}}\leq c_4, a\in[r]\} \cap \{ \tS: \Delta_{\min}^2 = p^{\gamma}\}$ denote the space for valid $\tS$ satisfying SNR condition~\eqref{eq:gamma}, and $\tP_{z ,\mtheta} \coloneqq \{ \mtheta\in\mathbb{R}^p_{+},\ 
{c_1p\over r}\leq |z^{-1}(a)|\leq {c_2 p\over r}, \onorm{\mtheta_{z^{-1}(a)}}_1=|z^{-1}(a)|, a \in [r] \}$ denote the space for valid $(z, \mtheta)$,  where $c_1, c_2, c_3, c_4$ are the constants in parameter space \eqref{eq:family}. If the signal exponent satisfies $\gamma < -(K-1)$, then, for any true core tensor $\tS \in \tP_{\tS}(\gamma)$, no estimator $\hat z_{\text{stat}}$ achieves exact recovery in expectation; that is, when $\gamma< -(K-1)$, we have 
    \begin{equation}\label{eq:minminmax}
   \liminf_{p \rightarrow \infty}\inf_{\tS \in  \tP_{\tS}(\gamma)}  \inf_{\hat z_{\text{stat}}}\sup_{ (z, \mtheta) \in \tP_{z, \mtheta}} \bbE \left[ p\ell(\hat z_{\text{stat}}, z) \right]\geq 1.
\end{equation}
{ 
    Further, we define the parameter space $\tP'(\gamma') \coloneqq \tP\cap \{ \Delta_{\mX}^2 = p^{\gamma'} \}$, where $\Delta_{\mX}^2$ is the mean tensor minimal gap in \eqref{eq:delta_x}. When $\gamma' < -(K-1)$, we have 
    \begin{equation}
        \liminf_{p \rightarrow \infty} \inf_{\hat z_{\rm stat}} \sup_{(z, \tS, \mtheta) \in \tP'(\gamma')}  \bbE \left[ p\ell(\hat z_{\text{stat}}, z) \right]\geq 1.
    \end{equation}
    }

\item \textbf{MLE achievability.} 
Suppose that the signal exponent satisfies $\gamma >-(K-1)+c_0$ for an arbitrary constant {$c_0>0$}. Furthermore, assume that {$\mtheta$} is balanced and $\min_{i\in[p]}\theta(i)\geq c$ from some constant $c>0$. Then, {when $p \rightarrow \infty$, for fixed $r \geq 1$, }the MLE in~\eqref{eq:mle} achieves exact recovery in high probability; that is,
\begin{equation}
\ell(\hat z_{\text{MLE}}, z) \lesssim \text{SNR}^{-1}\exp\left(-{p^{K-1}\text{SNR}\over r^{K-1}}\right) \to 0,
\end{equation}
with probability going to 1. 
\end{itemize}
\end{thm}

The proofs for the two parts in Theorem~\ref{thm:stats} are in the Appendix B, Section~\ref{sec:statprove1} and Section~\ref{sec:statprove2}, respectively. The first part of Theorem~\ref{thm:stats} demonstrates impossibility of exact recovery whenever the core tensor $\tS$ satisfies SNR condition~\eqref{eq:gamma} with exponent $\gamma < -(K-1)$. The proof is information-theoretical, and therefore the results apply to all statistical estimators, including but not limited to MLE and trace maximization~\citep{ghoshdastidar2017uniform}. The minimax bound~\eqref{eq:minminmax} indicates the worst case impossibility for a particular core tensor $\tS$ with signal exponent $\gamma < -(K-1)$; i.e., under the assumptions of Theorem~\ref{thm:stats}, when $\gamma < -(K-1)$, we have 
\begin{equation}
    \liminf_{p \rightarrow \infty}  \inf_{\hat z_{\text{stat}}}\sup_{ (z,\tS, \mtheta) \in \tP(\gamma)} \bbE \left[ p\ell(\hat z_{\text{stat}}, z) \right]\geq 1.
\end{equation}
Such worst case impossibility is studied in related works \citep{han2022exact, gao2018community} while our lower bound~\eqref{eq:minminmax} provides a stronger impossibility statement for arbitrary core tensors with weak signals. The second part of Theorem~\ref{thm:stats} shows the exact recovery of MLE when $\gamma>-(K-1)+c_0$ for an arbitrary constant $c_0>0$. Combining the impossibility and achievability results, we conclude that the boundary $\gamma_{\text{stat}} \coloneqq -(K-1)$ is the critical value for statistical performance of dTBM with respect to our SNR.

\subsection{Computational critical value}\label{sec:complimit}
The computational critical value means the minimal SNR required for exact recovery with \emph{polynomial-time} computational cost. An important ingredient to establish the computational limits is the \emph{hypergraphic planted clique (HPC) conjecture} \citep{zhang2018tensor, brennan2020reducibility}. The HPC conjecture indicates the impossibility of fully recovering the planted cliques with polynomial-time algorithm when the clique size is less than the number of vertices in the hypergraph. The formal statement of HPC detection conjecture is provided in Definition~\ref{def:HPC} and Conjecture~\ref{hypo:HPC} as follows.  

\begin{defn}[Hypergraphic planted clique (HPC) detection]~\label{def:HPC}Consider an order-$K$ hypergraph $H = (V,E)$ where $V = [p]$ collects vertices and $E$ collects all the order-$K$ edges. Let $\tH_k(p, 1/2)$ denote the Erd\H{o}s-R\'{e}nyi $K$-hypergraph where the edge $(i_1,\ldots, i_K)$ belongs to $E$ with probability $1/2$. Further, we let $\tH_K(p, 1/2, \kappa)$ denote the hyhpergraph with planted cliques of size $\kappa$. Specifically, we generate a hypergraph from $\tH_k(p, 1/2)$, pick $\kappa$ vertices uniformly from $[p]$, denoted $K$, and then connect all the hyperedges with vertices in $K$. Note that the clique size $\kappa$ can be a function of $p$, denoted $\kappa_p$. The order-$K$ HPC detection aims to identify whether there exists a planted clique hidden in an Erd\H{o}s-R\'{e}nyi $K$-hypergraph. The HPC detection is formulated as the following hypothesis testing problem
\begin{equation}
    H_0:\ H \sim \tH_K(p,1/2) \quad \text{versus} \quad H_1: \ H \sim \tH_K(p,1/2, \kappa_p).
\end{equation}
\end{defn}

\begin{conjecture}[HPC conjecture]\label{hypo:HPC} Consider the HPC detection problem in Definition~\ref{def:HPC} with $K \geq 2$. Suppose the sequence $\{\kappa_p\}$ such that $\limsup_{p \rightarrow \infty} \log \kappa_p/ \log \sqrt{p} \leq (1 - \tau)$ for any $\tau > 0$. Then, for every sequence of polynomial-time test $\{ \varphi_p\}: H \mapsto \{0,1\}$ we have 
\begin{equation}
    \liminf_{p \rightarrow \infty} \bbP_{H_0} \of{ \varphi_p(H) =1 } +  \bbP_{H_1} \of{ \varphi_p(H) =0} > \frac{1}{2}.
\end{equation}
\end{conjecture}

Under the HPC conjecture, we establish the SNR lower bound that is necessary for any \emph{polynomial-time} estimator to achieve exact clustering.

\begin{thm}[Computational critical value]\label{thm:comp} 
Consider general Gaussian dTBMs under the parameter space $\tP$ with $K\geq 2$. Then, we have the following computational phase transition.

\begin{itemize}[wide]
    \item \textbf{Impossibility.} Assume HPC conjecture holds {and $r \geq 2$}. If the signal exponent satisfies $\gamma < -K/2$, then, no \emph{polynomial-time estimator} $\hat z_{\text{comp}}$ achieves exact recovery in expectation as $p\to \infty$; that is, when $\gamma<-K/2$, we have 
\begin{align}
   \liminf_{p\to \infty}\sup_{(z, \tS, \mtheta) \in \tP(\gamma)}  \bbE \left[ p\ell(\hat z_{\text{comp}}, z) \right]\geq 1.
\end{align}
\item \textbf{Polynomial-time algorithm achievability.} Suppose that we have {fixed $r \geq 1$}, and the signal exponent satisfies $\gamma  >  -K/2+c_0$ for an arbitrary constant {$c_0>0$}. Furthermore, assume that the degree $\mtheta$ is balanced, lower bounded in that $\min_{i\in[p]}\theta_i\geq c$ for some constant $c>0$, and satisfies the locally linear stability in Definition~\ref{def:stable} {in the neighborhood $\tN(z,\varepsilon)$ for all $\varepsilon \leq E_0$ and some $E_0 \gtrsim \log^{-1}p $.}
Then, {as $p \rightarrow \infty$}, there exists a polynomial-time algorithm $\hat z_{\text{ploy}}$ that achieves exact recovery in high probability; that is,
\begin{equation}
 \ell(\hat z_{\text{poly}}, z) \lesssim \text{SNR}^{-1}\exp \of{ - \frac{p^{K-1}\text{SNR}}{r^{K-1}} } \to 0,
\end{equation}
with probability going to 1. 
\end{itemize}
\end{thm}

The proofs for the two parts in Theorem~\ref{thm:comp} are in the Appendix B, Section~\ref{sec:compprove1} and Section~\ref{sec:statprove2}, respectively. The first part of Theorem~\ref{thm:comp} indicates the impossibility of exact recovery by polynomial-time algorithms when $\gamma < -K/2$, and the second part shows the existence of such algorithm when {$\gamma > -K/2+c_0$  for an arbitrary constant $c_0 > 0$} under extra technical assumptions. In Section~\ref{sec:alg}, we will present an efficient polynomial-time algorithm in this setting. Therefore, we conclude that $\gamma_{\text{comp}}:=-K/2$ is the critical value for computational performance of dTBM with respect to our SNR. 

\begin{rmk}[Statistical-computational gaps]
Now, we have established the phase transition of exact clustering under order-$K$ dTBM by combining Theorems~\ref{thm:stats} and \ref{thm:comp}. Figure~\ref{fig:phase_axis} summarizes our results of critical SNRs when $K \geq 2$. In the weak SNR region $\gamma < -(K-1)$, no statistical estimator succeeds in degree-corrected higher-order clustering. In the strong SNR region $\gamma  > -K/2$, our proposed algorithm precisely recovers the clustering in polynomial time. In the moderate SNR regime, $-(K-1)\leq \gamma \leq -K/2$, the degree-corrected clustering problem is statistically easy but computationally hard. Particularly, dTBM reduces to matrix degree-corrected model when $K =2$, and the statistical and computational bounds show the same critical value. When $K =1$, dTBM reduces to the degree-corrected sub-Gaussian mixture model (GMM) with model
\begin{equation}
    \mY = \mTheta \mM \mS + \mE,
\end{equation}
where $\mY \in \bbR^{p \times d}$ collects $n$ data points in $\bbR^d$, $\mS \in \bbR^{r \times d}$ collects the $d$-dimensional centroids for $r$ clusters, and $\mTheta \in \bbR^{p \times p}, \mM \in \{0,1\}^{p \times r}, \mE \in \bbR^{p \times d}$ have the same meaning as in dTBM. \cite{lu2016statistical} implies that polynomial-time algorithms are able to achieve the statistical minimax lower bound in GMM. Therefore, we conclude that the statistical-computational gap emerges only for higher-order tensors with $K \geq 3$. The result reveals the intrinsic distinctions among (vector) one-dimensional clustering, (matrix) biclustering, and (tensor) higher-order clustering. 
\end{rmk}

\section{Polynomial-time algorithm under mild SNR}\label{sec:alg}
\begin{figure}[ht!]
\centering
\includegraphics[width=.9\columnwidth]{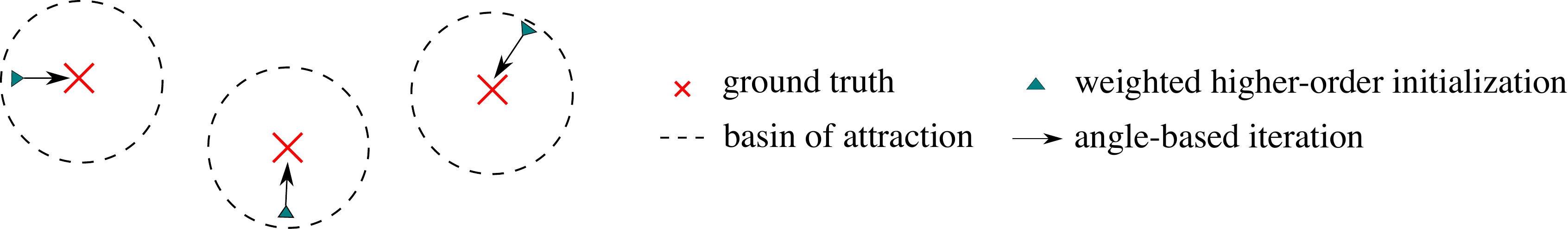}
\caption{Illustration of our global-to-local algorithm.}\label{fig:demo}
\end{figure}

In this section, we present an efficient polynomial-time clustering algorithm under mild SNR. The procedure takes a global-to-local approach. See Figure~\ref{fig:demo} for illustration. The global step finds the basin of attraction with polynomial misclustering error, whereas the local iterations improve the initial clustering to exact recovery. Both steps are critical to obtain a satisfactory algorithm output. In what follows, we first use the symmetric tensor as a working example to describe the algorithm procedures to gain insight. {Our theoretical analysis focuses on dTBMs with symmetric mean tensor and independent sub-Gaussian noises} 
such as Gaussian and uniform observations.   {The extensions for Bernoulli observations and other practical issues are in Sections~\ref{subsec:ber} and \ref{subsec:exten}.}

To construct algorithm guarantees,  we introduce the misclustering loss between an estimator $\hat z$ and the true $z$:
\begin{align}
     L(\hat z, z) = \frac{1}{p}  \sum_{i \in [p]} \theta(i) \sum_{b \in [r]}  \ind &\offf{ \hat z(i) = b } \cdot \onorm{ \off{ \mS_{ z(i):}  }^s - \off{ \mS_{b:}  }^s  }^2, \label{eq:defnofL}
\end{align}
where the superscript $\cdot^{s}$ denotes the normalized vector; i.e., $\ma^s:=\ma/\onorm{\ma}$ if $\ma\neq 0$ and $\ma^{s}=0$ if $\ma=0$ for any vector $\ma$. The following lemma indicates the close relationship between the loss $L(\hat z, z)$ and error $\ell(\hat z, z)$. The loss $L(\hat z, z)$ serves as an important intermediate quantity to control the misclustering error.

\begin{lem}[Relationship between misclustering error and loss]\label{lem:mis} {Consider the dTBM under the parameter space $\tP$.} Suppose $\min_{i \in [p]}\theta(i) > c$ for some constant $c > 0$. We have $\ell(\hat z, z) \Delta_{\min}^2 \leq L(\hat z, z)$.
\end{lem}

\subsection{Weighted higher-order initialization}

We start with weighted higher-order clustering algorithm as initialization.   {We take an order-3 tensor and the clustering on the first mode as illustration for insight.} Consider noiseless case with $\tX = \bbE\tY$ and $\mX = \text{Mat}(\tX)$. 
By model~\eqref{eq:model_tensor}, for all $i \in [p]$, we have
\begin{equation}\label{eq:kmeans}
    \theta(i)^{-1} \mX_{i:} = \off{\text{Mat}( \tS \times_2 \mTheta \mM \times_3  \mTheta \mM )}_{z(i):}. 
\end{equation}
This implies that, all node $i$ belonging to the $a$-th community (i.e., $z(i)=a$) share the same normalized mean vector $\theta(i)^{-1} \mX_{i:}$, and vice versa. Intuitively, one can apply $k$-means clustering to the vectors $\{ \theta(i)^{-1} \mX_{i:} \}_{i\in[p]}$, which leads to main idea of our Sub-algorithm~\hyperref[alg:main]{1}.

Specifically, our initialization consists of the denoising step and the clustering step. The denoising step (lines 1-2 in Sub-algorithm~\hyperref[alg:main]{1}) estimates $\tX$ from $\tY$ by a double projection spectral method.  
The first projection performs HOSVD~\citep{de2000multilinear} via $\mU_{\text{pre},k} = \text{SVD}_{r} \of{ \text{Mat}_k(\tY) }, k \in [3]$, where $\text{SVD}_r(\cdot)$ returns the top-$r$ left singular vectors. The second projection performs HOSVD on the projected $\tY$ onto the multilinear Kronecker space $\mU_{\text{pre},k} \otimes  \mU_{\text{pre},k}$; i.e.,
\begin{equation}\label{eq:two-step_factor}
    \hat \mU_1 = \text{SVD}_{r} \of{\text{Mat}_1\of{ \tY\times_2  \mU_{\text{pre},2} \mU_{\text{pre},2}^T \times_3   \mU_{\text{pre},3} \mU_{\text{pre},3}^T }}.
\end{equation}
and similar for $ \hat \mU_2 ,  \hat \mU_3 $.
The final denoised tensor $\hat \tX$ is defined by
\begin{equation}\label{eq:two-step_est}
    \hat \tX = \tY \times_1 \hat \mU_1 \hat 
\mU^T_1 \times_2 \hat \mU_2 \hat \mU^T_2 \times_3 \hat \mU_3 \hat \mU^T_3. 
\end{equation}
The double projection improves usual matrix spectral methods in order to alleviate the noise effects for $K\geq 3$~\citep{han2022exact}. The clustering step (lines 3-5 in Sub-algorithm~\hyperref[alg:main]{1}) performs the weighted $k$-means clustering. 
We write $\hat \mX=\text{Mat}_1(\hat \tX)$, and normalize the rows into $\hat \mX^s_{i:}=\onormSize{}{\hat \mX_{i:}}^{-1}\hat \mX_{i:}$ as a surrogate of $\theta(i)^{-1} \mX_{i:}$. Then, a weighted $k$-means clustering is performed on the normalized rows with weights equal to $\onormSize{}{\hat \mX_{i:}}^2$. The choice of weights is to bound the $k$-means objective function by the Frobenius-norm accuracy of $\hat \tX$. Unlike existing clustering algorithm~\citep{ke2019community}, we apply the clustering on the unfolded tensor $\hat \mX$ rather than on the factors $\hat \mU_k$. This strategy relaxes the singular-value gap condition~\citep{gao2018community, han2022exact}.
We assign degenerate rows with purely zero entries to an arbitrarily random cluster; these nodes are negligible in high-dimensions because of the lower bound on $\onormSize{}{\text{Mat}(\tS)_{a:}}$ in~\eqref{eq:family}. The final result gives the initial cluster assignment $z^{(0)}$. Full procedures for clustering are provided in Sub-algorithm~\hyperref[alg:main]{1}.

\begin{algorithm*}[th!]
\caption*{\bf Algorithm: Multiway spherical clustering for degree-corrected tensor block model }
\vspace{.15cm}
\begin{algorithmic}[1]
\Algphase{Sub-algorithm 1: Weighted higher-order initialization}
\INPUT Observation $\tY \in \bbR^{p\times \cdots \times p}$, cluster number $r$, relaxation factor $\eta > 1$ in $k$-means clustering.

\State {
  Compute factor matrices $ \mU_{\text{pre},k} = \text{SVD}_{r} (\text{Mat}_k(\tY)), k \in [K]$ and the $(K-1)$-mode projections 
\begin{equation}
    \tX_{\text{pre},k} = \tY \times_1   \mU_{\text{pre},1} \mU_{\text{pre},1}^T \times_2 \cdots  \times_{k-1} \mU_{\text{pre},k-1} \mU_{\text{pre},k-1}^T \times_{k+1} \mU_{\text{pre},k+1} \mU_{\text{pre},k+1}^T\times_{k+2} \cdots  \times_K \mU_{\text{pre},K} \mU_{\text{pre}, K}^T.
\end{equation}
}
\State {
Compute factor matrices $\hat \mU_k = \text{SVD}_{r}(\text{Mat}_k(\tX_{\text{pre},k})), k \in [K]$ and the denoised tensor
\begin{equation}
    \hat \tX = \tY \times_1 \hat \mU_1 \hat \mU^T_1 \times_2 \cdots \times_K \hat \mU_K \hat \mU^T_K.
\end{equation}
} 
\For{  {$k \in [K]$}}
\State {    Let $\hat \mX = \text{Mat}_k(\hat \tX)$ and $S_0=\{i \in [p]: \onormSize{}{\hat \mX_{i:}} = 0\}$. Set $\hat z(i)$ randomly in $[r]$ for $i \in S_0$.}
\State{  For all $i\in S_0^c$, compute normalized rows
$\hat \mX_{i:}^s :=\onormSize{}{\hat \mX_{i:}}^{-1} \hat \mX_{i:}.$
}
\State { Solve the clustering $\hat z_k \colon [p]\to[r]$ and centroids $\{\hat \mx_j\}_{j\in[r]}$ using weighted $k$-means, such that}
\begin{align}
    &\sum_{i \in S_0^c }  \onormSize{}{\hat \mX_{i:}}^2 \onormSize{}{\hat \mX_{i:}^s - \hat \mx_{\hat z_k(i)} }^2 
    \leq 
    \eta \min_{\substack{\bar \mx_j, j\in[r], \bar z_k(i),i\in S_0^{c}}} \sum_{i \in S^c } \onormSize{}{\hat \mX_{i:}}^2 \onormSize{}{ \hat \mX_{i:}^s -   \bar \mx_{\bar z_k(i)}}^2.
\end{align}
\EndFor

\OUTPUT {   Initial clustering $z^{(0)}_k \leftarrow \hat z_k, k \in [K]$.}

\Algphase{Sub-algorithm 2: Angle-based iteration}\label{alg:2}
\INPUT Observation $\tY \in \bbR^{p \times \cdots \times p}$, initialization $z^{(0)}_k \colon [p]\to[r], k \in [K]$ from Sub-algorithm 1, iteration number $T$.
\For {$t = 0$ to $T-1$}
\State Update the block tensor $\tS^{(t)}$ via
$\tS^{(t)} (a_1,...,a_K)= \text{Ave} \{\tY(i_1,\ldots,i_K): z^{(t)}_k(i_k) = a_k, k \in [K]\}.$
\For{ $k \in [K]$}
\State {   Calculate the reduced tensor $\tY^{\text{d}}_k \in \bbR^{r \times \cdots \times r \times p \times r \times \cdots \times r}$ via
\begin{equation}
    \tY^{\text{d}}_k(a_1,\ldots,a_{k-1}, i ,a_{k+1},\ldots,a_K) 
    = \text{Ave}\{\tY(i_1,\ldots,i_{k-1},i,i_{k+1},\ldots,i_K): z^{(t)}(i_j) = a_j, j \neq k \}
\end{equation}}

\State {   Let $\mY_k^{\text{d}} = \text{Mat}_k(\tY^{\text{d}})$ and $J_0 = \{ i\in[p]: \onorm{\mY^{\text{d}}_{i:}} = 0\}$. Set $z_k^{(t+1)}(i)$ randomly in $[r]$ for $i \in J_0$.}

\State {   Let $\mS^{(t)}_k = \text{Mat}_k(\tS^{(t)})$. For all $i \in J_0^c$, update the cluster assignment by
\begin{equation}
    z(i)^{(t+1)}_k = \argmax_{a \in [r]} \cos \left( \mY^{\text{d}}_{k,i:},\ \mS^{(t)}_{k, a:} \right).
\end{equation}}
\EndFor
\EndFor

\OUTPUT {    Estimated clustering $z^{(T)}_k: [p] \mapsto [r], k \in [K]$.}

\end{algorithmic}
\end{algorithm*}

We now establish the misclustering error rate of initialization.

\begin{thm}[Error for weighted higher-order initialization]\label{thm:initial} Consider the general sub-Gaussian dTBM with {fixed $r \geq 1$, $K \geq 2$}, i.i.d.\ noise under the parameter space $\tP$, and Assumption~\ref{assmp:min_gap}. Assume $\min_{i\in[p]}\theta(i) \geq c$ for some constant $c>0$. Let $\Delta_{\mX}$ denote the minimal gap in mean tensor defined in~\eqref{eq:delta_x}, and let $ z^{(0)}_k$ denote the output of Sub-algorithm~\hyperref[alg:main]{1}. With probability going to 1, {as $p \rightarrow \infty$}, we have
\begin{equation}
    \ell(z^{(0)}_k, z) \lesssim {\sigma^2 r^K p^{-K/2}\over \Delta_{\mX}^2}.
\end{equation}
Further, assume that $\mtheta$ is balanced as~\eqref{eq:degree}. We have
\begin{equation}\label{eq:ini}
 \ell(z_k^{(0)}, z) \lesssim {r^K p^{-K/2}\over \text{SNR}} \quad \text{and} \quad L(z^{(0)}_k, z) \lesssim  {\sigma^2 r^K p^{-K/2}},
\end{equation}
{with probability going to 1 as $p \rightarrow \infty$.}
\end{thm}

\begin{rmk}[Comparison to previous results] For fixed SNR, our initialization error rate with $K=2$ agrees with the initialization error rate $\tO(p^{-1})$ in matrix models~\citep{gao2018community}. Furthermore, in the special case of non-degree TBMs with $\mtheta = \mathbf{1}_p$, we achieve the same initial misclustering error $\tO(p^{-K/2})$ as in non-degree models~\citep{han2022exact}. Theorem~\ref{thm:initial} implies the advantage of our algorithm in achieving both accuracy and model flexibility. 
\end{rmk}

\begin{rmk}[Failure of conventional tensor HOSVD] If we use conventional HOSVD for tensor denoising; that is, we use $\mU_{\text{pre},k}$ in place of $\hat \mU_k$ in line 2, then the misclustering rate becomes $\tO(p^{-1})$ for all $K\geq 2$. This rate is substantially worse than our current rate~\eqref{eq:ini}.
\end{rmk}

\begin{rmk}[Singular-value gap-free clustering] Note that our clustering directly applies to the estimated mean tensor $\hat \tX$ rather than the leading tensor factors $\hat \mU_k$. Applying clustering to the tensor factors suffers from the non-identifiability issue due to the infinitely many orthogonal rotations when the number of blocks $r \geq 3$ in the absence of singular-value gaps. 
Such ambiguity causes the trouble for effective clustering~\citep{abbe2020entrywise}. In contrast, our initialization algorithm applies the clustering to the overall mean tensor $\hat \tX$. This strategy avoids the non-identifiability issue regardless of the number of blocks and singular-value gaps.  
\end{rmk}

\subsection{Angle-based iteration}\label{subsec:angle}
Our Theorem~\ref{thm:initial} has shown the polynomially decaying error rate from our initialization. Now we improve the error rate to exponential decay using local iterations. 
We propose an angle-based local iteration to improve the outputs from Sub-algorithm~\hyperref[alg:main]{1}. 
To gain the intuition, consider an one-dimensional degree-corrected clustering problem with data vectors $\mx_i = \theta(i) \ms_{z(i)} + \mepsilon_i, i \in [p]$, where $\ms_i$'s are known cluster centroids, $\theta(i)$'s are unknown positive degrees, and $z\colon [p] \mapsto [r]$ is the cluster assignment of interest. The angle-based $k$-means algorithm estimates the assignment $z$ by minimizing the angle between data vectors and centroids; i.e., 
 \begin{equation}\label{eq:angle_kmeans}
     z(i) = \argmax_{a \in [r]} \cos ( \mx_i,\ \ms_{a} ),\ \text{ for all }i \in [p].
 \end{equation}
The classical Euclidean-distance based clustering~\citep{han2022exact} fails to recover $z$ in the presence of degree heterogeneity, even under noiseless case. In contrast, the proposed angle-based $k$-means algorithm
achieves accurate recovery without the explicit estimation of $\mtheta$. 

Our Sub-algorithm~\hyperref[alg:main]{2} shares the same spirit as in the angle-based $k$-means. We still take the order-3 tensor for illustration. Specifically, Sub-algorithm~\hyperref[alg:main]{2} updates estimated core tensor and cluster assignment in each iteration. We use superscript $\cdot^{(t)}$ to denote the estimate from the $t$-th iteration, where $t=1, 2, \ldots.$ For core tensor, we consider the following update strategy  
 \[
 \tS^{(t)}(a_1,a_2,a_3)=\text{Ave}\{\tY(i_1,i_2,i_3)\colon z^{(t)}_k(i_k)=a_k, k\in[3]\}.
\]
Intuitively, $\tS^{(t)}$ becomes closer to the true core $\tS$ as $z^{(t)}_k$ is more precise. For cluster assignment, we first aggregate the slices of $\tY$ and obtain the reduced tensor $\tY^{\text{d}}_1 \in \bbR^{p \times r \times r}$ on the first mode with given $z^{(t)}_k$, where
\[
\tY^{\text{d}}_1(i,a_2,a_3)=\text{Ave}\{\tY(i,i_2,i_3)\colon z^{(t)}_k(i_k)=a_k, k \neq 1\}.
\]
Similarly, we also obtain $\tY^{\text{d}}_2, \tY^{\text{d}}_3$.
We use $\mY_k^d$ and $ \mS^{(t)}_k$ to denote the $\text{Mat}_k(\tY^{\text{d}})$ and $\text{Mat}_k(\tS^{(t)})$. The rows $\mY^d_{k,i:}$ and $\mS^{(t)}_{k,a:}$ correspond to the $\mx_i$ and $\ms_a$ in the one-dimensional clustering~\eqref{eq:angle_kmeans}. Then, we obtain the updated assignment by
\[
z_k(i)^{(t+1)} = \argmax_{a \in [r]} \cos \of{ \mY^{\text{d}}_{k,i:}, \mS^{(t)}_{k,a:} },\ \text{ for all } i \in [p],
\]
provided that $\mS^{(t)}_{k,a:}$ is a non-zero vector. Otherwise, if $\mS^{(t)}_{k,a:}$ is a zero vector, then we make the convention to assign $z^{(t+1)}_k(i)$ randomly in $[r]$. Full procedures for our angle-based iteration are described in Sub-algorithm~\hyperref[alg:main]{2}. 

We now establish the misclustering error rate of iterations under the stability assumption. 
\begin{defn}[Locally linear stability] \label{def:stable}
Define the $\varepsilon$-neighborhood of $z$ by $\tN(z,\epsilon)=\{\bar z\colon \ell(\bar z, z)\leq \epsilon\}$. Let $\bar z\colon[p]\to [r]$ be a clustering function. We define two vectors associated with $\bar z$,
\begin{align}
    \mp(\bar z) =(|\bar z^{-1}(1)|, \ldots,|\bar z^{-1}(r)|)^T, \quad
    \mp_{\mtheta}(\bar z) =(\onormSize{}{\mtheta_{\bar z^{-1}(1)}}_1,\ldots,\onormSize{}{\mtheta_{\bar z^{-1}(r)}}_1)^T.
\end{align}

We call the degree is $\varepsilon$-locally linearly stable if and only if 
\begin{equation}\label{eq:local}
    \sin(\mp(\bar z),\ \mp_{\mtheta}(\bar z))\lesssim \varepsilon \Delta_{\min} ,\quad \text{for all } \bar z\in \tN(z, \varepsilon).
\end{equation}
\end{defn}

Roughly speaking, the vector $\mp(\bar z)$ represents the raw cluster sizes, and $\mp_{\theta}(\bar z)$ represents the relative cluster sizes weighted by degrees. 
The local stability holds trivially for $\varepsilon=0$ based on the construction of parameter space~\eqref{eq:family}. The condition~\eqref{eq:local} controls the impact of node degree to the $\mp_{\theta}(\cdot)$ with respect to the misclustering rate $\varepsilon$ and angle gap.  Intuitively, the condition~\eqref{eq:local} controls the skewness of degree so that the angle between raw cluster size and degree-weighted cluster size is well controlled. The stability assumption is proposed for technical convenience, and we relax this condition in numerical studies; see Section~\ref{sec:simulation}.

\begin{thm}[Error for angle-based iteration]\label{thm:refinement} Consider the general sub-Gaussian dTBM with {fixed $r \geq 1$, $K \geq 2$}, independent noise under the parameter space $\tP$, and Assumption~\ref{assmp:min_gap}. {Assume that the locally linear stability of degree holds in the neighborhood $\tN(z,\varepsilon)$ for all $\varepsilon \leq E_0$ and some $E_0 \gtrsim \log^{-1}p $.}
Let $\{z^{(0)}_k\}_{k=1}^K$ be the initialization for Sub-algorithm~\hyperref[alg:main]{2} and $z^{(t)}_k$ be the $t$-th iteration output on the $k$-th mode. Suppose $\min_{i \in [p]}\theta(i) \geq c $ for some {constant} $c > 0$, the $\text{SNR} \geq \tilde C p^{-(K-1)}\log p$ for some sufficiently large positive constant $\tilde C$, and the initialization satisfies 
\begin{equation}
    L(z^{(0)}_k, z) \lesssim \frac{\Delta_{\min}^2}{r \log p}, \quad k \in [K].
\end{equation}
 With probability going to 1 {as $p \rightarrow \infty$}, there exists a contraction parameter $\rho \in (0,1)$ such that 
\begin{align}\label{eq:final}
    \ell(z, \hat z^{(t+1)}_k) \lesssim &\ \KeepStyleUnderBrace{
   \text{SNR}^{-1}
    \exp\of{- \frac{p^{K-1}\text{SNR}}{r^{K-1}}}}_{\substack{\text{statistical error}} }+ \KeepStyleUnderBrace{ \rho^t \ell(z, z^{(0)}_k). }_{\substack{\text{computational error}}}
\end{align}
\normalsize
\end{thm}
From the conclusion~\eqref{eq:final}, we find that the iteration error is decomposed into two parts: statistical error and computational error. The statistical error is unavoidable with noisy data regardless $t$, whereas the computational error decays in an exponential rate as the number of iterations $t \rightarrow \infty$. 

\begin{cor}[Exact recovery of dTBM with weighted higher-order initialization]  Let the initialization $\{z^{(0)}_k\}_{k = 1}^K$ be the output from Sub-algorithm~\hyperref[alg:main]{1}. Assume $\text{SNR} \gtrsim p^{-K/2} \log p$. Combining {all parameter assumptions and the results in} Theorems~\ref{thm:initial} and \ref{thm:refinement}, with probability going to 1 {as $p \rightarrow \infty$}, our estimate $z^{(T)}_k$ achieves exact recovery within polynomial iterations; more precisely,
\begin{equation}
     z^{(T)}_k = \pi_k \circ z, \quad \text{for all }T\gtrsim \log_{1/\rho} p\  \text{and}\ k \in [K ].
\end{equation}
for some permutation $\pi_k \in \Pi$. 
\end{cor}
Therefore, our combined algorithm is \textit{computationally efficient} as long as SNR $\gtrsim p^{-K/2} \log p$. Note that, ignoring the logarithmic term, the minimal SNR requirement, $p^{-K/2}$, coincides with the computational critical value in Theorem~\ref{thm:comp}. Therefore, our algorithm is optimal regarding the signal requirement and lies in the sharpest \emph{computationally efficient} regime in Figure~\ref{fig:phase_axis}. 

\subsection{Extension to Bernoulli observations}\label{subsec:ber}

Bernoulli or network observations are common in multiple fields. Our iteration Theorem~\ref{thm:refinement} holds for Bernoulli models, but our initialization Theorem~\ref{thm:initial} does not. Moreover, our current dTBM is insufficient to address sparsity with decaying mean tensor. Here, we provide extra discussions for Bernoulli initialization and strategies under sparse settings.

\begin{itemize}[wide]
    \item \textit{Extension to dense binary dTBMs.} The main difficulty to establish initialization guarantees for Bernoulli observations lies in the denoising step (lines 1-2 in Sub-algorithm~\hyperref[alg:main]{1}). We now provide a high-level explanation for the technical difficulty when applying Theorem~\ref{thm:initial} to Bernoulli observations.

    The derivation of Theorem~\ref{thm:initial} relies on the upper bound of the estimation error for the mean tensor in Lemma~\ref{lem:two-step_esterror}; i.e., with high probability
\begin{equation}\label{eq:prop1}
    \onormSize{}{\hat \tX - \tX}_F^2 \lesssim p^{K/2},
\end{equation}
where $\tX = \bbE\tY$ and $\hat \tX$ is defined in Step 2 of Sub-algorithm~\hyperref[alg:main]{1}. Unfortunately, the inequality~\eqref{eq:prop1} holds only for i.i.d.\ sub-Gaussian observations, while Bernoulli observations are generally not identically distributed.  

One possible remedy is to apply singular value decomposition to the \emph{square unfolding}~\citep{mu2014square}, $\Mat_{sq}(\cdot)$, of Bernoulli tensor $\tY \in \{0,1\}^{p_1 \times \cdots \times p_K}$. Specifically, the square matricization $\mat_{sq} (\tY) \in \{0,1\}^{p^{\floor{K/2}} \times p^{\ceil{K/2}} }$ has entries $[\Mat_{sq}(\tY)](j_1, j_2) = \tY(i_1, \ldots, i_K)$, where
\begin{align}
    j_1 &= i_1 + p_1(i_2 - 1) + \cdots + p_1 \cdots p_{\floor{K/2}-1} (i_{\floor{K/2}} - 1),\\
    j_2 &= i_{\ceil{K/2}} + p_{ \ceil{K/2}}(i_{\ceil{K/2} + 1} - 1) + \cdots + p_{ \ceil{K/2}} \cdot p_{K-1} (i_K-1).
\end{align}
The matrix $\Mat_{sq}(\tY)$ is asymmetric. We interpret $\Mat_{sq}(\tY)$ as the adjacency matrix for a bipartite network with connections between two groups of nodes. The two groups of nodes in the bipartite network have $p_1 \cdots p_{\floor{K/2}}$ and $ p_{\ceil{K/2}} \cdots p_K$ nodes, respectively. The entry $[\Mat_{sq}(\tY)](j_1, j_2)$ refers to the presence of connection between the nodes indexed by combinations $(i_1, \ldots, i_{\floor{K/2}})$ and $(i_{\ceil{K/2}}, \ldots, i_K)$. We summarize the procedure in Sub-algorithm 3.

\begin{algorithm}[h!]
\caption*{\bf Sub-algorithm 3: Weighted higher-order initialization for Bernoulli observation}
\vspace{.15cm}
\begin{algorithmic}[1] 
\INPUT Bernoulli tensor $\tY \in \{0,1\}^{p\times \cdots \times p}$, cluster number $r$, relaxation factor $\eta > 1$ in $k$-means clustering.
\State  Let the matrix $\mat_{sq} (\tY) \in \{0,1\}^{p^{\floor{K/2}} \times p^{\ceil{K/2}} }$ denote the nearly square unfolded tensor. Compute the estimate $\tX'$, where
\begin{equation}\label{eq:matrixsvd}
    \hat \tX' = \argmin_{\text{rank}(\mat_{sq}(\tX)) \leq r^{\ceil{K/2}}} \onormSize{}{ \mat_{sq}(\tX) -  \mat_{sq}(\tY)}_F^2.
\end{equation}
\normalsize
\State Implement lines 3-5 of Sub-algorithm~\hyperref[alg:main]{1} with $\hat \tX$ replaced by $\hat \tX'$ in~\eqref{eq:matrixsvd}.
\OUTPUT Initial clustering $z^{(0)}_k \leftarrow \hat z_k, k \in [K]$.
\end{algorithmic}
\end{algorithm}

\begin{prop}[Error for Bernoulli initialization]\label{prop:ber} Consider the Bernoulli dTBM in the parameter space $\tP$ {with fixed $r \geq 1, K \geq 2$}. Assume that Assumption~\ref{assmp:min_gap} holds, $\mtheta$ is balanced, and $\min_{i\in[p]}\theta(i) \geq c$ for some constant $c>0$. Let $ z^{(0)}_k$ denote the output of Sub-algorithm~3. With probability going to 1 {as $p \rightarrow \infty$}, we have
\begin{equation}\label{eq:ini_b}
 \ell(z^{(0)}_k, z_k) \lesssim \frac{r^K p^{- \floor{K/2} }}{\text{SNR}}, \quad \text{and} \quad L(z^{(0)}_k, z_k) \lesssim  {\sigma^2 r^K p^{-\floor{K/2}}}.
\end{equation}
\end{prop}

\begin{rmk}[Comparison with Gaussian model] The Bernoulli bound $\tO(p^{- \floor{K/2} })$ in Proposition~\ref{prop:ber} is relatively looser than the Gaussian bound $\tO(p^{- K/2 })$ in Theorem~\ref{thm:initial}. The gap between Bernoulli and Gaussian error decreases as the order $K$ increases. Nevertheless, combining with angle iteration Sub-algorithm~\hyperref[alg:main]{2}, Bernoulli clustering still achieves exponential error rate $\exp \of{ - p^{(K-1)}}$ at a price of a larger SNR.
The investigation of the gap between upper bound $p^{- \floor{K/2}}$ and the lower bound $p^{- K/2}$ for Bernoulli tensors will be left as future work. In numerical experiments, we will use our original initialization, Sub-algorithm 1, to verify the robustness to Bernoulli observations.
\end{rmk}

\begin{rmk}[Comparison with previous methods]  Previous work \citep{ke2019community} develops a spectral clustering method for Bernoulli dTBM. \cite{ke2019community} adopts a different signal notion based on the singular gap in the core tensor, denoted as $\Delta_{\rm singular}$. By \citet[Theorem 1]{ke2019community}, the spectral method achieves exact recovery with $\Delta_{\rm singular} \gtrsim p^{-1/2}$. However, we are not able to infer the exact recovery of spectral method by our angle-base SNR condition. Consider an order-2 dTBM with $p > 2, \sigma^2 = 1$, $\mtheta = \mathbf{1}_p$, equal size assignment $|z^{-1}(a)| = p/r$ for all $ a \in [r]$, and core matrix equal to the 2-dimensional identity matrix $\mS = \mI_2$. The singular gap under this setting is $\Delta_{\rm singular} = \min\{ \lambda_1 - \lambda_2, \lambda_2 \} = 0$, where $\lambda_1 \geq \lambda_2$ are singular values of $\mS$. In contrast, our angle gap $\Delta_{\min}^2 = 2$ satisfies the SNR condition in Theorem~\ref{thm:refinement}. Then, our algorithm achieves the exact recovery, but the spectral method in \cite{ke2019community} fails.

Hence, for fair comparison, we compare the best performance of our algorithm and \cite{ke2019community} under the strongest signal setting of each model. Since both methods contain an iteration procedure, we set the iteration number to infinity to avoid the computational error. Considering the largest angle-based SNR $\asymp 1$ in Theorem~\ref{thm:refinement}, our Bernoulli clustering achieves exponential error rate of order $\exp(-p^{(K-1)})$; considering the largest singular gap $\Delta_{\rm singular} \asymp 1$ in Theorem 1 of \cite{ke2019community}, the spectral clustering has a polynomial error rate of order $p^{-2}$. Our algorithm still shows a better theoretical accuracy than the competitive work for Bernoulli observations. 
\end{rmk}

\item \textit{Extension to sparse binary dTBMs.} The sparsity is often a popular feature in hypergraphs~\citep{florescu2016spectral,ke2019community, ahn2018hypergraph}. Specifically, the sparse binary dTBM assumes that, the entries of $\tY$ follow independent Bernoulli distributions with the mean
     \begin{equation}\label{eq:sparse_dtbm}
    \bbE\tY = \alpha_p \tS \times_1 \mTheta \mM \times_2  \cdots \times_K \mTheta \mM,
\end{equation}
where the extra scalar parameter $\alpha_p \in (0,1]$ is function of $p$ that controls the sparsity. A smaller $\alpha_p$ indicates a higher level of sparsity. Our current work focuses on dense dTBM with $\alpha_p=1$.  While sparse dTBM is an interesting application, the algorithm and its analysis require different techniques. Below, we discuss possible modifications of the algorithm.

The sparsity affects our initialization guarantee in our Theorem~\ref{thm:initial}. In our initialization, the spectral denoising step (lines 1-2 in Sub-algorithm~\hyperref[alg:main]{1}) implements matrix SVD to unfolded tensors. 
However, SVD-based methods are believed to fail in extremely sparse SBM due to the localization phenomenon in the singular vectors \citep{florescu2016spectral}. Inspired by \cite{florescu2016spectral}, we adopt the diagonal-deleted HOSVD (D-HOSVD) \citep{ke2019community} as the initialization in our higher-order clustering. 

The sparsity also affects the iteration guarantee in our Theorem~\ref{thm:refinement}. The decaying mean tensor leads to a worse statistical error of order $\tO( -\alpha_p p^{K-1})$ on $\hat \tX$. The theoretical analyses for sparse binary dTBM and algorithms are left as future directions. Instead, we add numerical experiments to evaluate the robustness of our algorithm and the improvement of D-HOSVD initialization in the sparse dTBM; see Appendix A.

\end{itemize}
\subsection{Practical issues}~\label{subsec:exten}

{\bf Computational complexity.} Our two-stage algorithm has a computational cost polynomial in tensor dimension $p$. Specifically, the complexity of Sub-algorithm~\hyperref[alg:main]{1} is $\tO(K p^{K+1} + Krp^K)$, where the first term is contributed by the double projection and the calculation of $\hat \tX$, and the second term comes from normalization and the $k$-means. The cost of each update in Sub-algorithm~\hyperref[alg:main]{2} is $\tO(p^K + pr^K)$, where $p^K$ comes from the calculation of $\tS^{(t)}$ and $\tY^{\text{d}}_k$, and $pr^K$ comes from the normalization of $\tY^{\text{d}}_k$, the calculation of $\tS^{(t)}$, and the cluster assignment update in Step~13.

{\bf Hyper-parameter selection.} In our theoretical analysis, we have assumed the true cluster number $r$ is given to our algorithm. In practice, the cluster number $r$ is often unknown, and we now propose a method to choose $r$ from data. We impose the Bayesian information criterion (BIC) and choose the cluster number that minimizes BIC; i.e., under the symmetric Gaussian dTBM~\eqref{eq:model_tensor},
\begin{align}\label{eq:BIC}
 \hat r = \argmin_{r \in \mathbb{Z}_+ } \of{p^K \log(\onormSize{}{\hat \tX - \tY}_F^2)  + p_e(r) K \log p},
\end{align}
with $\hat \tX = \hat \tS(r) \times_1 \hat \mTheta(r) \hat \mM(r) \times_2 \cdots \times_K \hat \mTheta(r) \hat \mM(r),$ where the triplet $(\hat z(r), \hat \tS(r),\hat \mtheta(r))$ are estimated parameters with cluster number $r$, and $p_e (r)= r^K + p(\log r + 1) - r$ is the effective number of parameters. Note that we have added the argument $(r)$ to related quantities as functions of $r$. In particular, the estimate $\hat \mtheta(r)$ in~\eqref{eq:BIC} is obtained by first calculating the reduced tensor $\hat \tY^{\text{d}}$ with $\hat z(r)$, and then normalizing the row norms $\onormSize{}{\hat \mY_{i:}^{\text{d}}}$ to 1 in each cluster; i.e., 
\begin{equation}
 \hat \mtheta(r)=(\hat \theta(1,r),\ldots,\hat \theta(p,r))^T,
\end{equation}
with $\hat \theta(i,r) = {\onormSize{}{\hat \mY^{\text{d}}(r)_{i:}}}/{\sum_{j: \hat z(j,r) =  \hat z(i,r)} \onormSize{}{\hat \mY^{\text{d}}(r)_{j:}}}$, $\hat \mY^{\text{d}}(r) = \mat(\hat \tY^{\text{d}}(r))$, $\hat \tY^{\text{d}}(r)(i,a_2,\ldots,a_K) =  \text{Ave}\{\tY(i,i_2,\ldots,i_K): \hat z(i_k, r) = a_k, k \neq 1 \}$, and $\hat z(i,r)$ denotes the community label for the $i$-th node with given cluster number $r$. We evaluate the performance of the BIC criterion in Section~\ref{subsec:num_theory}.

\section{Comparison with non-degree tensor block model} \label{sec:tbm}

We discuss the connections and differences between dTBM and TBM~\citep{han2022exact} from three aspects: signal notions, theoretical results, and algorithms. Without loss of generality, let $\sigma^2=1$. 

\begin{itemize}[wide]
    \item \textit{Signal notion.} The signal levels in both TBM~\citep{han2022exact} and our dTBM are functions of the core tensor $\tS$. We emphasize that the signal notions are different between the two models. In particular, the Euclidean-based signal notion in TBM~\cite{han2022exact} fails to accurately describe the phase transition in our dTBM due to the possible heterogeneity in degree $\mtheta$. To compare, we denote our angle-based signal notion in~\eqref{eq:gamma} and the Euclidean-based SNR in \cite{han2022exact} as $\Delta_{\text{ang}}^2$ and $\Delta_{\text{Euc}}^2$, respectively:
    \begin{align}
        \Delta_{\text{ang}}^2 =  2(1 - \max_{a \neq b\in [r]}\cos \of{\mS_{a:},\  \mS_{b:}} ), \quad 
 \Delta_{\text{Euc}}^2 = \min_{a \neq b \in [r]} \onormSize{}{\mS_{a:} - \mS_{b:}}^2.
    \end{align}
By Lemma~\ref{lem:norm_diff} in the Appendix B, we have 
\begin{equation}\label{eq:signalcompare}
     \Delta_{\text{ang}}^2  \max_{a \in [r]}\onormSize{}{\mS_{a:}}^2 \leq \Delta_{\text{Euc}}^2.
\end{equation}
The above inequality indicates that the condition $\Delta_{\text{Euc}}^2 \leq p^{\gamma}$ is sufficient but not necessary for $\Delta_{\text{ang}}^2 \leq p^{\gamma}$. In fact, if we were to use $\Delta_{\text{Euc}}^2$ for both models, then the phase transition of dTBM can be arbitrarily worse than that for TBM.

Here, we provide an example to illustrate the dramatical difference between TBM and dTBM with the same core tensor.  

\begin{example}[Comparison with Euclidean-based signal notion] \label{example:euc_alg} Consider a biclustering model with $\mtheta=\mathbf{1}$ and an order-2 core matrix 
\begin{equation}
    \mS = \begin{pmatrix} p^{(\gamma+1)/2 } + 2  & 2 p^{(\gamma+1)/2} + 4\\
    2 & 4
    \end{pmatrix},\quad \text{with}\ \gamma \leq -1.
\end{equation}
The core matrix $\mS$ lies in the parameter spaces of TBM and our dTBM. Here, the constraint $\gamma \leq -1$ is added to ensure the bounded condition of $\mS$ in our parameter space in \eqref{eq:family}. The angle-based and Euclidean-based signal levels of $\mS$ are 
\begin{equation}
    \Delta_{\text{ang }}^2(\mS) = 0 \ \left(\leq p^{\gamma}\right), \quad \Delta_{\text{Euc}}^2(\mS) = 5 p^{\gamma + 1} \ \left(\geq p^{\gamma}\right).
\end{equation}
We conclude that TBM with $\mS$ achieves exact recovery with a polynomial-time algorithm; see \citet[Theorem 4]{han2022exact}. By contrast, the dTBM with the same $\mS$ and input $r=2$ violets the identifiability condition, and thus fails to be solved by all estimators; see our Theorem~\ref{thm:unique}. 
\end{example}
    
    \item \textit{Theoretical results.} In both works, we study the phase transition of TBM and dTBM with respect to the Euclidean and angle-based SNRs. We briefly summarize the results in \cite{han2022exact} and compare with ours. 

    \textit{Statistical critical value:}
    \small
    \begin{align}
        \text{Ours:}& \ \Delta_{\text{ang}}^2 \lesssim p^{-(K-1)} \Rightarrow \text{statistically impossible;} \quad \Delta_{\text{ang}}^2 \gtrsim   p^{-(K-1)} \Rightarrow \text{MLE achieves exact recovery;} \\
        \text{Han's:}& \ \Delta_{\text{Euc}}^2 \lesssim p^{-(K-1)} \Rightarrow \text{statistically impossible;} \quad \Delta_{\text{Euc}}^2 \gtrsim   p^{-(K-1)} \Rightarrow \text{MLE achieves exact recovery}.
    \end{align}
    \normalsize
    
     \textit{Computational critical value:}
     \small
    \begin{align}
        \text{Ours:}& \ \Delta_{\text{ang}}^2 \lesssim p^{-K/2} \Rightarrow \text{computationally impossible;} \quad \Delta_{\text{ang}}^2 \gtrsim   p^{-K/2} \Rightarrow \text{ computationally efficient;} \\
        \text{Han's:}& \ \Delta_{\text{Euc}}^2 \lesssim p^{-K/2} \Rightarrow \text{computationally impossible;} \quad \Delta_{\text{Euc}}^2 \gtrsim   p^{-K/2} \Rightarrow \text{computationally efficient}.
    \end{align}
    \normalsize
    
    

    The above comparison reveals four major differences.

    First, none of our results in Section~\ref{sec:limits} are corollaries of \cite{han2022exact}. Both models show the similar conclusion but under different conditions. While the TBM impossibility~\citep{han2022exact} provides a necessary condition for our dTBM impossibility, we find that such a condition is often loose. There exists a regime of $\tS$ in which TBM problems are computationally efficient but dTBM problems are statistically impossible; see Example~\ref{example:euc_alg}. This observation has motivated us to develop the new signal notion $\Delta^2_{\text{ang}}$ for sharp dTBM phase transition conditions.  
     
 Second, to find the phase transition, we need to show both the impossibility and achievability when SNR is below and above the critical value, respectively. While the TBM impossibility can serve as a loose condition of our dTBM impossibility, more efforts are required to show the achievability. In particular, since TBM is a more restrictive model than dTBM, the achievability in \cite{han2022exact} does not imply the achievability of dTBM in a larger parameter space. The latter requires us to develop new MLE and polynomial algorithms for dTBM achievability.

 Third, from the perspective of proofs, we develop new dTBM-specific techniques to handle the extra degree heterogeneity. In our Theorem~\ref{thm:stats}, we construct a special non-trivial degree heterogeneity to establish the lower bound for arbitrary core tensor with small angle gap, while, TBM~\citep{han2022exact} considers the constructions without degree parameter.
In our Theorem~\ref{thm:comp}, we construct a rank-2 tensor to relate HPC conjecture to $\Delta^2_{\text{ang}}$, while TBM~\citep{han2022exact} constructs a rank-1 tensor to relate HPC conjecture to $\Delta^2_{\text{Euc}}$. The asymptotic non-equivalence between $\Delta^2_{\text{ang}}$ and $\Delta^2_{\text{Euc}}$ renders our proof technically more involved.
 
 Last, we discuss the statistical impossibility statements. Our Theorem~\ref{thm:stats} implies the statistical impossibility whenever the core tensor $\tS$ leads to an angle-based SNR below the critical value, while, Theorem 6 in \cite{han2022exact} implies the worst case statistical impossibility for a particular core tensor $\tS$ with Euclidean-based SNR below the statistical limit. Hence, our Theorem~\ref{thm:stats} shows a stronger statistical impossibility for dTBM than that presented in TBM \citet[Theorem 6]{han2022exact}. However, inspecting the proof of \cite{han2022exact}, the proof of Theorem 6 indeed implies a stronger TBM impossibility statement for arbitrary core tensor; i.e., when $\gamma < -(K-1)$
 \begin{equation}
     \liminf_{p \rightarrow \infty} \inf_{ \tS \in \tP_{\tS, {\rm TBM}} \cap \{\Delta_{\rm Euc}^2 = p^{\gamma} \} }\inf_{\hat z_{\rm stats}}  \sup_{z \in \tP_{z, {\rm TBM}}}  \bbE[p \ell(\hat z_{\rm stats}, z)] \geq 1,
 \end{equation}
 \normalsize
 where $\tP_{\tS, {\rm TBM}}$ and $\tP_{z, {\rm TBM}}$ refer to the space for core tensor $\tS$ and assignment $z$ under TBM, respectively. Again, in terms of the strong statistical impossibility, both models show the similar conclusion but under different conditions. Since two impossibilities consider different core tensor regimes with non-equivalent $\Delta^2_{\text{ang}}$ and $\Delta^2_{\text{Euc}}$, we emphasize that different proof techniques are required to obtain these similar conclusions. See our proof sketch in Section~\ref{subsec:impossible}, Appendices~\ref{sec:statprove1} and \ref{sec:compprove1} for detail technical differences.

 \item \textit{Algorithms.} Both \cite{han2022exact} and our work propose the two-step algorithm, which combines warm initialization and iterative refinement to achieve exact recovery. This local-to-global strategy is not new in clustering literature~\citep{gao2022iterative, chien2019minimax}. The highlight of our algorithm is the angle-based update in lines 10-14, Sub-algorithm~\hyperref[alg:main]{2}, which is specifically designed for dTBM to avoid the estimation of $\mtheta$. This angle-based update brings new proof challenges. We develop polar-coordinate based techniques to establish the error rate for the proposed algorithm. 
\end{itemize}

\section{Numerical studies}\label{sec:simulation}

 We evaluate the performance of the weighted higher-order initialization and angle-based iteration in this section. We report average errors and standard deviations across 30 replications in each experiment. Clustering accuracy is assessed by clustering error rate (CER, i.e., one minus rand index). The CER between $(\hat z, z)$ is equivalent to misclustering error $\ell(\hat z, z)$ up to constant multiplications~\citep{meilua2012local}, and a lower CER indicates a better performance.

We generate order-3 tensors with \emph{assortative}~\citep{gao2018community} core tensors to control SNR; i.e., we set $\tS_{aaa} = s_1$ for $a \in [r]$ and others be $s_2$, where $s_1 > s_2 > 0$. Let $\alpha = s_1/s_2$. We set $\alpha$ close to 1 such that $1-\alpha=o(p)$. In particular, we have $\alpha = 1 + \Omega(p^{\gamma/2})$ with $\gamma<0$ by Assumption~\ref{assmp:min_gap} and definition~\eqref{eq:gamma}. Hence, we easily adjust SNR via varying $\alpha$. The assortative setting is proposed for simulations, and our algorithm is applicable for general tensors in practice. The cluster assignment $z$ is randomly generated with equal probability across  $r$ clusters for each mode. Without further explanation, we generate degree heterogeneity $\mtheta$ from absolute normal distribution by $\theta(i) = |X_i| + 1 - 1/\sqrt{2\pi}$ with $|X_i| \stackrel{\text{i.i.d.}}\sim N(0,1), i \in [p]$ and normalize $\mtheta$ to satisfy \eqref{eq:family}. Also, we set $\sigma^2 = 1$ for Gaussian data without further specification. 

\subsection{Verification of theoretical results}\label{subsec:num_theory}

\begin{figure}[htb]
    \centering
    \includegraphics[width=.85\columnwidth]{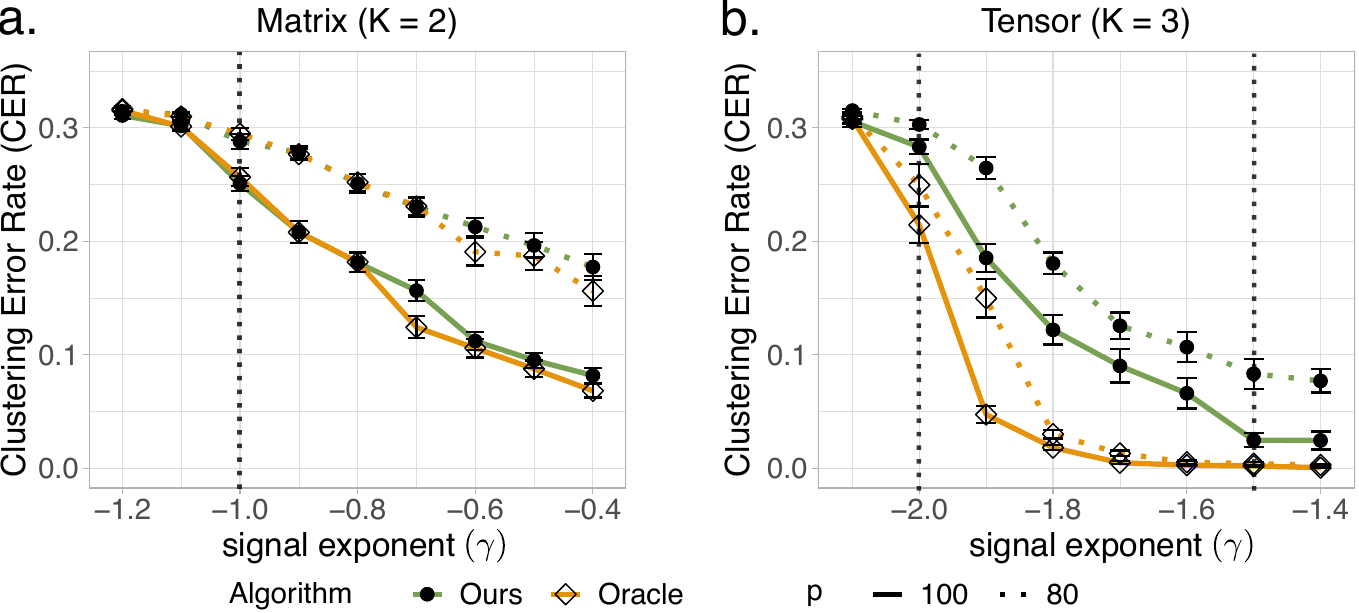}
    \caption{SNR phase transitions for clustering in dTBM with $p = \{80, 100\}, r = 5$ under (a) matrix case with $\gamma \in [-1.2, -0.4]$ and (b) tensor case with $ \gamma \in [-2.1, -1.4]$.
    }
    \label{fig:phase}
\end{figure}

The first experiment verifies  statistical-computational gap described in Section~\ref{sec:limits}. Consider the Gaussian model with $p = \{80, 100\}$, $r = 5$. We vary $\gamma $ in $ [-1.2, -0.4]$ and $[-2.1, -1.4]$ for matrix ($K=2$) and tensor $(K = 3)$ clustering, respectively. Note that finding MLE under dTBM is computationally intractable. We approximate MLE using an oracle estimator, i.e., the output of Sub-algorithm~\hyperref[alg:main]{2} initialized from true assignment. Figure~\ref{fig:phase}a shows that both our algorithm and oracle estimator start to decrease around the critical value $\gamma_{\text{stat}}  = \gamma_{\text{comp}}  = -1$ in matrix case. In contrast, Figure~\ref{fig:phase}b shows a significant gap in the phase transitions between the algorithm estimator and oracle estimator in tensor case. The oracle error rapidly decreases to 0 when $\gamma_{\text{stat}} = -2$, whereas the algorithm estimator tends to achieve exact clustering when $\gamma_{\text{comp}} = -1.5$. Figure~\ref{fig:phase} confirms the existence of the statistical-computational gap in our Theorems~\ref{thm:stats} and~\ref{thm:comp}. 

\begin{figure}[htp!]
    \centering
     \includegraphics[width=.85\columnwidth]{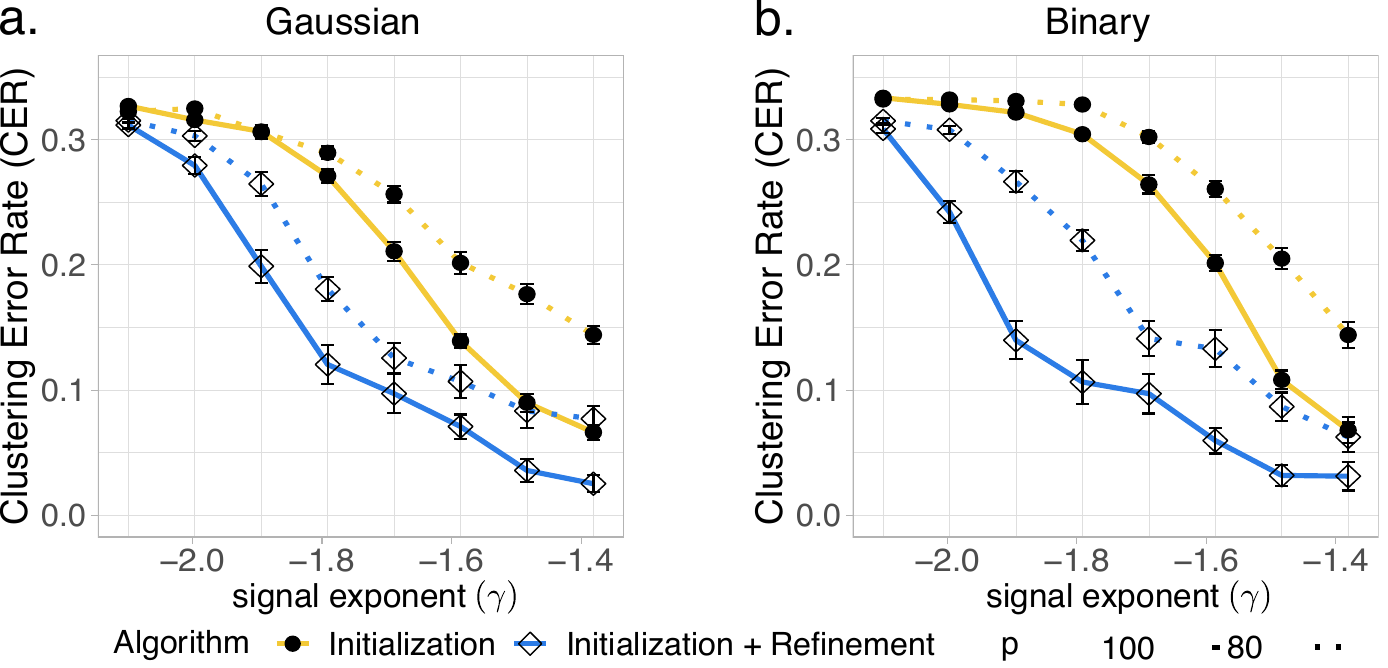}
    \caption{CER versus signal exponent $(\gamma)$ for initialization only and for combined algorithm. We set $p = \{80, 100\}, r = 5, \gamma \in [-2.1, -1.4]$ under (a) Gaussian models and (b) Bernoulli models. }
    \label{fig:ini_re}
\end{figure}

The second experiment verifies the performance guarantees of two algorithms: (i) weighted higher-order initialization; (ii) combined algorithm of weighted higher-order initialization and angle-based iteration. We consider both the Gaussian and Bernoulli models with $p = \{80, 100\}$, $r = 5$, $\gamma \in [-2.1, -1.4]$. Figure~\ref{fig:ini_re} shows the substantial improvement of combined algorithm over initialization, especially under weak and intermediate signals. This phenomenon agrees with the error rates in Theorems~\ref{thm:initial} and \ref{thm:refinement} 
and confirms the necessity of the local iterations.

\begin{table}[hbt]
\centering
    \begin{tabularx}{\columnwidth}{c *{8}{Y}}
    \toprule
    Settings & \multicolumn{2}{c}{\small  $p = 50, \sigma^2 = 0.25$} & \multicolumn{2}{c}{\small $p = 50, \sigma^2 = 1$} & \multicolumn{2}{c}{\small  $p = 80, \sigma^2 = 0.25$} & \multicolumn{2}{c}{\small  $p = 80, \sigma^2 = 1$}\\
    \cmidrule(lr){2-3} \cmidrule(l){4-5} \cmidrule(l){6-7} \cmidrule(l){8-9}
         True cluster number $r$ &  \small 2 & \small 4 &  \small 2  &  \small 4 &  \small2  & \small 4 &  \small 2 & \small 4  \\
         \midrule
         Estimated cluster number $\hat r$ &   \small 2(0)  &  \small 3.9(0.2)&  \small 2(0)    &   \small 3.1(0.5) &  \small 2(0)    &  \small 4(0)   &  \small 2(0)    &  \small 3.9(0.3)   \\
     \bottomrule
    \end{tabularx}
    \caption{Estimated cluster number given by BIC criterion under the low noise level $(\sigma^2 = 0.25)$ and high noise level $(\sigma^2 = 0.5)$ settings. Numbers in parentheses are standard deviations of $\hat r$ over 30 replications.}
    \label{tab:select}
\end{table}

The third experiment evaluates the empirical performance of the BIC criterion to select unknown cluster number. We generate the data from an order-3 Gaussian model with $p = \{50,80\}$, $r = \{2,4\}$, and noise level $\sigma^2 \in \{ 0.25,1\}$. Table~\ref{tab:select} shows that our BIC criterion well chooses the true $r$ under most settings.  Note that the BIC slightly underestimates the true cluster number $(r = 4)$ with smaller dimension and higher noise $(p = 50, \sigma^2=1)$, and the accuracy immediately increases with larger dimension $p = 80$. The improvement follows from the fact that a larger dimension $p$ indicates a larger sample size in the tensor block model. Therefore, we conclude that BIC criterion is a reasonable way to tune the cluster number.

\subsection{Comparison with other methods}\label{subsec:comp}

We compare our algorithm with following higher-order clustering methods:
\begin{itemize}[wide,topsep=-3pt,itemsep=0pt,parsep=1pt]
    \item \textbf{\small HOSVD}: HOSVD on data tensor and $k$-means on the rows of the factor matrix;
    \item \textbf{\small HOSVD+}: HOSVD on data tensor and $k$-means on the $\ell_2$-normalized rows of the factor matrix;
    \item \textbf{\small HLloyd}~\citep{han2022exact}: High-order clustering algorithm developed for non-degree tensor block models;
    \item \textbf{\small SCORE}~\citep{ke2019community}: Tensor-SCORE for clustering developed for sparse binary tensors.
\end{itemize}

Among the four alternative algorithms, the \textbf{\small SCORE} is the closest method to ours. We set the tuning parameters of \textbf{\small SCORE} as in previous literature \citep{ke2019community}. The methods \textbf{\small SCORE} and \textbf{\small HOSVD+} are designed for degree models, whereas \textbf{\small HOSVD} and \textbf{\small HLloyd} are designed for non-degree models. We conduct two experiments to assess the impacts of (i) signal strength and (ii) degree heterogeneity, based on Gaussian and Bernoulli models with $ p = 100, r = 5$. We refer to our algorithm as \textbf{\small dTBM} in the comparison. 

\begin{figure}[h!]
    \centering
    \includegraphics[width=.85\columnwidth]{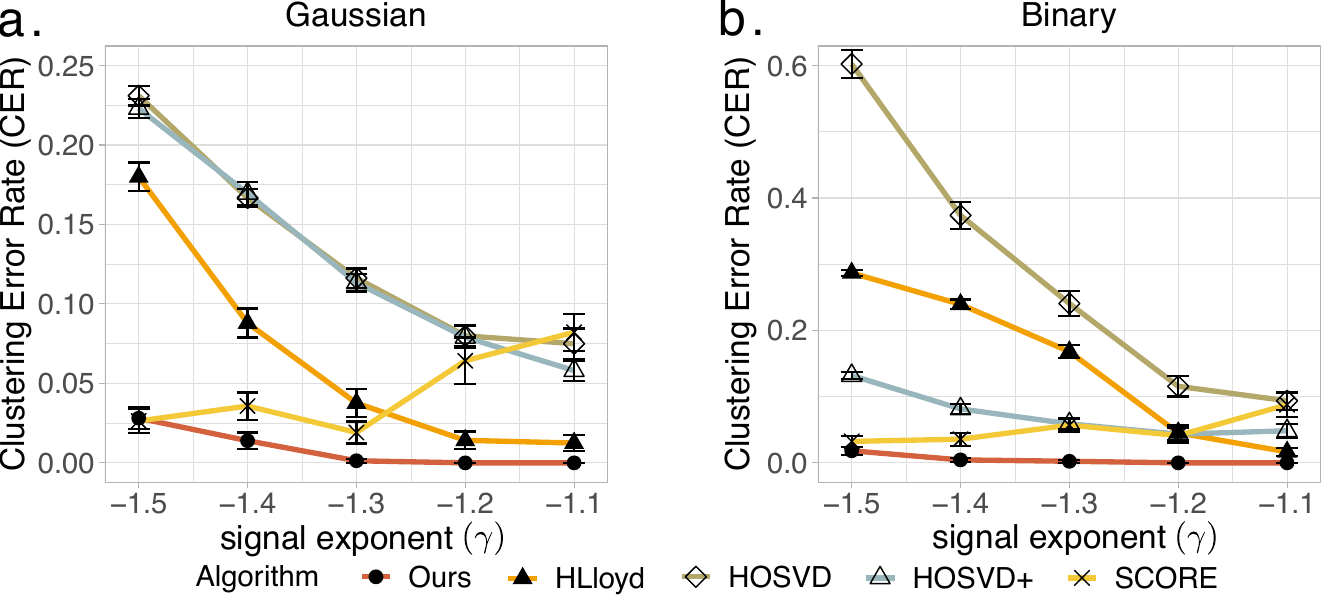}
    \caption{CER versus signal exponent (denoted $\gamma$) for different methods. We set $p = 100, r = 5, \gamma \in [-1.5, -1.1]$ under (a) Gaussian and (b) Bernoulli models.}
    \label{fig:comp_gamma}
\end{figure}

We investigate the effects of signal to clustering performance by varying $\gamma \in [-1.5, -1.1]$. Figure~\ref{fig:comp_gamma} shows that our method \textbf{\small dTBM} outperforms all other algorithms. The sub-optimality of \textbf{\small SCORE} and \textbf{\small HOSVD+} indicates the necessity of local iterations on the clustering. Furthermore,  Figure~\ref{fig:comp_gamma} shows the inadequacy of non-degree algorithms in the presence of mild degree heterogeneity. 
The experiment demonstrates the benefits of addressing heterogeneity in higher-order clustering tasks.   

\begin{figure}[htp!]
    \centering
    \includegraphics[width=.85\columnwidth]{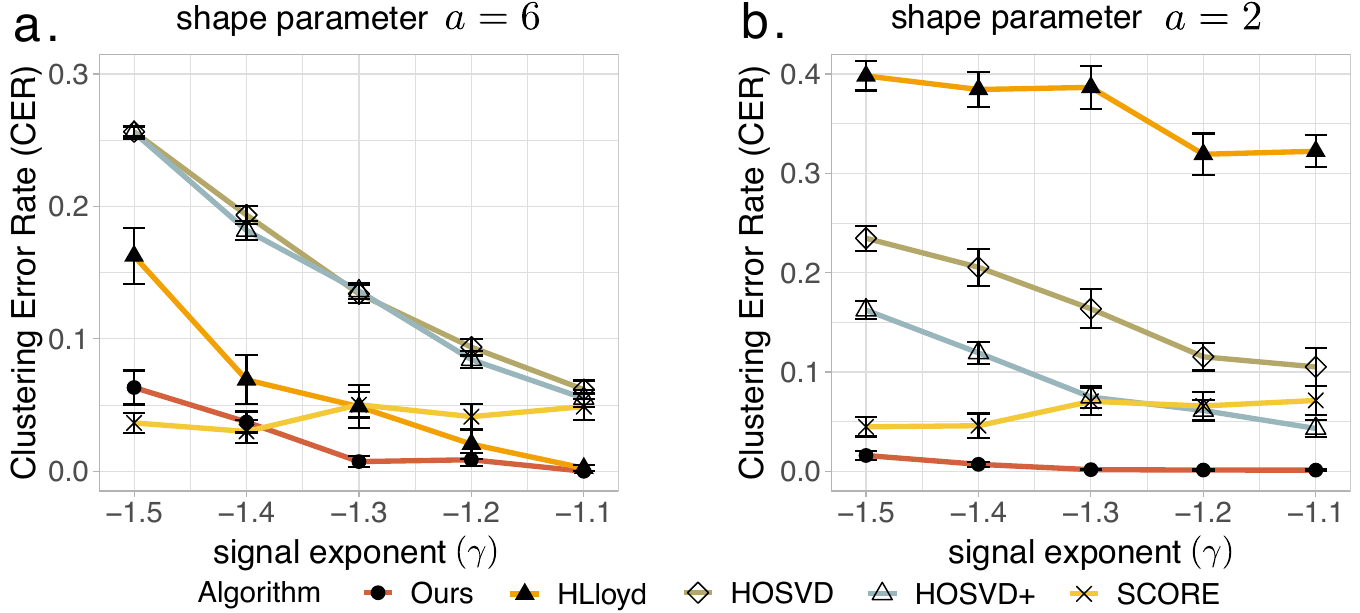}
    \caption{CER comparison versus signal exponent (denoted $\gamma$) under (a) low (shape parameter $a = 6$)  (b) high (shape parameter $a = 2$) degree heterogeneity. We set $p = 100, r = 5, \gamma \in [-1.5, -1.1]$ under Gaussian model.}
    \label{fig:comp_gamma_theta}
\end{figure}

The only exception in Figure~\ref{fig:comp_gamma} is the slightly better performance of \textbf{\small HLloyd} over \textbf{\small HOSVD+} under Gaussian model. However, we find the advantage of \textbf{\small HLloyd} disappears with higher degree heterogeneity. We perform extra simulations to verify the impact of degree effects. We use the same setting as in the first experiment in the Section~\ref{subsec:comp}, except that we now generate the degree heterogeneity $\mtheta$ from Pareto distribution prior to normalization. The density function of Pareto distribution is $f(x|a,b) = a b^a x^{-(a+1)} \ind\{ x \geq b \}$, where $a$ is called \emph{shape} parameter. We vary $a \in \{2,6\}$ and choose $b$ such that $\bbE X = a(a - 1)^{-1}b = 1$ for $X$ following Pareto$(a,b)$. Note that a smaller $a$ leads to a larger variance in $\mtheta$ and hence a larger degree heterogeneity. We consider the Gaussian model under low $(a = 6)$ and high $(a = 2)$ degree heterogeneity. Figure~\ref{fig:comp_gamma_theta} shows that the errors for non-degree algorithms (\textbf{\small HLloyd}, \textbf{\small HOSVD}) increase with degree heterogeneity. In addition, the advantage of \textbf{\small HLloyd} over \textbf{\small HOSVD+} disappears with higher degree heterogeneity. 

\begin{figure}[h!]
    \centering
    \includegraphics[width=.85\columnwidth]{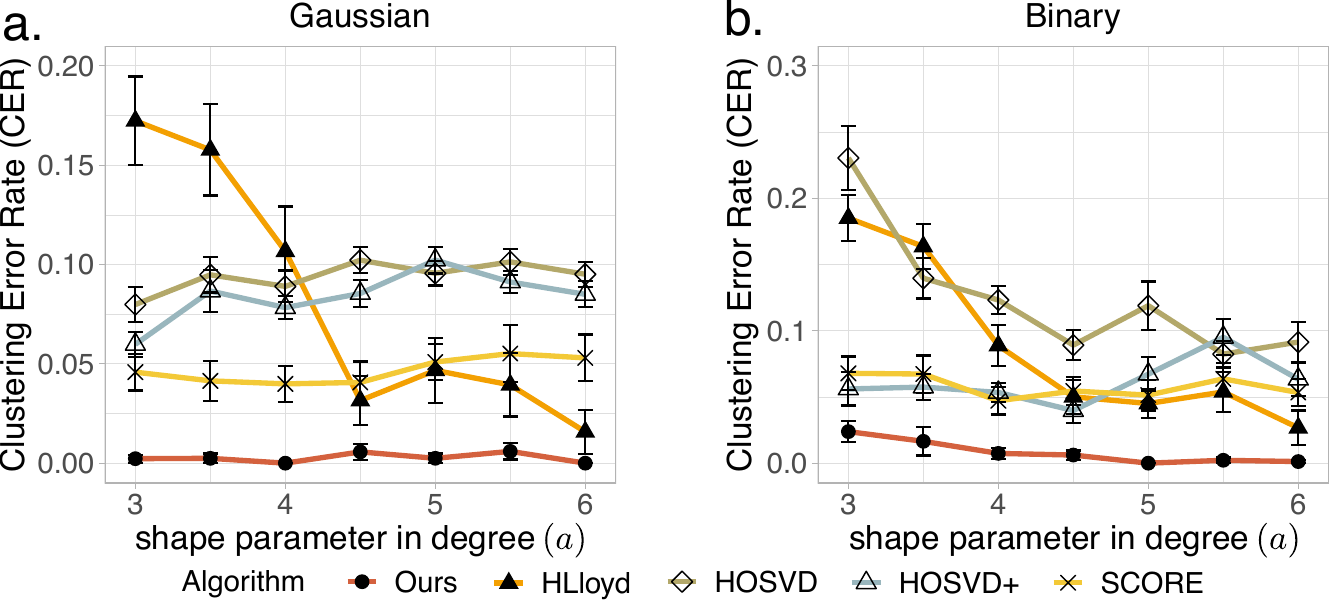}
    \caption{CER versus shape parameter in degree (denoted $a\in[3,6]$) for different methods. We set $p = 100, r = 5, \gamma = -1.2$ under (a) Gaussian and (b) Bernoulli models.}
    \label{fig:comp_theta}
\end{figure}

The last experiment investigates the effects of degree heterogeneity to clustering performance. We fix the signal exponent $\gamma = -1.2$ and vary the extent of degree heterogeneity. In this experiment, we generate $\mtheta$ from Pareto distribution prior to normalization. We vary the shape parameter $a \in [3,6]$ in the Pareto distribution to investigate a range of degree heterogeneities. Figure~\ref{fig:comp_theta} demonstrates the stability of degree-corrected algorithms (\textbf{\small dTBM}, \textbf{\small SCORE}, \textbf{\small HOSVD+}) over the entire range of degree heterogeneity under consideration. In contrast, non-degree algorithms (\textbf{\small HLloyd}, \textbf{\small HOSVD}) show poor performance with large heterogeneity, especially in Bernoulli cases. This experiment, again, highlights the benefit of addressing degree heterogeneity in higher-order clustering.

\section{Real data applications}\label{sec:real}
\subsection{Human brain connectome data analysis}

The Human Connectome Project (HCP) aims to construct the structural and functional neural connections in human brains~\citep{van2013wu}. We preprocess the original dataset following \cite{desikan2006automated} and partition the brain into 68 regions. The cleaned dataset includes brain networks for 136 individuals. Each brain network is represented by a 68-by-68 binary symmetric matrix, where the entry with value 1 indicates the presence of connection between node pairs, while the value 0 indicates the absence. We use $\tY \in \{0,1\}^{68 \times 68 \times 136}$ to denote the binary tensor. Individual attributes such as gender and sex are recorded.

\begin{figure*}[htb]
    \centering
    \includegraphics[width = .7\textwidth]{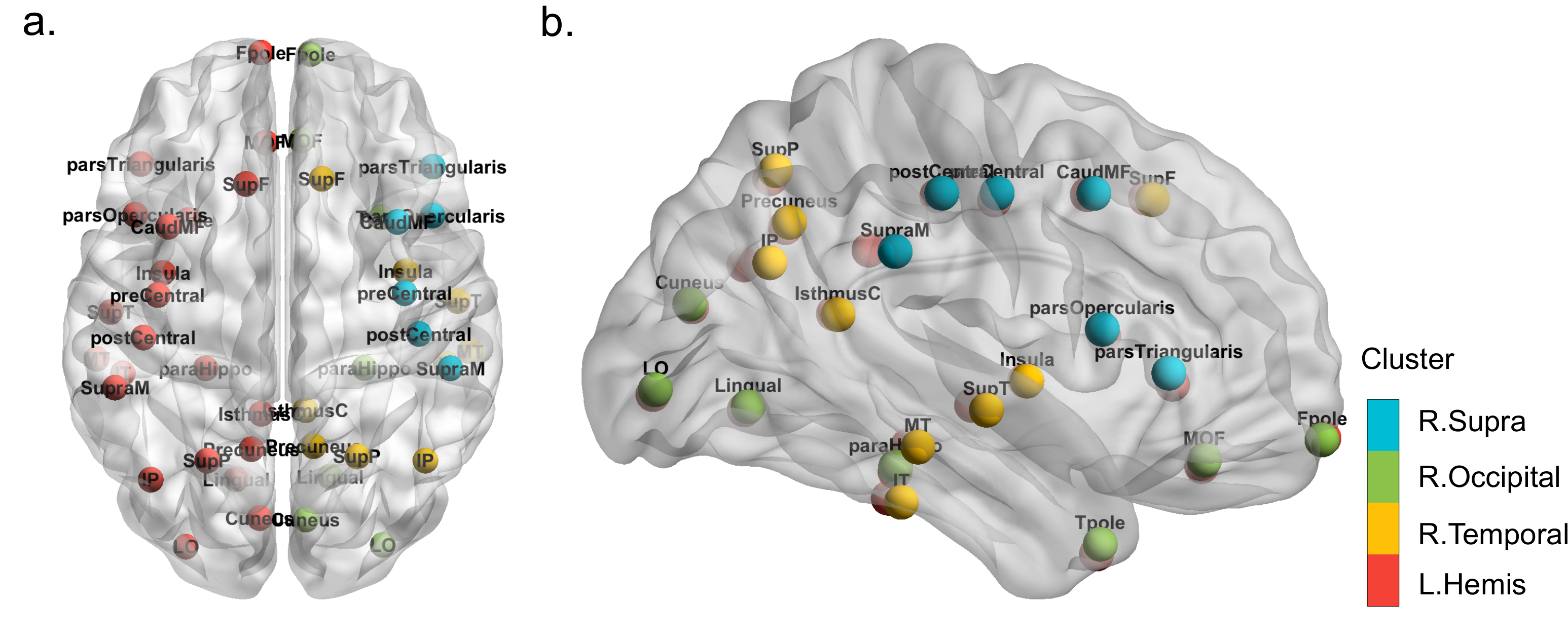}
    \caption{Illustration of brain node clustering results for HCP data with (a) top and (b) side views. }
    \label{fig:cluster_brain}
\end{figure*}

We apply our general asymmetric algorithm to the HCP data with the numbers of clusters on three modes $r_1 = r_2 = 4$ and $r_3 = 3$. The selection of $r_1$ and $r_2$ follows the human brain anatomy and the symmetry in the brain network, and the $r_3$ is specified following previous analysis~\citep{hu2022generalized}. Because of the symmetry in the data, the estimated brain node clustering results are the same on the first and second modes. Figure~\ref{fig:cluster_brain} shows that brain connection exhibits a strong spatial separation structure. Specifically, the first cluster, named \emph{L.Hemis}, involves all the nodes in the left hemisphere. The nodes in the right hemisphere are further separated into three clusters led by the middle-part tissues in Temporal and Parietal lobes (\emph{R.Temporal}), the back-part tissues in Occipital lobe (\emph{R.Occipital}), and the front-part tissues in Frontal and Parietal lobes (\emph{R.Supra}). This clustering result is reasonable since the left and right hemispheres often play different roles in human brains. 

\begin{figure*}[htb]
    \centering
    \includegraphics[width = 1\textwidth]{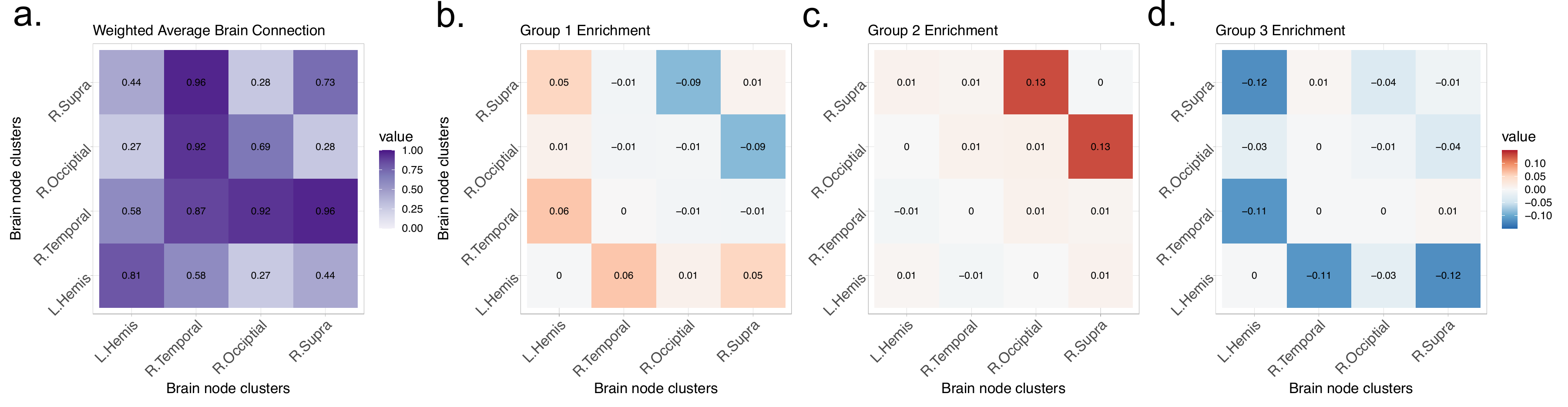}
    \caption{Mode 3 slices of estimated core tensor $\hat \tS$. (a) Average estimated slice weighted by the group size; (b)-(d) Group-specified enrichment, i.e., the difference between each slice of $\hat \tS$ and the averaged slice. }
    \label{fig:ests}
\end{figure*}

\begin{figure*}[htb]
    \centering
    \includegraphics[width = 1\textwidth]{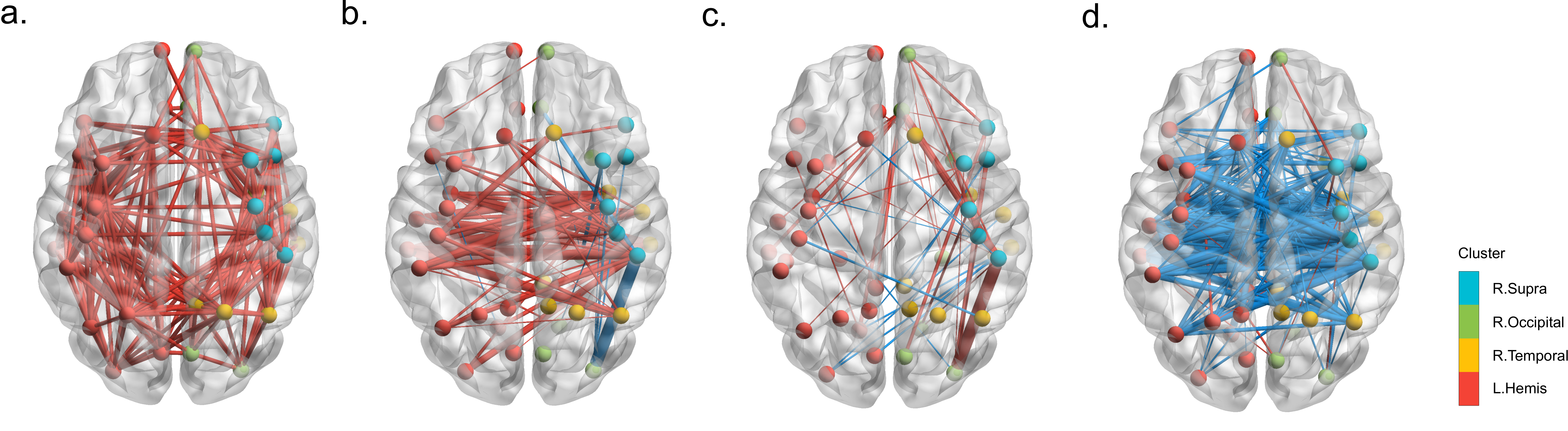}
    \caption{Observed brain connections in the population and each group of individuals. (a) Average brain network; (b)-(d) Group-specified brain network enrichments in Groups 1-3. Red edges represent the positive enrichment and blue edges represent the negative enrichment.}
    \label{fig:brain_conn}
\end{figure*}

Figure~\ref{fig:ests} illustrates the estimated core tensor $\hat \tS$ with estimated clustering, and Figure~\ref{fig:brain_conn} visualizes the average brain connections and the connection enrichment in contrast to average networks in each group. In general, we find that the inner-hemisphere connection has stronger connection compared to inter-hemisphere connections (Figure~\ref{fig:ests}a). Also, the back and front parts (\emph{R.Occipital}, \emph{R.Supra}) are shown to have more interactions with temporal tissues than inner-cluster connections. In addition, the group 1 with 54\% females shows an enrichment on the inter-hemisphere connections (Figure~\ref{fig:ests}b), while group 4 with only 36\% females exhibits a reduction (Figure~\ref{fig:ests}d). This result agrees with previous findings in \cite{hu2022generalized}. The enrichment on the back-front connection is also recognized in group~3 (Figure~\ref{fig:ests}c). The interpretive patterns in our results demonstrate the usefulness of our clustering methods in the human brain connectome data application. 

\subsection{Peru Legislation data analysis}

We also apply our method to the legislation networks in the Congress of the Republic of Peru \citep{lee2017time}. Because of the frequent political power shifts in the Peruvian Congress during 2006-2011, we choose to focus on the data for the first half of 2006-2007 year. The dataset records the co-sponsorship of 116 legislators from top 5 parties and 802 bill proposals. We reconstruct legislation network as an order-3 binary tensor $\tY \in \{0,1\}^{116 \times 116 \times 116}$, where $\tY_{ijk} = 1$ if the legislators $(i,j,k)$ have sponsored the same bill, and $\tY_{ijk} = 0$ otherwise. The true party affiliations of legislators are provided and serve as the ground truth. We apply various higher-order clustering methods to $\tY$ with $r = 5$. Table~\ref{tab:peru} shows that our \textbf{\small dTBM} achieves the best performance compared to others. The second best method is the two-stage algorithm \textbf{\small HLloyd}, followed by the spectral methods \textbf{\small SCORE} and \textbf{\small HOSVD+}. This result is consistent with our simulations under strong signal and moderate degree heterogeneity. The comparison suggests that our method \textbf{\small dTBM} is more appealing in real-world applications.

\begin{table}[ht]
\renewcommand{\arraystretch}{1.3}
    \centering
    \begin{tabular}{c |c  c cc c} 
    \hline 
        Method & \textbf{\small dTBM} 
        &\textbf{\small HOSVD}
        &\textbf{\small HOSVD+} & \textbf{\small HLloyd} &  \textbf{\small SCORE}\\
         CER & \textbf{0.116}
         &  0.22 
         &0.213 & 0.149 &0.199\\
         \hline
    \end{tabular}
    \caption{Clustering errors (measured by CER) for various methods in the analysis of Peru Legislation dataset.}
    \label{tab:peru}
\end{table}

\section{Proof Sketches}\label{sec:mainproof}

In this section, we provide the proof sketches for the main Theorem~\ref{thm:stats} (Impossibility), Theorem~\ref{thm:comp} (Impossibility), and Theorems~\ref{thm:initial}-\ref{thm:refinement}. Detail proofs and extra theoretical results are provided in Appendix B.

\subsection{Proof sketches of Theorems~\ref{thm:stats} and \ref{thm:comp} (Impossibility)} \label{subsec:impossible}

The proofs of impossibility in Theorems~\ref{thm:stats} and \ref{thm:comp} share the same proof idea with \citet[Theorems 6 and 7]{han2022exact} and \citet[Theorem 2]{gao2018community}. In both proofs of statistical and computational impossibilities, the key idea is to construct a particular set of parameters to lower bound the minimax rate. Specifically, for statistical impossibility in Theorem~\ref{thm:stats}, we construct a particular $(z^*_{\rm stats}, \mtheta^*_{\rm stats})  \in \tP_{z, \mtheta} $ such that for all $\tS^* \in \tP_{\tS}(\gamma)$
\begin{align}
    & \inf_{\hat z_{\rm stats}} \sup_{(z, \mtheta) \in \tP_{z, \mtheta}} \bbE[p \ell(\hat z_{\rm stat}, z)] \geq \inf_{\hat z_{\rm stats}} \bbE[p \ell(\hat z_{\rm stat}, z_{\rm stats}^*)| (z^*_{\rm stats}, \tS^*,\mtheta^*_{\rm stats})] \geq 1;\label{eq:stats_impo_sketch}  
\end{align}
for computational impossibility in Theorem~\ref{thm:comp}, we construct a particular $(z^*_{\rm comp}, \tS^*_{\rm comp}, \mtheta^*_{\rm comp}) \in \tP(\gamma)$ such that
\begin{align}
    &\inf_{\hat z_{\rm comp}} \sup_{(z, \tS, \mtheta) \in \tP(\gamma)} \bbE[p \ell(\hat z_{\rm comp}, z)] \geq \inf_{\hat z_{\rm comp}} \bbE[p \ell(\hat z_{\rm comp}, z_{\rm comp}^*)| (z^*_{\rm comp}, \tS^*_{\rm comp},\mtheta^*_{\rm comp})] \geq 1.
\end{align}

The constructions of $(z^*_{\rm stats}, \mtheta^*_{\rm stats})$ and $(z^*_{\rm comp}, \tS^*_{\rm comp}, \mtheta^*_{\rm comp})$ are the most critical steps. With good constructions, the lower bound ``$\geq 1$" can be verified by classical statistical conclusions (e.g.\ Neyman-Pearson Lemma) or prior work (e.g.\ HPC Conjecture). 

A notable detail in the proof of statistical impossibility is the arbitrariness of $\tS^*$. The first infimum over $ \tP_{\tS}(\gamma)$ in the minimax rate~\eqref{eq:minminmax} requires that the lower bound~\eqref{eq:stats_impo_sketch} holds for any $\tS^* \in \tP_{\tS}(\gamma)$. The arbitrary choice of $\tS^*$ brings extra difficulties in the parameter construction, and consequently a non-trivial $\mtheta^*_{\rm stats} \neq \mathbf{1}$ is chosen to address the arbitrariness. Previous TBM construction in the proof of \citet[Theorem 6]{han2022exact} with $\mtheta^*_{\rm stats} = \mathbf{1}$ is no longer applicable in our case. Meanwhile, our construction $(z^*_{\rm comp}, \tS^*_{\rm comp}, \mtheta^*_{\rm comp})$ leads to a rank-2 mean tensor to relate the HPC Conjecture while TBM \citet[Theorem 7]{han2022exact} constructs a rank-1 mean tensor. Hence, we emphasize that dTBM-specific techniques are required to obtain our impossibility results, though the proof idea is common for minimax lower bound analysis. 

\subsection{Proof sketch of Theorem~\ref{thm:initial}}
The proof of Theorem~\ref{thm:initial} is inspired by the proof idea of \citet[Lemma 1]{gao2018community}. The extra difficulties are the angle gap characterization and multilinear algebra property in tensors; we address both challenges in our proof. Specifically, we control the misclustering error by the estimation error of $\hat \tX$ calculated in Step 2 of Sub-algorithm~\hyperref[alg:main]{1}.  We prove the following inequality
\begin{align}
    \ell(z^{(0)},z) &\lesssim \frac{1}{p}\min_{\pi \in \Pi} \sum_{i: z^{(0)}(i) \neq \pi(z(i))} \theta(i)^2 \lesssim \frac{\sigma^2 r^{K-1}}{ \Delta_{\min}^2 p^K} \onormSize{}{\hat \tX - \tX}_F^2 \lesssim \frac{r^K p^{-K/2}}{\text{SNR}}, \label{eq:proof_4}
\end{align}
where $\tX = \bbE \tY$ is the true mean. The first inequality in~\eqref{eq:proof_4} holds with the assumption $\min_{i \in [p]} \theta(i) \geq c>0$ in Theorem~\ref{thm:initial}. 
{The second inequality relies on the key Lemma~\ref{lem:angle_gap_x}, which indicates}
\begin{equation}\label{eq:proof_4_gap}
\min_{z(i) \neq z(j)} \onormSize{}{[\mX_{i:}]^s - [\mX_{j:}]^s} \gtrsim \Delta_{\min},
\end{equation}
where $\mX = \mat(\tX)$. The most challenging part in the proof of Theorem~\ref{thm:initial} lies in the derivation of inequality~\eqref{eq:proof_4_gap} (or the proof of Lemma~\ref{lem:angle_gap_x}), in which the proof of \cite{gao2018community} is no longer applicable due to different angle gap assumption in our dTBM. To address the angle gap notion, we develop the extra padding technique in Lemma~\ref{lem:pad} and balance assumption~\eqref{eq:degree}. Last, we finish the proof of Theorem~\ref{thm:initial} by showing the third inequality of~\eqref{eq:proof_4} using \citet[Proposition 1]{han2022exact}. 

\subsection{Proof sketch of Theorem~\ref{thm:refinement}}\label{sec:thm5}
The proof of Theorem~\ref{thm:refinement} is inspired by the proof idea of \citet[Theorem 2]{han2022exact}. We develop extra polar-coordinate based techniques with angle gap characterization to address the nuisance degree heterogeneity. Recall the intermediate quantity, misclustering loss, defined in~\eqref{eq:defnofL}
\begin{align}
    L^{(t)} : =L(z,z^{(t)}) = \frac{1}{p}  \sum_{i \in [p]} \theta(i) \sum_{b \in [r]}  \ind \offf{ z^{(t)}(i) = b } \onorm{ \off{ \mS_{ z(i):}  }^s - \off{ \mS_{b:}  }^s  }^2.
\end{align}
\normalsize
We show that $L^{(t)}$ provides an upper bound for the misclustering error of interest via the inequality $\ell^{(t)}\leq {L^{(t)}\over \Delta^2_{\min}}$ in Lemma~\ref{lem:mis}. Therefore, it suffices to control $L^{(t)}$. Further, we introduce the oracle estimators for core tensor under the true cluster assignment via 
\begin{equation}
    \tilde \tS = \tY \times_1 \mW^T \times_2 \cdots \times_K \mW^T, 
\end{equation}
where $\mW = \mM \of{ \text{diag}(\mone_{p}^T \mM) }^{-1}$ is the weighted true membership matrix. Let $ \mV = \mW^{\otimes (K-1)}$ denote the Kronecker product of $(K-1)$ copies of $\mW$ matrices, and we define the $t$-th iteration quantities $\mW^{(t)}, \mV^{(t)}$ corresponding to $\mM^{(t)}$ (or equivalently $z^{(t)}$). To evaluate $L^{(t+1)}$, we prove the bound
\begin{align}
    &\ind \offf{ z^{(t+1)}(i) = b } = \ind \offf{       \onormSize{}{ [ \mY_{i:} \mV^{(t)}  ]^s - [\mS_{b:}^{(t)}]^s }^2 \leq \onormSize{}{ [ \mY_{i:} \mV^{(t)}  ]^s - [\mS_{z(i):}^{(t)}]^s }^2} \leq A_{ib} + B_{ib}, \label{eq:proof_5_event}
\end{align}
where $\mY = \mat(\tY)$, $ \mS = \mat(\tS)$, $\mS^{(t)} = \mat(\tS^{(t)})$ and
\begin{align}
        A_{ib} &= \ind \offf{\ang{ \mE_{i:} \mV, \off{  \tilde \mS_{z(i):} }^s - \off{  \tilde \mS_{b:} }^s } \lesssim -  \onorm{ \off{ \mS_{z(i):}  }^s - \off{ \mS_{b:}  }^s  }^2 },\\
        B_{ib} &= \ind \offf{\onorm{ \off{ \mS_{z(i):}  }^s - \off{ \mS_{b:}  }^s  }^2 \lesssim F_{ib}^{(t)} + G_{ib}^{(t)} + H_{ib}^{(t)} }.
\end{align}
\normalsize
The terms $F_{ib}^{(t)}, G_{ib}^{(t)}, H_{ib}^{(t)}$ are controlled by $z^{(t)}, \tS^{(t)}$; see the detailed definitions in \eqref{eq:f}, \eqref{eq:g}, \eqref{eq:h}. Note that the event $A_{ib}$ only involves the oracle estimator independent of $t$, while all the terms related to the $t$-th iteration are in $B_{ib}$. Thus, the inequality~\eqref{eq:proof_5_event} decomposes the misclustering loss in the $(t+1)$-th iteration into the oracle loss and the loss in $t$-th iteration. This decomposition leads to the separation of statistical error and computational error in the final upper bound of Theorem~\ref{thm:refinement}.

Specifically, we prove the contraction inequality
\begin{align}
     &L^{(t+1)} \leq M \xi + \rho L^{(t)}, 
     \text{ with } \xi = \frac{1}{p}  \sum_{i \in [p]} \theta(i) \sum_{b \in [r]}  A_{ib} \onorm{ \off{ \mS_{ z(i):}  }^s - \off{ \mS_{b:}  }^s  }^2, \label{eq:proof_5_ineq}
\end{align}
where $M$ is a positive constant, $\rho \in (0,1)$ is the contraction parameter, and we call $\xi$ the oracle loss. Controlling the probability of event $B_{ib}$ and obtaining the $\rho L^{(t)}$ term in the right hand side of~\eqref{eq:proof_5_ineq} are the most challenging parts in the proof of Theorem~\ref{thm:refinement}. Note that the true and estimated core tensors are involved via their normalized rows such as $\mS_{a:}^s, \tilde \mS_{a:}^s, [\mS^{(t)}_{a:}]^s$. The Cartesian coordinate based analysis in \cite{han2022exact} is no longer applicable in our case. Instead, we use the polar-coordinate based analysis and the geometry property of trigonometric functions to derive the high probability upper bounds for $F_{ib}^{(t)}, G_{ib}^{(t)}, H_{ib}^{(t)}$. 

Further, by sub-Gaussian concentration, we prove the high probability upper bound for oracle loss
\begin{equation}\label{eq:proof_5_xi}
    \xi  \lesssim {\text{SNR}^{-1}}\exp\of{- \frac{p^{K-1}\text{SNR}}{r^{K-1}}}.
\end{equation}
Combining the decomposition~\eqref{eq:proof_5_ineq} and the oracle bound~\eqref{eq:proof_5_xi}, we finish the proof of Theorem~\ref{thm:refinement}.

The proof of MLE error shares the similar idea as Theorems~\ref{thm:initial}-\ref{thm:refinement}. We first show a weaker polynomial rate for MLE and then improve the rate from polynomial to exponential through the iterations. The only difference is that the MLE remains the same over iterations due to its global optimality. See Appendix B, Section~\ref{sec:statprove2} for the detailed proof.

\section*{Acknowledgment}

This research is supported in part by NSF CAREER DMS-2141865, DMS-1915978, DMS-2023239, EF-2133740, and funding from the Wisconsin Alumni Research foundation. We thank Zheng Tracy Ke, Anru Zhang, Rungang Han, Yuetian Luo for helpful discussions and for sharing software packages.

\bibliographystyle{apalike}
\bibliography{tensor_wang}

\newpage

\appendix
\section*{Appendices}
\section*{A \ \  Additional numerical experiments}
\begin{figure}[htb]
    \centering
    \includegraphics[width = .85\columnwidth]{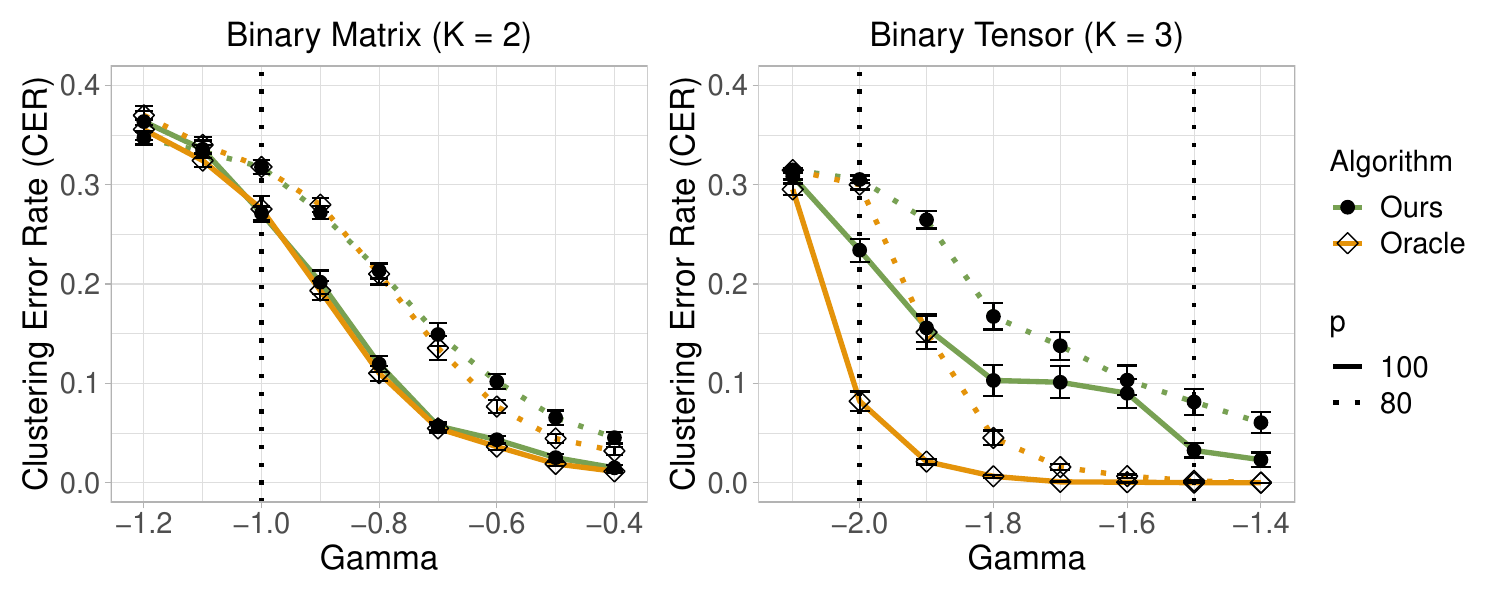}
    \caption{SNR phase transitions for Bernoulli dTBM with $p = \{80, 100\}, r = 5$ under (a) matrix case with $\gamma \in [-1.2, -0.4]$ and (b) tensor case with $ \gamma \in [-2.1, -1.4]$.}
    \label{fig:phase_binary}
\end{figure}

\setcounter{section}{1}
\textbf{Bernoulli phase transition.} The first additional experiment verifies the statistical-computational gap in Section~\ref{sec:limits} under the Bernoulli model. Consider the Bernoulli model with $p = \{80, 100\}$, $r = 5$. We vary $\gamma $ in $ [-1.2, -0.4]$ and $[-2.1, -1.4]$ for matrix ($K=2$) and tensor $(K = 3)$ clustering, respectively. We  approximate MLE using an oracle estimator, i.e., the output of Sub-algorithm~\hyperref[alg:main]{2} initialized from the true assignment. Figure~\ref{fig:phase_binary} shows a similar pattern as Figure~\ref{fig:phase}. The algorithm and oracle estimators have no gap in the matrix case, while an error gap emerges between the critical values $\gamma_{\text{stat}} = -2$ and $\gamma_{\text{comp}} = -1.5$ in the tensor case. Figure~\ref{fig:phase} suggests the statistical-computational gap in Bernoulli models.

{\bf Sparsity.} The second additional experiment evaluates the algorithm performances under the sparse binary dTBM~\eqref{eq:sparse_dtbm}. We fix the signal exponent $\gamma = -1.2$ and vary the sparsity parameter $\alpha_p \in [0.05, 0.9]$. A smaller $\alpha_p$ leads to a higher probability of zero entries in the observation. In addition to the three algorithms mentioned in Section~\ref{subsec:comp} (denoted {\bf \small Initialization}, {\bf \small dTBM}, and {\bf \small SCORE}), we consider other three algorithms based on the discussion in Section~\ref{subsec:ber}: 
  \begin{itemize} 
  \item \textbf{\small D-HOSVD}, the diagonal-deleted HOSVD in \cite{ke2019community}; 
  \item \textbf{\small D-HOSVD + Angle}, the combined algorithm of our angle-based iteration with initialization from \textbf{\small D-HOSVD};
  \item \textbf{\small SCORE + Angle}, the combined algorithms of our angle-based iteration with initialization from \textbf{\small SCORE}.
  \end{itemize}

\begin{figure}[htp!]
    \centering
    \includegraphics[width=.85\columnwidth]{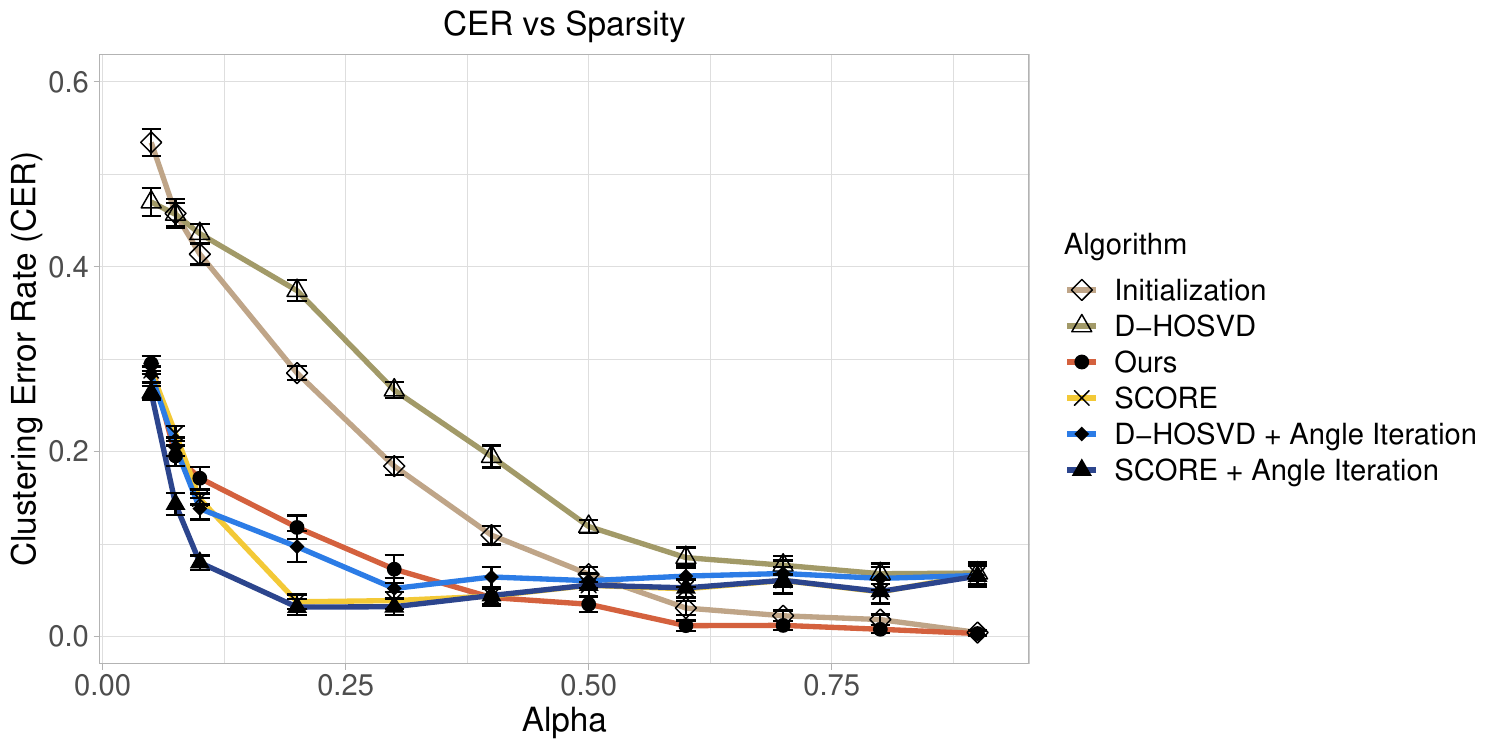}
    \caption{CER comparison versus sparsity parameter $\alpha_p$ in $[0.05, 0.9]$. We set $p = 100, r = 5$ and $\gamma = -1.2$ under sparse binary dTBM.}
    \label{fig:sparse}
\end{figure}

Figure~\ref{fig:sparse} shows a slightly larger error in \textbf{\small dTBM} than that in \textbf{\small SCORE}, \textbf{\small D-HOSVD + Angle}, and \textbf{\small SCORE + Angle} under the sparse setting with $\alpha_p < 0.3$. 
The small gap between  \textbf{\small dTBM} and other sparse-specific methods implies the robustness of our algorithm. In addition, comparing \textbf{\small SCORE} versus \textbf{\small SCORE + Angle} (or \textbf{\small D-HOSVD} versus \textbf{\small D-HOSVD + Angle}) indicates the benefit of our angle iterations under the sparse dTBM. In the intermediate and dense cases with $\alpha_p \geq 0.3$, our proposed \textbf{\small dTBM} has a clear improvement over others, which again verifies the success of our algorithm in dense settings.

\section*{B \ \  Proofs}
\setcounter{section}{2}

We provide the proofs for all the theorems in our main paper. In each sub-section, we first show the proof of main theorem and then collect the useful lemmas in the end.   {We combine the proofs of MLE achievement in Theorem~\ref{thm:stats} and polynomial-time achievement in Theorem~\ref{thm:refinement} in the last section due to the similar idea. }

\subsection{Notation}
Before the proofs, we first introduce the notation used throughout the appendix and the general dTBM without symmetric assumptions. The parameter space and minimal gap assumption are also extended for the general asymmetric dTBM.

{\bf Preliminaries.}
\begin{enumerate}[wide]
    \item For mode $ k \in [K]$, denote mode-$k$ tensor matricizations by
    \begin{align}
        &\mY_k = \mat_k \of{ \tY }, \quad \mS_k = \mat_k \of{\tS}, \quad \mE_k = \mat_k \of{ \tE}, \quad \mX_k = \mat_k \of{\tX}.
    \end{align}
    \item For a vector $\ma$, let $\ma^{s} \coloneqq \ma/\onorm{\ma}$ denote the normalized vector. We make the convention that $\ma^s = {\bf 0}$ if $\ma = {\bf 0}$. 
    \item For a matrix $\mA \in \bbR^{n \times m} $, let $\mA^{\otimes K}:=\mA\otimes \cdots \otimes \mA\in \bbR^{n^K \times m^K}$ denote the Kronecker product of $K$ copies of matrices $\mA $.
    \item For a matrix $\mA$, let $\onormSize{}{\mA}_\sigma$ denote the spectral norm of matrix $\mA$, which is equal to the maximal singular value of $\mA$; let $\lambda_k(\mA)$ denote the $k$-th largest singular value of $\mA$; let $\onormSize{}{\mA}_F$ denote the Frobenius norm of matrix $\mA$.
\end{enumerate}

{\bf Extension to general asymmetric dTBM.} 
 
 The general order-$K$ $(p_1, \ldots, p_K)$-dimensional dTBM with $r_k$ communities and degree heterogeneity $\mtheta_k = \entry{\theta_k(i)} \in \bbR_+^{p_k}$ is represented by
\begin{equation}\label{eq:general_dtbm}
    \tY = \tX+ \tE, \ \text{where}\ \tX=\tS \times_1 \mTheta_1 \mM_1 \times_2 \cdots \times_K \mTheta_K \mM_K,
\end{equation}
where $\tY \in \bbR^{p_1 \times \cdots \times p_K}$ is the data tensor, $\tX\in \mathbb{R}^{p_1\times \cdots \times p_K}$ is the mean tensor, $\tS \in \bbR^{r_1 \times \cdots \times r_K}$ is the core tensor, $\tE \in \bbR^{p_1 \times \cdots \times p_K}$ is the noise tensor consisting of independent zero-mean sub-Gaussian entries with variance bounded by $\sigma^2$, $\mTheta_k = \text{diag}(\mtheta_k)$, and $\mM_k\in \{0,1\}^{p_k \times r_k}$ is the membership matrix corresponding to the assignment $z_k: [p_k] \mapsto [r_k]$, for all $k \in [K]$. 

For ease of notation, we use $\{z_k\}$ to denote the collection $\{z_k\}_{k=1}^K$, and $\{\mtheta_k\}$ to denote the collection $\{\mtheta_k\}_{k=1}^K$. Correspondingly, we consider the parameter space for the triplet $\of{\{z_k\}, \tS, \{\mtheta_k\}}$,
\begin{align}
  \tP(\{r_k\}) = \Big\{ & \of{\{z_k\}, \tS, \{\mtheta_k\}}: \mtheta_k \in\mathbb{R}^p_{+}, {c_1 p_k\over r_k} |z_k^{-1}(a)| \leq {c_2 p_k\over r_k}, \\ & c_3 \leq  \onorm{\mS_{k,a:}} \leq c_4 ,\onormSize{}{\mtheta_{k,z_k^{-1}(a)}}_1=|z_k^{-1}(a)|, \text{for all }a \in [r_k], k\in[K]\Big\}. \label{eq:general_family}
\end{align}
\normalsize

We call the degree heterogeneity $\{\mtheta_k\}$ is balanced if for all $k \in [K]$,
\begin{equation}\label{eq:general_balanced}
    {\min_{a\in[r]} \onormSize{}{\mtheta_{k, z_k^{-1}(a)}}=\left(1+o(1)\right)\max_{a\in[r]}\onormSize{}{\mtheta_{k,z_k^{-1}(a)}}}.
\end{equation}

We also consider the generalized Assumption~\ref{assmp:min_gap} on angle gap.
\begin{assumption}[Generalized angle gap]\label{assmp:general_minimal_gap} Recall $\mS_k=\mat_k(\tS)$. We assume the minimal gap between normalized rows of $\mS_k$ is bounded away from zero for all $k\in[K]$; i.e.,
\begin{equation}
     \Delta_{\min} \coloneqq \min_{k \in [K]} \min_{a \neq b \in [r_k]} \onorm{ \mS_{k,a:}^s - \mS_{k, b:}^s } > 0.
\end{equation}
\end{assumption}
Similarly, let $\text{SNR} = \Delta_{\min}^2/\sigma^2$ with the generalized minimal gap $\Delta_{\min}^2$ defined in Assumption~\ref{assmp:general_minimal_gap}. We define the regime
\begin{align}
    \tP(\gamma) = \tP (\{r_k\}) \cap\{\tS \text{ satisfies $\text{SNR} = p^{\gamma}$ and $p_k \asymp p, k \in [K]$} \}.
\end{align}
\normalsize

\subsection{Proof of Theorem~\ref{thm:unique}}

\begin{proof}[Proof of Theorem~\ref{thm:unique}] 

To study the identifiability, we consider the noiseless model with $\tE = 0$. Assume that there exist two parameterizations satisfying
\begin{align}
    \tX&=\tS\times_1\Theta_1 \mM_1 \times_2 \cdots \times_K \Theta_K \mM'_K =\tS'\times_1\Theta'_1 \mM'_1 \times_2 \cdots \times_K \Theta'_K \mM'_K,\label{eq:another}
\end{align}
where $\of{ \{z_k\}, \tS, \{ \mtheta_k \} } \in \tP(\{r_k\})$ and $\of{ \{z'_k\}, \tS', \{ \mtheta'_k \} } \in \tP(\{r_k'\})$ are two sets of parameters. We prove the sufficient and necessary conditions separately.

\begin{enumerate}[wide]
    \item[$(\Leftarrow)$] For the necessity, it suffices to construct two distinct parameters up to cluster label permutation, if the model~\eqref{eq:general_dtbm} violates Assumption~\ref{assmp:general_minimal_gap}. {Note that $\Delta_{\min}^2 = 1$ when there exists $k \in [K]$ such that $r_k = 1$. Hence, we consider the case that  $r_k \geq 2$ for all $k \in [K]$.}  Without loss of generality, we assume $\onorm{ \mS_{1,1:}^s - \mS_{1,2:}^s } = 0$.

{
  
By constraints in parameter space~\eqref{eq:general_family}, neither $\mS_{1,1:}$ nor $\mS_{1,2:}$ is a zero vector.
}
There exists a positive constant $c$ such that $\mS_{1,1:} = c \mS_{1,2:}$. Thus, there exists a core tensor $\tS_0 \in \bbR^{r_1 -1 \times \cdots \times r_K}$ such that 
\begin{equation}
    \tS = \tS_0 \times_1 \mC \mR,
\end{equation}
where $\mC = \text{diag}(1, c, 1,...,1) \in \bbR^{r_1 \times r_1}$ and 
\begin{equation}
    \mR = \begin{pmatrix}
    1& 0\\
    1&0 \\
    0 & \mone_{r_1-2}
    \end{pmatrix} \in \bbR^{r_1 \times (r_1 -1)}.
\end{equation}
Let $\mD = \text{diag}(1+c, 1,...,1) \in \bbR^{r_1 -1 \times r_1 -1}$. Consider the parameterization $\mM'_1 = \mM_1 \mR,\tS' = \tS_0 \times_1 \mD,$ and 
\begin{equation}
     \theta'_{1}(i) = \begin{cases}
    \frac{1}{1+c} \theta_{1}(i) & i \in z_{1}^{-1}(1),\\
     \frac{c}{1+c} \theta_{1}(i) & i \in z_{1}^{-1}(2),\\
     \theta_{1}(i) & \text{ otherwise},
    \end{cases}
\end{equation}
and $\mM'_k = \mM_k, \mtheta'_k = \mtheta_k$ for all $k = 2, \ldots, K$. Then we have constructed a
triplet $\of{ \{z'_k\}, \tS', \{ \mtheta'_k \} }$ that is distinct from $\of{ \{z_k\}, \tS, \{ \mtheta_k \} }$ up to label permutation. 

\item[$(\Rightarrow)$] For the sufficiency, it suffices to show that all possible triplets $\of{ \{z'_k\}, \tS', \{ \mtheta'_k \} }$ are identical to $\of{ \{z_k\}, \tS, \{ \mtheta_k \} }$ up to label permutation if the model~\eqref{eq:general_dtbm} satisfies Assumption~\eqref{assmp:general_minimal_gap}. We show the uniqueness of the three parameters, $\{\mM_k\}, \{\tS\}, \{\mtheta_k\}$ separately.

First, we show the uniqueness of $\mM_k$ for all $k \in [K]$. 
{When $r_k = 1$, all possible $\mM_k$'s are equal to the vector ${\bf 1}_{p_k}$, and the uniqueness holds trivially. Hence, we consider the case that $r_k \geq 2$. Without loss of generality, we consider $k=1$ with $r_1 \geq 2$ and show the uniqueness of the first mode membership matrix;}  
i.e., $\mM'_1 = \mM_1\mP_1$ where $\mP_1$ is a permutation matrix. The conclusion for $k\geq 2$ can be showed similarly and thus omitted. 

Consider an arbitrary node pair $(i,j)$. If $z_1(i) = z_1(j)$, then we have $\onormSize{}{\mX_{1, z_1(i):}^s - \mX_{1, z_1(j):}^s } = 0$ and thus $\onormSize{}{ (\mS')_{1, z'_1(i):}^{s} - (\mS')_{1, z'_1(j):}^{s} } = 0$ by Lemma~\ref{lem:angle}. Then, by Assumption~\eqref{assmp:general_minimal_gap}, we have $z'_1(i) = z'_1(j)$. Conversely, if $z_1(i) \neq z_1(j)$, then we have $ \onorm{\mX_{1,i:}^s - \mX_{1,j:}^s} \neq 0$ and thus $\onorm{ (\mS')_{1, z'_1(i):}^{s} - (\mS')_{1, z'_1(j):}^{s} } \neq 0$ by Lemma~\ref{lem:angle}. Hence, we have $z'_1(i) \neq z'_1(j)$. Therefore, we have proven that $z'_1$ is identical $z_i$ up to label permutation.

Next, we show the uniqueness of $\mtheta_k$ for all $k \in [K]$ provided that $z_k = z_k'$. Similarly, consider $k=1$ only, and omit the procedure for $k\geq 2$. 

Consider an arbitrary $j \in [p_1]$ such that $z_1(j) = a$. Then for all the nodes $i \in  z_1^{-1}(a)$ in the same cluster of $j$, we have 
\begin{equation}
    \frac{\mX_{1,z_1(i):}}{\mX_{1,z_1(j):}} = \frac{\mX'_{1,z_1(i):}}{\mX'_{1,z_1(j):}}, \text{ which implies }  \frac{\theta_1(j)}{\theta_1(i)} = \frac{\theta'_1(j)}{\theta'_1(i)}.\label{eq:theta_uniq}
\end{equation}
Let $\theta'_1(j) = c\theta_1(j)$ for some positive constant $c$. By equation~\eqref{eq:theta_uniq}, we have $\theta'_1(i) = c \theta_1(i)$ for all $ i \in  z_1^{-1}(a)$. By the constraint $(\{z_k\}, \tS', \{\mtheta'_k\}) \in \tP(\{r_k\})$, we have 
\begin{equation}
    \sum_{j \in z_1^{-1}(a)} \theta'_1(j) = c \sum_{j \in z_1^{-1}(a)} \theta_1(j) = 1,
\end{equation}
which implies $c = 1$. Hence, we have proven $\mtheta_1 = \mtheta'_1$ provided that $z_1 = z'_1$.

Last, we show the uniqueness of $\tS$; i.e., $\tS'=\tS\times_1 \mP^{-1}_1\times_2\cdots \times_K \mP^{-1}_K$, where $\mP_k$'s are permutation matrices for all $k\in[K]$.  Provided $z'_k = z_k, \mtheta'_k = \mtheta_k$, we have $\mM'_k = \mM_k \mP_k$ and $\mTheta'_k = \mTheta_k$ for all $k \in [K]$. 

Let $\mD_k = \off{ (\mTheta'_k \mM'_k)^T (\mTheta'_k \mM'_k) }^{-1} (\mTheta'_k \mM'_k)^T, k \in [K]$. By the parameterization~\eqref{eq:another}, we have 
\begin{align}
    \tS' &= \tX \times_1 \mD_1 \times_2 \cdots \times_K \mD_K \\
    &= \tS \times_1 \mD_1 \mTheta_1 \mM_1 \times_1 \cdots \times_K \mD_K \mTheta_K \mM_K \\
    &= \tS \times_1 \mP^{-1}_1 \times_2 \cdots \times_K \mP^{-1}_K.
\end{align}

\end{enumerate}

Therefore, we finish the proof of Theorem~\ref{thm:unique}.
\end{proof}

{\bf Useful Lemma for the Proof of Theorem~\ref{thm:unique}} 

\begin{lem}[Motivation of angle-based clustering]\label{lem:angle} Consider the signal tensor $\tX$ in the general asymmetric dTBM~\eqref{eq:general_dtbm} with $(\{z_k\},\tS,\{\mtheta_k\})\in \tP(\{r_k\})$ and $r_k \geq 2, k \in [K]$. Then, for any $k \in [K]$ and index pair $(i,j)\in[p_k]^2$, we have 
\begin{align}
     & \onorm{ \mS_{k,z_k(i):}^s -  \mS_{k,z_k(j):}^s } = 0 \quad  \text{if and only if} \quad  \onorm{  \mX_{k, z_k(i):}^s -  \mX_{k,z_k(j):}^s } = 0.
\end{align}
\end{lem}

\begin{proof}[Proof of Lemma~\ref{lem:angle}] Without loss of generality, we prove $k = 1$ only and drop the subscript $k$ in $\mX_k, \mS_k$ for notational convenience. 
By tensor matricization, we have
\begin{equation}
    \mX_{j:} = \theta_1(j) \mS_{z_1(j):} \off{\mTheta_2 \mM_2 \otimes \cdots \otimes \mTheta_K \mM_K}^T.
\end{equation}     
Let $\tilde \mM = \mTheta_2 \mM_2 \otimes \cdots \otimes \mTheta_K \mM_K$. Notice that for two vectors $\ma, \mb$ and two positive constants $c_1, c_2 >0$, we have
\begin{equation}
\onorm{\ma^s - \mb^s} = \onorm{(c_1 \ma)^s - (c_2\mb)^s}.
\end{equation}
Thus it suffices to show the following statement holds for any index pair $(i,j)\in[p_1]^2$,
\begin{align}
    &\onorm{ \mS_{z_1(i):}^s - \mS_{z_1(j):}^s} = 0 \quad \text{if and only if}  \quad \onorm{ \off{\mS_{z_1(i):} \tilde \mM^T }^s - \off{\mS_{z_1(j):}\tilde \mM^T}^s} = 0.
\end{align}
\begin{enumerate}[wide]
    \item[$(\Leftarrow)$] Suppose $\onorm{ \off{\mS_{z_1(i):} \tilde \mM^T }^s - \off{\mS_{z_1(j):}\tilde \mM^T}^s} = 0$. There exists a positive constant $c$ such that $\mS_{z_1(i):} \tilde \mM^T= c \mS_{z_1(j):} \tilde \mM^T$. Note that
\begin{equation}
    \mS_{z_1(i):} = \mS_{z_1(i):} \tilde \mM^T \off{ \tilde \mM \of{ \tilde \mM^T  \tilde \mM}^{-1}},
\end{equation}
where $ \tilde \mM^T  \tilde \mM$ is an invertiable diagonal matrix with positive diagonal elements. Thus, we have $ \mS_{z_1(i):} = c  \mS_{z_1(j):}$, which implies $ \onorm{  \mS_{z_1(i):}^s -  \mS_{z_1(j):}^s } = 0 $.

\item[$(\Rightarrow)$] Suppose $ \onorm{ \mS_{z_1(i):}^s - \mS_{z_1(j):}^s } = 0 $. There exists a positive constant $c$ such that $\mS_{z_1(i):} = c \mS_{z_1(j):}$, and thus $\mS_{z_1(i):} \tilde \mM^T = c \mS_{z_1(j):} \tilde \mM^T$, which implies $\onorm{\left[\mS_{z_1(i):} \tilde \mM^T\right]^s- \left[\mS_{z_1(j):} \tilde \mM^T\right]^s}=0$.
\end{enumerate}
Therefore, we finish the proof of Lemma~\ref{lem:angle}.
\end{proof}

\subsection{Proof of Lemma~\ref{lem:angle_gap_x} and Lemma~\ref{lem:mis}}

\begin{proof}[Proof of Lemma~\ref{lem:angle_gap_x}] 

Note that the vector $\mS_{z(i):}$ can be folded to a tensor $\tS' =\entry{\tS'_{a_2,\ldots,a_K}} \in \bbR^{r^{K-1}}$; i.e., $\text{vec}(\tS') = \mS_{z(i):}$. Define weight vectors $\mw_{a_2, \cdots, a_K}$ corresponding to the elements in $\tS'_{a_2,\ldots,a_K}$ by
\begin{equation}
    \mw_{a_2 \cdots a_K} = [ \mtheta_{z^{-1} (a_2)}^T \otimes \cdots \otimes \mtheta_{z^{-1} (a_K)}^T] \in \bbR^{|z^{-1} (a_2)| \times \cdots \times |z^{-1} (a_K)|},
\end{equation}
for all $a_k \in [r], k = 2,\ldots, K$, where $\otimes$ denotes the Kronecker product. Therefore, we have  $\mX_{i:} = \theta(i) \pad_{\mw}(\mS_{z(i):})$ where $\mw = \{ \mw_{a_2, \cdots, a_K}\}_{a_k \in [r], k\in [K]/\{1\}}$. Specifically, we have $\onormSize{}{\mw_{a_2, \ldots, a_K}}^2 = \prod_{k = 2}^K \onormSize{}{ \mtheta_{z^{-1}(a_k)}}^2$, and by the balanced assumption~\eqref{eq:degree} we have
\begin{equation}\label{eq:pad_balance}
    \max_{(a_2,\ldots, a_K)} \onormSize{}{\mw_{a_2,\ldots, a_K}}^2 = (1 + o(1))  \min_{(a_2,\ldots, a_K)} \onormSize{}{\mw_{a_2,\ldots, a_K}}^2.
\end{equation}
\normalsize

Consider the inner product of $\mX_{i:}$ and $\mX_{j:}$ for $z(i) \neq z(j)$. By the definition of weighted padding operator~\eqref{eq:paddef} and the balanced assumption~\eqref{eq:pad_balance}, we have 
\begin{align}
    \ang{\mX_{i:},\mX_{j:}}
    &\ = \theta(i) \theta(j) \ang{ \pad_{\mw}(\mS_{z(i):}), \pad_{\mw}(\mS_{z(j):} )}\\
    &\ = \theta(i) \theta(j) \min_{(a_2,\ldots, a_K)} \onormSize{}{\mw_{a_2,\ldots, a_K}}^2 \ang{\mS_{z(i):}, \mS_{z(j):}}(1 + o(1)).
\end{align}
Therefore, when $p$ large enough, the inner product $\ang{\mX_{i:},\mX_{j:}} $ has the same sign as $\ang{\mS_{z(i):}, \mS_{z(j):}}$. 

{
 Then, we have 
\begin{align}
    \cos (\mS_{z_1(i):}, \mS_{z_1(j):}) &= \frac{ \ang{  \mS_{z_1(i):}, \mS_{z_1(j):} }}{ \onormSize{}{\mS_{z_1(i):}} \onormSize{}{\mS_{z_1(j):}} }= (1+o(1))\frac{ \ang{  \mX_{i:}, \mX_{j:} }}{ \onormSize{}{ \mX_{i:}} \onormSize{}{ \mX_{j:}} }= (1 + o(1)) \cos(\mX_{i:}, \mX_{j:}),
\end{align}
where the second inequality follows by the balance assumption on $\mtheta$.

Further, notice that $\onormSize{}{\mv_1^s - \mv_2^s}^2 =  2(1 - \cos(\mv_1, \mv_2))$. For all $i,j$ such that $z(i) \neq z(j)$, when $p \rightarrow \infty$, we have
\begin{equation}
    \onormSize{}{\mX_{i:}^s - \mX_{j:}^s} \asymp \onormSize{}{\mS_{z_1(i):}^s - \mS_{z_1(j):}^s} \gtrsim  \Delta_{\min}.
\end{equation}

}



\end{proof}

\begin{proof}[Proof of Lemma~\ref{lem:mis}]
By the definition of minimal gap in Assumption~\ref{assmp:min_gap}, we have 
\begin{align}
     L^{(t)} &= \frac{1}{p}  \sum_{i \in [p]} \theta(i) \sum_{b \in [r]}  \ind \offf{ z^{(t)}(i) = b } \onormSize{}{ [ \mS_{ z(i):}  ]^s - [ \mS_{b:}  ]^s  }^2 \\
     &\geq \frac{1}{p}  \sum_{i \in [p]} \theta(i) \sum_{b \in [r]}  \ind \offf{ z^{(t)}(i) = b } \Delta_{\min}^2 \\
     & \geq c \ell^{(t)} \Delta_{\min}^2,
\end{align}
    where the last inequality follows from the assumption $\min_{i \in [p]} \theta(i) \geq c>0$.
\end{proof}

\subsection{Proof of Theorem~\ref{thm:stats} (Impossibility)}\label{sec:statprove1}

\begin{proof}[Proof of Theorem~\ref{thm:stats} (Impossibility)]Consider the general asymmetric dTBM~\eqref{eq:general_dtbm} in the special case that $p_k = p$ and $r_k = r$ for all $ k\in [K]$ {with $K\geq 2$, $2 \leq r\lesssim p^{1/3}$ as $p \rightarrow \infty$}. For simplicity, we show the minimax rate for the estimation on the first mode $\hat z_1$; the proof for other modes are essentially the same. 
   
   To prove the minimax rate~\eqref{eq:minminmax}, it suffices to take an arbitrary $\tS^* \in  \tP_{\tS}(\gamma)$ wih $\gamma < -(K-1)$ and construct $(z^*_k, \mtheta^*_k)$ such that 
   \begin{equation}
       \inf_{\hat z_1} \bbE \left[ p\ell(\hat z_1, z_1^*) | (z^*_k, \tS^*,  \mtheta^*_k)  \right]\geq 1.
   \end{equation}
   
   We first define a subset of indices $T_k \subset [p_k], k \in [K]$ in order to avoid the complication of label permutation. Based on \citet[Proof of Theorem 6]{han2022exact}, we consider the restricted family of $\hat z_k$'s for which the following three conditions are satisfied:
   \begin{align}
        &\text{(a)}\ \hat z_k(i)=z_k(i) \text{ for all }i\in T_k; \quad \text{(b)} \ |T^c_k|\asymp {p\over r}; \\
        &\text{(c)}\ \min_{\pi\in \Pi}\sum_{i\in[p]}\ind\{\hat z_k(i) \neq \pi\circ z_k (i)\} = \sum_{i\in[p]}\ind\{\hat z_k(i) \neq  z_k (i)\},
   \end{align}
for all $k \in [K]$.
   Now, we consider the construction:
   \begin{enumerate}
       \item[(i)] $\{z_k^*\}$ satisfies properties (a)-(c) with misclassification sets $T_k^c$ for all $k \in [K]$;
       \item [(ii)] $\{\mtheta_k^*\}$ such that $\mtheta_k^*(i) \leq \sigma r^{(K-1)/2} p^{-(K-1)/2}$ for all $i \in T_k^c, k \in [K]$ and $\max_{k \in [K], a \in [r]} \onormSize{}{\mtheta_{k, z^{*, -1}_k(a)}}^2_2$ $ \asymp p/r$.
   \end{enumerate}
   
   Combining the inequalities (39) and (40) in the proof of Theorem 2 in \cite{gao2018community}, we have 
   \begin{align}
         \inf_{\hat z_1 } &\bbE \left[ \ell(\hat z_1, z_1^*) | (z^*_k, \tS^*, \mtheta^*_k)  \right] \geq\\
         & \frac{C}{r^3 |T_1^c|} \sum_{i \in T_1^c} \inf_{\hat z_1  (i)} \{ \bbP[\hat z_1(i) = 1| z_1^*(i) = 2,z^*_k, \tS^*, \mtheta^*_k] +   \bbP[\hat z_1(i) = 2| z_1^*(i) = 1, z^*_k, \tS^*, \mtheta^*_k] \},\label{eq:inf_lower}
   \end{align}
   where $C$ is some positive constant,  $\hat z_1$ on the left hand side denote the generic assignment functions in $\tP(\gamma)$, and the infimum on the right hand side is taken over the generic assignment function family of $\hat z_1(i)$ for all nodes $i \in T_1^c$. Here, the factor $r^3=r\cdot r^2$ in~\eqref{eq:inf_lower} comes from two sources: $r^2\asymp {r\choose 2}$ comes from the multiple testing burden for all pairwise comparisons among $r$ clusters; and another $r$ comes from the number of elements $|T^c_k|\asymp p/r$ to be clustered. 
   
   Next, we need to find the lower bound of the rightmost side in~\eqref{eq:inf_lower}. 
We consider the hypothesis test based on model~\eqref{eq:general_dtbm}. First, we reparameterize the model under the construction (i)-(ii).
\begin{equation}
    \mx_a^* = \off{\mat_1 \of{ \tS^*\times_2 \mTheta_2^* \mM_2^*\times_3\cdots \times_K \mTheta^*_K \mM_K^* }}_{a:},
\end{equation}
for all $a \in [r]$, where $\mx_a^*$'s are centroids in $\bbR^{p^{K-1}}$. Without loss of generality, we consider the lower bound for the summand in~\eqref{eq:inf_lower} for $i=1$. The analysis for other $i\in T^c_1$ are similar. For notational simplicity, we suppress the subscript $i$ and write $\my, \theta^*, z$ in place of $\my_1, \mtheta_1^*(1)$ and $z_1(1)$, respectively. The equivalent vector problem for assessing the summand in~\eqref{eq:inf_lower} is
\begin{equation}\label{eq:z}
\my=\theta^* \mx_{z}^*+\me,
\end{equation}
where $z\in \{1,2\}$ is an unknown parameter, $\theta^* \in \bbR_+$ is the given heterogeneity degree, $\mx_1^*,\mx_2^*\in\bbR^{p^{K-1}}$ are given centroids, and $\me\in\bbR^{p^{K-1}}$ consists of i.i.d.\ $N(0,\sigma^2)$ entries.  Then, we consider the hypothesis testing under the model~\eqref{eq:z}:
\begin{equation}\label{eq:test}
 H_0: z = 1, \my = \theta^* \mx_1^* + \me \  \leftrightarrow \  H_1: z = 2, \my = \theta^* \mx_2^* + \me,
\end{equation}
\normalsize
   
   The hypothesis testing~\eqref{eq:test} is a simple versus simple testing, since the assignment $z$ is the only unknown parameter in the test. 
   By Neyman-Pearson lemma, the likelihood ratio test is optimal with minimal Type I + II error. Under Gaussian model, the likelihood ratio test of \eqref{eq:test} is equivalent to the least square estimator $\hat z_{LS} = \argmin_{a = \{1,2\}} \onormSize{}{\my - \theta^* \mx_a^*}_F^2$. 
   
   Let $\mS = \mat_1(\tS)$. Note that 
   \begin{align}
       \onormSize{}{ \theta^* \mx_1^*  - \theta^* \mx_2^*}_F &\leq  \theta^* \onormSize{}{ \mS^*_{1:} - \mS^*_{2:} }_F \prod_{k = 2}^K \lambda_{\max}(\mTheta_k^*\mM_k^*)  \\
       & \leq \theta^* \onormSize{}{ \mS^*_{1:} - \mS^*_{2:} }_F   \max_{k \in [K]/\{1\}, a \in [r]} \onormSize{}{\mtheta_{k, z^{*, -1}_k(a)}}_2^{K-1} \\
       & \leq  \sigma r^{(K-1)/2} p^{-(K-1)/2} 2 c_4 p^{(K-1)/2} r^{-(K-1)/2}\\
       & \leq 2 c_4 \sigma, \label{eq:bound}
   \end{align}
   where $\lambda_{\max}(\cdot)$ denotes the maximal singular value, the second inequality follows from Lemma~\ref{lem:singular_thetam}, and the third inequality follows from property (ii) and the boundedness constraint in $\tP_{\tS}(\gamma)$ such that $\onormSize{}{ \mS^*_{1:} - \mS^*_{2:} }_F  \leq \onormSize{}{\mS_{1:}^*}_F + \onormSize{}{\mS_{2:}^*}_F \leq 2c_4$.

   Hence, we have 
   \begin{align}
       & \inf_{\hat z_1(1)} \{ \bbP[\hat z_1(1) = 1| z_1^*(1) = 2, z^*_k, \tS^*, \mtheta^*_k]  +   \bbP[\hat z_1(1) = 2| z_1^*(1) = 1, z^*_k, \tS^*, \mtheta^*_k] \} \\
       & \quad = 2 \bbP[ \hat z_{LS} = 1 |  z_1^*(1) = 2, z^*_k, \tS^*, \mtheta^*_k ] \\
      & \quad  = 2 \bbP[ \onormSize{}{\my - \theta^* \mx_1^*}_F^2 \leq \onormSize{}{\my - \theta^* \mx_2^*}_F^2 | z_1^*(1) = 2, z^*_k, \tS^*, \mtheta^*_k   ]\\
      & \quad = 2 \bbP [ 2 \langle \me, \theta^* \mx_1^* - \theta^*\mx_2^* \rangle \geq  \onormSize{}{\theta^* \mx_1^* - \theta^*\mx_2^*}_F^2 ]\\
       & \quad =  2 \bbP[ N(0,1) \geq \theta^* \onormSize{}{\mx_1^* - \mx_2^*}_F /(2\sigma) ] \\
       & \quad \geq 2 \bbP[ N(0,1) \geq c_4 ]   \geq c, \label{eq:inf}
   \end{align}
where the first equation holds by symmetry, the third equation holds by rearrangement, the fourth equation holds from the fact that $\ang{\me, \theta^* \mx_1^* - \theta^*\mx_2^* } \sim N(0, \sigma \onormSize{}{\theta^* \mx_1^* - \theta^*\mx_2^*}_F)$, and $c$ is some positive constant in the last inequality.
   
   Plugging the inequality~\eqref{eq:inf} into the inequality~\eqref{eq:inf_lower} for all $i \in T_1^c$, then, we have 
   \begin{equation}
       \liminf_{p \rightarrow \infty}  \inf_{\hat z_1 } \bbE \left[ p \ell(\hat z_1, z_1^*) | z^*_k, \mtheta^*_k, \tS^*  \right]  \geq \liminf_{p \rightarrow \infty} \frac{ C cp }{r^3} \geq C c,
   \end{equation}
   where the last inequality follows by the condition $r = o(p^{1/3})$. By the discrete nature of the misclustering error, we obtain our conclusion
   \begin{equation}
       \liminf_{p \rightarrow \infty} \inf_{\tS^* \in  \tP_{\tS}(\gamma)}  \inf_{\hat z_{\text{stat}} }\sup_{ (z^*, \mtheta^*) \in \tP_{z, \mtheta}} \bbE \left[ p\ell(\hat z_{\text{stat}}, z) \right]  \geq 1. 
   \end{equation}

{

Last, with constructed $z^*_k, \mtheta^*_k$ satisfying properties (i) and (ii) and $\gamma' < -(K-1)$, we construct a core tensor $\tS^*$ such that $\Delta_{\mX^*}^2 \leq p^{-(K-1)}$. Based on the property (ii) and the boundedness constraint of $\tS^*$ in $\tP$, we still have $\onormSize{}{ \theta^* \mx_1^*  - \theta^* \mx_2^*}_F \leq 2c_4 \sigma$. Hence, we obtain the desired result 
\begin{align}
    &\liminf_{p \rightarrow \infty} \inf_{\hat z_1} \sup_{(z, \tS, \mtheta) \in \tP'(\gamma')} \bbE \left[ p\ell(\hat z_1, z_1) \right]  \geq   \liminf_{p \rightarrow \infty} \inf_{\hat z_{\rm stat}}  \bbE \left[ p\ell(\hat z_1, z^*_1) |  z^*_k, \tS^*, \mtheta^*_k  \right]  \geq 1.
\end{align}


}

   \end{proof}

\subsection{Proof of Theorem~\ref{thm:comp} (Impossibility)}\label{sec:compprove1}

\begin{proof}[Proof of Theorem~\ref{thm:comp} (Impossibility)]
The idea of proving computational hardness is to show the computational lower bound for a special class of degree-corrected tensor clustering model with $K\geq 2$ {and $ r \geq 2$}. We construct the following special class of higher-order degree-corrected tensor clustering  model. For a given signal level $\gamma\in\bbR$ and noise variance $\sigma$, define a rank-2 symmetric tensor $\tS\in\bbR^{3\times \cdots\times 3}$ subject to
\begin{equation}\label{eq:S}
\tS=\tS(\gamma) = \begin{bmatrix}
		1\\
		1\\
		1
		\end{bmatrix}^{\otimes K}
+\sigma p^{-\gamma/2} \begin{bmatrix}
		1\\
		-1\\
		0
	\end{bmatrix}^{\otimes K}.
\end{equation}
Then, we consider the signal tensor family
\begin{align}
\tP_{\text{shifted}}(\gamma)=\{\tX\colon & \tX=\tS\times_1\mM_1\times_2 \cdots \times_K\mM_K,\  \text{$\mM_k\in\{0,1\}^{p\times 3}$ is a membership matrix that}
\\
& \text{satisfies $|\mM_k(\colon,i)|\asymp p$ for all $i\in[3]$ and $k\in[K]$}\}.
\end{align}
\normalsize
We claim that the constructed family satisfies the following two properties:
\begin{enumerate}[wide,label=(\roman*)]
    \item For every $\gamma\in \mathbb{R}$, $\tP_{\text{shifted}}(\gamma)\subset \tP(\gamma)$, where $\tP(\gamma)$ is the degree-corrected cluster tensor family~\eqref{eq:gammafamily}.
    \item For every $\gamma\in \mathbb{R}$, $\{\tX-1\colon \tX\in \tP_{\text{shifted}}(\gamma)\}\subset \tP_{\text{non-degree}}(\gamma)$, where $\tP_{\text{non-degree}}(\gamma)$ denotes the sub-family of rank-one tensor block model constructed in proof of \citet[Theorem 7]{han2022exact}. 
\end{enumerate}
The verification of the above two properties is provided in the end of this proof. 

Now, following the proof of~\citet[Theorem 7]{han2022exact}, when $\gamma<-K/2$, every polynomial-time algorithm estimator $(\hat \mM_k)_{k\in[K]}$ obeys
\begin{align}\label{eq:p}
\liminf_{p\to \infty} \sup_{\tX\in \tP_{\text{non-degree}}(\gamma)}\mathbb{P}(\exists k \in[K],\  \hat \mM_k \neq \mM_k)\geq 1/2,
\end{align}
under the HPC Conjecture~\ref{hypo:HPC}.
The inequality~\eqref{eq:p} implies
\[
\liminf_{p\to \infty} \sup_{\tX\in \tP_{\text{non-degree}}(\gamma)}\max_{k\in[K]}\bbE[p\ell(z_k, \hat z_k)]\geq 1.
\]
Based on properties (i)-(ii), we conclude that
\[
\liminf_{p\to \infty} \sup_{\tX\in \tP(\gamma)}\max_{k\in[K]}\bbE[p\ell(z_k, \hat z_k)]\geq 1.
\]
We complete the proof by verifying the properties (i)-(ii). For (i), we verify that the angle gap for the core tensor $\tS$ in~\eqref{eq:S} is on the order of $\sigma p^{-\gamma/2}$. Specifically, write $\mone=(1,1,1)$ and $\me=(1,-1,0)$. We have
\[
\Mat(\tS)=
\begin{bmatrix}
\text{Vec}(\mone^{\otimes K-1})+\sigma p^{-\gamma/2}\text{Vec}\left(\me^{\otimes (K-1)}\right) \\
\text{Vec}(\mone^{\otimes K-1})-\sigma p^{-\gamma/2}\text{Vec}\left(\me^{\otimes (K-1)}\right) \\
\text{Vec}(\mone^{\otimes K-1})
\end{bmatrix}.
\]
Based on the orthogonality $\langle \mone, \me\rangle=0$, the minimal angle gap among rows of $\Mat(\tS)$ is
\begin{align}
\Delta^2_{\min}(\tS)&\asymp \tan^2(\Mat(\tS)_{1:}, \Mat(\tS)_{3:})
=\left(\onorm{\me}_2 \over \onorm{\mone}_2\right)^{2(K-1)} \sigma^2 d^{-\gamma}\asymp \sigma^2 d^{-\gamma}.
\end{align}
Therefore, we have shown that $\tP_{\text{shifited}}(\gamma)=\tP(\gamma)$. Finally, the property (ii) follows directly by comparing the definition of $\tS$ in~\eqref{eq:S} with that in the proof of \citet[Theorem 7]{han2022exact}. 
\end{proof}

\subsection{Proof of Theorem~\ref{thm:initial} and Proposition~\ref{prop:ber}} \label{sec:initial_prove}

\begin{proof}[Proof of Theorem~\ref{thm:initial}] We prove Theorem~\ref{thm:initial} under the dTBM~\eqref{eq:model_tensor} with symmetric mean tensor, parameters $(z, \tS, \mtheta)$, {fixed $r\geq 1, K \geq 2$}, and i.i.d.\ noise. {For the case $r = 1$, we have $L(z^{(0)},z) = 0, \ell(z^{(0)}, z) = 0$ trivially. Hence, we focus on the proof of the first mode clustering $z^{(0)}_1$ with $r \geq 2$}; 
the proofs for the other modes can be extended similarly. We drop the subscript $k$ in the matricizations $\mM_k, \mX_k, \mS_k$ and in the estimate $z^{(0)}_1$. We firstly show the proof with balanced $\mtheta$.

\textbf{We firstly show the upper bound for misclustering error $\ell(z^{(0)}, z)$.}

First, by Lemma~\ref{lem:angle_gap_x}, there exists a positive constant such that $\min_{z(i) \neq z(j)} \onorm{ \mX^s_{i:} - \mX^s_{j:} } \geq c_0  \Delta_{\min}$. By the balance assumption on $\mtheta$ and Lemma~\ref{lem:upper_mis}, we have 
 \begin{equation}\label{eq:theta_bound}
          \min_{\pi \in \Pi} \sum_{i : z^{(0)}(i) \neq \pi(z(i))} \theta(i)^2  \leq  \sum_{i \in S_I} \theta(i)^2 + 4 \sum_{i \in S} \theta(i)^2 ,
    \end{equation}
    where 
    \begin{equation}
        S_0 = \{i: \onormSize{}{\hat \mX_{i:}} = 0\}, S = \{i \in S_0^c: \onormSize{}{\hat \mx_{z^{(0)}(i)} - \mX_{i:}^s } \geq c_0  \Delta_{\min}/2\}.
    \end{equation}
    \normalsize
    On one hand, note that for any set $P \in [p]$,
    \begin{align}
        \sum_{i \in P} \onormSize{}{\mX_{i:}}^2 &= \sum_{i \in P} \onormSize{}{ \theta(i) \mS_{z(i):} (\mTheta \mM)^{T, \otimes (K-1)} }^2 \\
        & \geq \sum_{i \in P} \theta(i)^2 \min_{a \in [r]} \onormSize{}{\mS_{a:}}^2 \lambda_{r}^{2(K-1)} (\mTheta \mM)\\
        & \gtrsim \sum_{i \in P} \theta(i)^2 p^{K-1} r^{-(K-1)}, 
    \end{align}
    where the last inequality follows Lemma~\ref{lem:singular_thetam}, the assumption that $\min_{i \in [p]} \theta(i) \geq c$, and the constraint $\min_{a \in [r]} \onormSize{}{\mS_{a:}} \geq c_3$ in the parameter space~\eqref{eq:family}. Thus, we have 
    \begin{equation}\label{eq:theta_p}
        \sum_{i \in P} \theta(i)^2 \lesssim \sum_{i \in P} \onormSize{}{\mX_{i:}}^2 p^{-(K-1)} r^{K-1}.
    \end{equation}
    
    On the other hand, note that 
    \begin{align}
         \sum_{i \in S} \onormSize{}{\mX_{i:}}^2
         & \leq 2 \sum_{i \in S} \onormSize{}{\hat \mX_{i:}}^2  +  2\sum_{i \in S} \onormSize{}{\hat \mX_{i:} - \mX_{i:}}^2 \label{eq:tr1}\\
         & \leq \frac{8}{c_0^2 \Delta_{\min}^2  } \sum_{i \in S} \onormSize{}{\hat \mX_{i:}}^2   \onormSize{}{\hat \mx_{z^{(0)}(i)} - \mX_{i:}^s }^2 
          +2 \onormSize{}{\hat \tX - \tX}_F^2 \label{eq:defs} \\
         & \leq \frac{16}{c_0^2 \Delta_{\min}^2  }\sum_{i \in S} \onormSize{}{\hat \mX_{i:}}^2  \off{ \onormSize{}{\hat \mx_{z^{(0)}(i)} - \hat \mX_{i:}^s }^2 + \onormSize{}{\hat \mX_{i:}^s - \mX_{i:}^s }^2 } + 2\onormSize{}{\hat \tX - \tX}_F^2 \label{eq:tr2} \\
         & \leq  \frac{16(1 + \eta)}{c_0^2 \Delta_{\min}^2 }\sum_{i \in S} \onormSize{}{\hat \mX_{i:}}^2  \onormSize{}{\hat \mX_{i:}^s - \mX_{i:}^s }^2  + 2\onormSize{}{\hat \tX - \tX}_F^2 \label{eq:km} \\
         & \leq \of{ \frac{16(1 + \eta)}{c_0^2 \Delta_{\min}^2 } + 2 }\onormSize{}{\hat \tX - \tX}_F^2 \label{eq:lem1}\\
         & \lesssim \of{ \frac{16(1 + \eta)}{c_0^2 \Delta_{\min}^2 } + 2 } \of{ p^{K/2}r + pr^2 + r^K } \sigma^2 \label{eq:lem_err},
    \end{align}
    where inequalities \eqref{eq:tr1} and \eqref{eq:tr2} follow from the triangle inequality, \eqref{eq:defs} follows from the definition of $S$, \eqref{eq:km} follows from the update rule of $k$-means in Step 6 of Sub-algorithm~\hyperref[alg:main]{1}, \eqref{eq:lem1} follows from Lemma~\ref{lem:norm_diff}, and the last inequality~\eqref{eq:lem_err} follows from Lemma~\ref{lem:two-step_esterror}. Also, note that 
    \begin{align}
        \sum_{i \in S_0} \onormSize{}{\mX_{i:}}^2 &=  \sum_{i \in S_0} \onormSize{}{\hat \mX_{i:}- \mX_{i:}}^2 \leq \onormSize{}{\hat \tX - \tX}_F^2 \lesssim \of{ p^{K/2}r + pr^2 + r^K } \sigma^2,\label{eq:s0}
    \end{align}
    where the equation follows from the definition of $S_0$. Therefore, combining the inequalities~\eqref{eq:theta_bound}, \eqref{eq:theta_p}, \eqref{eq:lem_err}, and \eqref{eq:s0}, we have 
    \begin{align}
        \min_{\pi \in \Pi} \sum_{i : z^{(0)}(i) \neq \pi(z(i))} \theta(i)^2 
        &\lesssim \of{\sum_{i \in S} \onormSize{}{\mX_{i:}}^2 + \sum_{i \in S_0} \onormSize{}{\mX_{i:}}^2  } p^{-(K-1)} r^{K-1} \\
        & \lesssim \frac{\sigma^2 r^{K-1}}{ \Delta_{\min}^2  p^{K-1} } \of{ p^{K/2}r + pr^2 + r^K }.\label{eq:theta_sum}
    \end{align}
    With the assumption that $\min_{i \in [p]} \theta(i) \geq c$, we finally obtain the result
    \begin{equation}
        \ell(z^{(0)}, z) \lesssim \frac{1}{p}\min_{\pi \in \Pi} \sum_{i : z^{(0)}(i) \neq \pi(z(i))} \theta(i)^2 \lesssim \frac{r^K p^{-K/2} }{ \text{SNR} },
    \end{equation}
    where the last inequality follows from the definition $\text{SNR} = \Delta_{\min}^2/\sigma^2$.
    
      {Without the balanced $\mtheta$, we have $\min_{z(i) \neq z(j)} \onorm{ \mX^s_{i:} - \mX^s_{j:} } \geq c_0  \Delta_{\mX}$. Replacing the definition of $S$ with $\Delta_{\mX}$, we obtain the desired result.}
    
    \textbf{Next, we show the bound for $L(z^{(0)}, z).$}
    
    Note that $\mX_{i:}^s$ have only $r$ different values. We let $\mX_a^s = \mX_{i:}^s$ for all $i$ such that $z(i) = a, a \in [r]$. 
    Notice that 
\begin{equation}
    \onormSize{}{\mX_{i:}}^2 \gtrsim p^{K-1} r^{-(K-1)} 
\end{equation}
and
\begin{equation}
    \onormSize{}{\mX_{i:} - \hat \mX_{i:}}^2 \leq \onormSize{}{\hat \tX - \tX}_F^2 \lesssim  p^{K/2}r + pr^2 + r^K .
\end{equation}
Therefore, when $p$ is large enough, we have 
\begin{align}
    \sum_{i \in [p]} \onormSize{}{\mX_{i:}}^2 \onormSize{}{\hat \mX_i^s - \hat \mx_{z^{(0)}(i)}}^2
    &\lesssim \sum_{i \in [p]} \of{\onormSize{}{\mX_{i:}}^2 - \onormSize{}{\mX_{i:} - \hat \mX_{i:}}^2} \onormSize{}{\hat \mX_{i:}^s - \hat \mx_{z^{(0)}(i)}}^2\\
   & \lesssim \sum_{i \in [p]} \onormSize{}{\hat \mX_{i:}}^2 \onormSize{}{\hat \mX_{i:}^s - \hat \mx_{z^{(0)}(i)}}^2 \\
   & \lesssim \eta  \sum_{i \in [p]} \onormSize{}{\hat \mX_{i:}}^2 \onormSize{}{\hat \mX_{i:}^s - \mX_{i:}^s}^2 \\
   & \lesssim \onormSize{}{\hat \tX - \tX}_F^2 \\
   & \lesssim p^{K/2}r + pr^2 + r^K .\label{eq:cor_1}
\end{align}

Hence, we have 
\begin{align}
    \sum_{i \in [p]} \onormSize{}{\hat \mX_{i:}^s - \hat \mx_{z^{(0)}(i)}}^2 &\lesssim \sum_{i \in [p]} \theta(i)^2 \onormSize{}{\hat \mX_i^s - \hat \mx_{z^{(0)}(i)}}^2 \\
    &\lesssim \frac{r^{K-1}}{p^{K-1}} \sum_{i \in [p]} \onormSize{}{\mX_{i:}}^2 \onormSize{}{\hat \mX_{i:}^s - \hat \mx_{z^{(0)}(i)}}^2\\
    & \lesssim \frac{r^{K-1}}{p^{K-1}} \of{ p^{K/2}r + pr^2 + r^K },\label{eq:cor_2}
\end{align}
where the first inequality follows from the assumption $\min_{i \in [p]} \theta(i) \geq c > 0$, the second inequality follows from the inequality~\eqref{eq:theta_p}, and the last inequality comes from the inequality~\eqref{eq:cor_1}.

Next, we consider the following quantity,
\begin{align}
    \sum_{i \in [p]} \theta(i) \onormSize{}{ \mX_{i:}^s - \hat \mx_{z^{(0)}(i)}}^2
    &\lesssim \sum_{i \in [p]} \theta(i)^2 \onormSize{}{\mX_{i:}^s - \hat \mX_{i:}^s}^2 + \sum_{i \in [p]} \theta(i)^2 \onormSize{}{\hat \mX_{i:}^s - \hat \mx_{z^{(0)}(i)} }^2\\
    & \lesssim \sum_{i \in [p]} \frac{ \theta(i)^2}{ \onormSize{}{\mX_{i:}}^2} \onormSize{}{\mX_{i:}  - \hat \mX_{i:}}^2 +  \sum_{i \in [p]} \theta(i)^2 \onormSize{}{\hat \mX_{i:}^s - \hat \mx_{z^{(0)}(i)} }^2\\
    & \lesssim \frac{r^{K-1}}{p^{K-1}} \of{ p^{K/2}r + pr^2 + r^K} ,\label{eq:cor_3}
\end{align}
where the first inequality follows from the assumption of $\theta(i)$ and triangle inequality, the second inequality follows from Lemma~\ref{lem:norm_diff}, and the last inequality follows from~\eqref{eq:cor_2}. In addition, with Theorem~\ref{thm:initial} and the condition SNR $\gtrsim p^{-K/2} \log p$, for all $a \in [r]$, we have
\begin{equation}
    |z^{-1}(a) \cap (z^{(0)})^{-1} (a)| \geq |z^{-1}(a)| - p\ell(z^{(0)} ,  z) \gtrsim \frac{p}{r} - \frac{p}{\log p} \gtrsim \frac{p}{r},
\end{equation}
\normalsize
when $p$ is large enough. Therefore, for all $a \in [r]$, we have
\begin{align}
       \onormSize{}{ \hat \mx_{a} - \mX_{a}^s  }^2  &= \frac{ \sum_{i \in z^{-1}(a) \cap (z^{(0)})^{-1} (a) } \onorm{ \mX_{i:}^s - \hat \mx_{z^{(0)}(i)}}^2 }{|z^{-1}(a) \cap (z^{(0)})^{-1} (a) |}  \\
       & \lesssim \frac{r}{p} \of{  \sum_{i \in [p]} \onormSize{}{ \mX_{i:}^s - \hat \mX_{i:}^s }^2 + \sum_{i \in [p]} \onormSize{}{\hat \mX_{i:}^s - \hat \mx_{z^{(0)}(i)}}^2 } \\
       & \lesssim  \frac{r^{K}}{p^{K}}\of{ p^{K/2}r + pr^2 + r^K}, \label{eq:cor_4}
\end{align}
\normalsize
where the last inequality follows from the inequality~\eqref{eq:cor_2}.

Finally, we obtain 
\begin{align}
     L^{(0)} &= \frac{1}{p}  \sum_{i \in [p]} \theta(i) \sum_{b \in [r]}  \ind \offf{ z^{(0)}(i) = b } \onormSize{}{ [ \mS_{  z(i):}  ]^s - [ \mS_{b:}  ]^s  }^2 \\
     & \lesssim \frac{1}{p}  \sum_{i \in [p], z^{(0)}(i) \neq z(i)} \theta(i)  \onormSize{}{ \mX_{i:}^s -  \mX_{z^{(0)}(i)}^s  }^2\\
     & \lesssim \frac{1}{p}  \sum_{i \in [p], z^{(0)}(i) \neq z(i)} \theta(i) \Big( \onormSize{}{\mX_{i:}^s - \hat \mx_{z^{(0)}(i)}}^2 + \onormSize{}{ \hat \mx_{z^{(0)}(i)} - \mX_{z^{(0)}(i)}^s  }^2 \Big) \\
     & \leq \bar C \frac{r^{K}}{p^{K}} \of{ p^{K/2}r + pr^2 + r^K}, \\
     & \leq \frac{\bar C\Delta_{\min}^2}{ \tilde C r \log p}
\end{align}
where the first inequality follows from Lemma~\ref{lem:angle_gap_x}, the third inequality follows from inequalities~\eqref{eq:cor_3} and \eqref{eq:cor_4}, and the last inequality follows from the assumption that SNR $\geq \tilde C p^{-K/2} \log p$.
\end{proof}

\begin{proof}[Proof of Proposition~\ref{prop:ber}]  Sub-algorithm~3 shares the same algorithm strategy as Sub-algorithm~\hyperref[alg:main]{1} but with a different estimation of the mean tensor, $\hat \tX'$. Hence, the proof of Proposition~\ref{prop:ber} follows the same proof idea with the proof of Theorem~\ref{thm:initial}. Replacing the estimation $\hat \tX$ by $\hat \tX'$ in the proof of Theorem~\ref{thm:initial}, we have 
\begin{align}
        & \min_{\pi \in \Pi} \sum_{i : z^{(0)}(i) \neq \pi(z(i))} \theta(i)^2  \lesssim \of{\sum_{i \in S} \onormSize{}{\mX_{i:}}^2 + \sum_{i \in S_0} \onormSize{}{\mX_{i:}}^2  } p^{-(K-1)} r^{K-1}. \label{eq:binary_theta}
    \end{align}
    \normalsize
By inequalities~\eqref{eq:lem1} and \eqref{eq:s0}, we have 
\begin{align}
   \sum_{i \in S} \onormSize{}{\mX_{i:}}^2 &\leq  \of{ \frac{16(1 + \eta)}{c_0^2 \Delta_{\min}^2 } + 2 }\onormSize{}{\hat \tX' - \tX}_F^2, \label{eq:binary_theta1}  \\
   \sum_{i \in S_0} \onormSize{}{\mX_{i:}}^2 &\leq  \onormSize{}{\hat \tX' - \tX}_F^2. \label{eq:binary_theta2}
\end{align}
Hence, it suffices to find the upper bound of the estimation error $\onormSize{}{\hat \tX' - \tX}_F^2$ to complete our proof. Note that the matricization $\Mat_{sq}(\tX) \in \bbR^{p^{\floor{K/2}} \times p^{\ceil{K/2}}}$ has $\text{rank}(\Mat_{sq}(\tX)) \leq r^{\ceil{K/2}}$, and Bernoulli random variables follow the sub-Gaussian distribution with bounded variance $\sigma^2 = 1/4$. Apply Lemma~\ref{lem:lowrank} to $\mY = \Mat_{sq}(\tY), \mX = \Mat_{sq}(\tX)$, and $\hat \mX = \Mat_{sq}(\hat \tX')$. Then, with probability tending to 1 as $ p \rightarrow \infty$, we have 
        \begin{equation}\label{eq:binary_est}
            \onormSize{}{\hat \tX' - \tX}_F^2 = \onormSize{}{\Mat_{sq}(\hat \tX') - \Mat_{sq}(\tX)}_F^2 \lesssim p^{\ceil{K/2}}.
        \end{equation}
Combining the estimation error~\eqref{eq:binary_est} with inequalities~\eqref{eq:binary_theta1}, \eqref{eq:binary_theta2}, and \eqref{eq:binary_theta}, we obtain 
\begin{equation}\label{eq:binary_theta_sum}
    \min_{\pi \in \Pi} \sum_{i : z^{(0)}(i) \neq \pi(z(i))} \theta(i)^2 \lesssim  \frac{\sigma^2 r^{K-1}}{ \Delta_{\min}^2  p^{K-1} } p^{\ceil{K/2}}.
\end{equation}
Replace the inequality~\eqref{eq:theta_sum} in the proof of Theorem~\ref{thm:initial} by inequality~\eqref{eq:binary_theta_sum}. With the the same procedures to obtain $\ell(\hat z^{(0)}, z)$ and $L(\hat z^{(0)}, z)$ for Theorem~\ref{thm:initial} , we finish the proof of Proposition~\ref{prop:ber}.
\end{proof}

{\bf Useful Definitions and Lemmas for the Proof of Theorem~\ref{thm:initial}} 

\begin{lem}[Basic inequality]\label{lem:norm_diff} For any two nonzero vectors $\mv_1,\mv_2$ of same dimension, we have 
\[
\sin(\mv_1 , \mv_2) \leq \onorm{{\mv_1^s}-{\mv_2^s}}\leq {2\onorm{\mv_1-\mv_2}\over \max\left(\onorm{\mv_1},\onorm{\mv_2}\right)}.
\]
\end{lem}
\begin{proof}[Proof of Lemma~\ref{lem:norm_diff}]
For the first inequality, let $\alpha \in [0,\pi]$ denote the angle between $\mv_1$ and $\mv_2$. We have 
\begin{align}
     \onorm{{\mv_1^s}-{\mv_2^s}} = \sqrt{2(1 - \cos \alpha)} = 2 \sin \frac{\alpha}{2} \geq \sin \alpha,
\end{align}
where the equations follow from the properties of trigonometric function and the inequality follows from the fact the $\cos \frac{\alpha}{2} \leq 1$ and $\sin \alpha = 2 \sin \frac{\alpha}{2} \cos \frac{\alpha}{2} > 0$ for $\alpha \in [0, \pi]$. 

For the second inequality, without loss of generality, we assume $\onorm{\mv_1}\geq \onorm{\mv_2}$. Then
\begin{align}
\onorm{{\mv_1^s}-{\mv_2^s}}&=\onorm{{\mv_1\over \onorm{\mv_1}}- {\mv_2\over \onorm{\mv_1}}+{\mv_2\over \onorm{\mv_1}}-{\mv_2\over \onorm{\mv_2}}}\\
&\leq {\onorm{\mv_1-\mv_2}\over \onorm{\mv_1}}+{\onorm{\mv_2}\onorm{\mv_1}-\onorm{\mv_2}\over \onorm{\mv_1}\onorm{\mv_2}}\\
&\leq {2\onorm{\mv_1-\mv_2}\over \onorm{\mv_2}}.
\end{align}

Therefore, Lemma~\ref{lem:norm_diff} is proved.
\end{proof}

\begin{defn}[Weighted padding vectors]\label{def:pad} For a vector $\ma = \entry{a_i} \in \bbR^d$, we define the padding vector of $\ma$ with the weight collection $\mw = \{\mw_i\colon \mw_i = \entry{w_{ik}} \in \bbR^{p_i}\}_{i = 1}^d$ as
\begin{equation}\label{eq:paddef}
    \pad_{\mw}(\ma) = [a_1 \circ \mw_1, \ldots, a_d \circ \mw_d]^T, 
\end{equation}
where $ a_i \circ \mw_i = [a_i w_{i1},\ldots, a_i w_{i p_i}]^T, \text{ for all } i \in [d].$
Here we also view $\pad_{\mw}(\cdot)\colon \bbR^{d}\mapsto \bbR^{\sum_{i\in[d]}p_i}$ as an operator. 
We have the bounds of the weighted padding vector
\begin{equation}\label{eq:pad_bound}
     \min_{i \in [d]} \onormSize{}{\mw_i}^2 \onormSize{}{\ma}^2 \leq \onormSize{}{\pad_{\mw}(\ma)}^2 \leq \max_{i \in [d]} \onormSize{}{\mw_i}^2 \onormSize{}{\ma}^2.
\end{equation}
Further, we define the inverse weighted padding operator $\pad^{-1}: \bbR^{\sum_{i\in[d]}p_i}\mapsto \bbR^d$ which satisfies 
\begin{equation}
    \pad^{-1}_{\mw}(\pad_{\mw}(\ma)) = \ma.
\end{equation}
\end{defn}

\begin{lem}[Angle for weighted padding vectors]\label{lem:pad} Suppose that we have two non-zero vectors $\ma, \mb \in \bbR^d$. Given the weight collection $\mw$, we have 
\begin{align}
    \frac{\min_{i \in [d]} \onormSize{}{\mw_i}}{\max_{i \in [d]} \onormSize{}{\mw_i} } \sin (\ma, \mb) & \stackrel{*}{\leq} \sin(\pad_{\mw}(\ma),  \pad_{\mw}(\mb))\stackrel{**}{\leq} \frac{ \max_{i \in [d]} \onormSize{}{\mw_i} }{\min_{i \in [d]} \onormSize{}{\mw_i}} \sin (\ma, \mb). \label{eq:lempad}
\end{align}
\end{lem}

\begin{proof}[Proof of Lemma~\ref{lem:pad}] We prove the two inequalities separately with similar ideas.

First, we prove the inequality ** in~\eqref{eq:lempad}. Decomposing $\mb$ yields
\begin{equation}
    \mb = \cos (\ma, \mb) \frac{\onormSize{}{\mb}}{\onormSize{}{\ma}} \ma + \sin (\ma, \mb) \frac{\onormSize{}{\mb}}{\onormSize{}{\ma^{\perp}}}\ma^{\perp},
\end{equation}
where $\ma^{\perp} \in \bbR^d$ is in the orthogonal complement space of $\ma$. By the Definition~\ref{def:pad}, we have 
\begin{equation}\label{eq:pad_second_ineq}
    \pad_{\mw}(\mb) = \cos (\ma, \mb) \frac{\onormSize{}{b}}{\onormSize{}{a}} \pad_{\mw}(\ma) + \sin (\ma, \mb) \frac{\onormSize{}{b}}{\onormSize{}{a^{\perp}}} \pad_{\mw}(\ma^{\perp}).
\end{equation}
Note that $\pad_{\mw}(\ma^{\perp})$ is not necessary equal to the orthogonal vector of $\pad(\ma)$; i.e., $\pad_{\mw}(\ma^{\perp}) \neq (\pad_{\mw}(\ma))^{\perp}$. By the geometry property of trigonometric functions, we obtain
\begin{align}
    \sin(\pad_{\mw}(\ma),  \pad_{\mw}(\mb))  &\leq \frac{  \onormSize{}{\mb} \onormSize{}{\pad_{\mw}(\ma^{\perp})} }{ \onormSize{}{\ma^{\perp}} \onormSize{}{ \pad_{\mw}(\mb)}} \sin (\ma, \mb) \leq  \frac{ \max_{i \in [d]} \onormSize{}{\mw_i} }{\min_{i \in [d]} \onormSize{}{\mw_i}} \sin (\ma, \mb),
\end{align}
where the second inequality follows by applying the property~\eqref{eq:pad_bound} to vectors $\mb$ and $\ma^{\perp}$. 

Next, we prove inequality * in~\eqref{eq:lempad}. With the decomposition of $\pad_{\mw}(\mb)$ and the inverse weighted padding operator, we have 
\begin{align}
    \mb &=  \cos(\pad_{\mw}(\ma), \pad_{\mw}(\mb)) \frac{\onormSize{}{\pad_{\mw}(\mb)}}{\onormSize{}{\pad_{\mw}(\ma)}} \ma  +  \sin(\pad_{\mw}(\ma), \pad_{\mw}(\mb))\frac{\onormSize{}{\pad_{\mw}(\mb)}}{\onormSize{}{(\pad_{\mw}(\ma))^{\perp}}} \pad_{\mw}^{-1}((\pad_{\mw}(\ma))^{\perp}).
\end{align}
\normalsize
Therefore, we obtain 
\begin{align}
    \sin(\ma, \mb) & \leq \frac{ \onormSize{}{\pad_{\mw}(\mb) } \onormSize{}{\pad_{\mw}^{-1}((\pad_{\mw}(\ma))^{\perp})  } }{ \onormSize{}{(\pad_{\mw}(\ma))^{\perp} } \onormSize{}{\mb} }  \sin(\pad_{\mw}(\ma), \pad_{\mw}(\mb))\\
    & \leq \frac{ \max_{i \in [d]} \onormSize{}{\mw_i} }{\min_{i \in [d]} \onormSize{}{\mw_i}} \sin(\pad_{\mw}(\ma), \pad_{\mw}(\mb)),
\end{align}
where the second inequality follows by applying the property~\eqref{eq:pad_bound} to vectors $\mb$ and $\pad_{\mw}^{-1}((\pad_{\mw}(\ma))^{\perp})$.
\end{proof}

\begin{lem}[Singular value of weighted membership matrix]\label{lem:singular_thetam} Under the parameter space~\eqref{eq:family} and assumption that $\min_{i \in [p]} \theta(i) \geq c$ {for some constant $c >0$}, the singular values of $\mTheta \mM$ are bounded as 
\begin{align}
    \sqrt{p/r} &\lesssim \sqrt{ \min_{a \in [r]}\onormSize{}{ \mtheta_{ z^{-1}(a)} }^2 }  \leq \lambda_{r}(\mTheta \mM) \leq \onormSize{}{\mTheta \mM}_{\sigma} \leq\sqrt{ \max_{a \in [r]}\onormSize{}{ \mtheta_{ z^{-1}(a)} }^2 } \lesssim  p /r . \label{eq:membership2}
\end{align}
\end{lem}

\begin{proof}[Proof of Lemma~\ref{lem:singular_thetam}] Note that 
\begin{equation}
    (\mTheta\mM)^T \mTheta \mM = \mD,
\end{equation}
with $\mD = \text{diag}(D_1,\ldots, D_r)$ where $D_a = \onormSize{}{\mtheta_{z^{-1}(a)}}^2, a \in [r]$.
By the definition of singular values, we have 
\begin{equation}
     \sqrt{\min_{a \in [r]}\onormSize{}{ \mtheta_{ z^{-1}(a)} }^2 } \leq \lambda_r(\mTheta \mM) \leq \onormSize{}{\mTheta\mM}_{\sigma} \leq \sqrt{ \max_{a \in [r]}\onormSize{}{ \mtheta_{ z^{-1}(a)} }^2}.
\end{equation}

Since that $\min_{i \in [p]}\theta(i) \geq c$ by the constraints in parameter space, we have  
\begin{equation}
   \min_{a \in [r]} \onormSize{}{ \mtheta_{ z^{-1}(a)} }^2 \geq c^2 \min_{a \in [r]} |z^{-1}(a)| \gtrsim \frac{p}{r},
\end{equation}
where the last inequality follows from the constraint in parameter space~\eqref{eq:family}. Finally, notice that 
\begin{equation}
    \sqrt{\max_{a \in [r]}\onormSize{}{ \mtheta_{ z^{-1}(a)} }^2 } \leq \max_{a \in [r]} \sqrt{ \onormSize{}{ \mtheta_{ z^{-1}(a)} }^2_1} \lesssim \frac{p}{r}.
\end{equation}

 Therefore, we complete the proof of Lemma~\ref{lem:singular_thetam}.
\end{proof}

\begin{lem}[Singular-value gap-free tensor estimation error bound]\label{lem:two-step_esterror}Consider an order-$K$ tensor $\tA = \tX + \tZ \in \bbR^{p \times \cdots \times p}$, where $\tX$ has Tucker rank $(r,...r)$ and $\tZ$ has independent sub-Gaussian entries with parameter $\sigma^2$. Let $\hat \tX$ denote the double projection estimated tensor in Step 2 of Sub-algorithm~\hyperref[alg:main]{1} in the main paper. Then with probability at least $1 - C \exp\of{- cp }$, we have
\begin{align}
    \onormSize{}{\hat \tX - \tX}_F^2 \leq C \sigma^2 \of{ p^{K/2}r + pr^2 + r^K },
\end{align}
where $C, c$ are some positive constants.
\end{lem}

\begin{proof}[Proof of Lemma~\ref{lem:two-step_esterror}]
See \citet[Proposition 1]{han2022exact}.
\end{proof}

\begin{lem}[Upper bound of misclustering error]\label{lem:upper_mis} Let $z: [p] \mapsto [r]$ be a cluster assignment such that $|z^{-1}(a)| \asymp p/r$ for all $a \in [r]$ {with $r \geq 2$}. Let node $i$ correspond to a vector $\mx_i  = \theta(i) \mv_{z(i)} \in \bbR^d$, where $\{\mv_a\}_{a = 1}^r$ are the cluster centers and $\mtheta = \entry{\theta(i)} \in \bbR^p_+$ is the positive degree heterogeneity.  Assume that  $\mtheta$ satisfies the balanced assumption~\eqref{eq:degree} such that ${\max_{a \in [r]} \onormSize{}{\mtheta_{z^{-1}(a)}}^2 \over \min_{a \in [r]} \onormSize{}{\mtheta_{z^{-1}(a)}}^2 }= 1 + o(1)$. Consider an arbitrary estimate $\hat z$ with $\hat \mx_i = \hat \mv_{\hat z(i)}$ for all $ i \in S$. Then, if
\begin{equation}\label{eq:upper_mis_cond}
    \min_{a \neq b \in [r]} \onormSize{}{\mv_a - \mv_b} \geq 2c,
\end{equation}
 for some constant $c >0$, we have 
\begin{equation}
    \min_{\pi \in \Pi} \sum_{i : \hat z(i) \neq \pi(z(i))} \theta(i)^2 \leq \sum_{i \in S_0} \theta(i)^2 + 4 \sum_{i \in S} \theta(i)^2,
\end{equation}
where $S_0$ is defined in {Step 4} of Sub-algorithm~\hyperref[alg:main]{1} and
\begin{equation}
   S = \{i \in S_0^c: \onormSize{}{\hat \mx_i - \mv_{z(i)} } \geq c \}.
\end{equation}

\end{lem}

\begin{proof}[Proof of Lemma~\ref{lem:upper_mis}] 

For each cluster $u\in[r]$, we use $C_u$ to collect the subset of points for which the estimated and true positions $\hat \mx_i, \mx_i$ are within distance $c$. Specifically, define
\begin{equation}
    C_u = \{ i \in z^{-1}(u) \cap S_0^c: \onormSize{}{\hat \mx_i - \mv_{z(i)}} < c \},
\end{equation}
and divide $[r]$ into three groups based on $C_u$ as 
\begin{align}
    R_1 &= \{ u \in [r]: C_u = \emptyset \},\\
    R_2 &= \{ u \in [r]: C_u \neq \emptyset, \text{ for all } i, j \in C_u, \hat z(i) = \hat z(j) \},\\
    R_3 &= \{ u \in [r]: C_u \neq \emptyset, \text{ there exist } i, j \in C_u, \hat z(i) \neq \hat z(j) \}.
\end{align}
Note that $\cup_{u \in [r]}C_u = S_0^c/S^c$ and $C_u \cap C_v = \emptyset$ for any $u \neq v$. Suppose there exist $ i \in C_u$ and $j \in C_v$ with $u \neq v \in [r]$ and $\hat z(i) = \hat z(j)$. Then we have 
\begin{equation}
    \onormSize{}{\mv_{z(i)} - \mv_{z(j)} } \leq  \onormSize{}{\mv_{z(i)} - \hat \mx_{i} } + \onormSize{}{\mv_{z(j)} - \hat \mx_j } < 2c,
\end{equation}
which contradicts to the assumption~\eqref{eq:upper_mis_cond}. Hence, the estimates $\hat z(i) \neq \hat z(j)$ for the nodes $ i \in C_u$ and $j \in C_v$ with $u \neq v$. By the definition of $R_2$, the nodes in $\cup_{u \in R_2} C_u$ have the same assignment with $z$ and $\hat z$. Then, 
we have 
\begin{equation}
    \min_{\pi \in \Pi} \sum_{i : \hat z(i) \neq \pi(z(i))} \theta(i)^2 \leq \sum_{i \in S_0} \theta(i)^2 + \sum_{i \in S} \theta(i)^2 + \sum_{ i \in \cup_{u \in R_3} C_u} \theta(i)^2.
\end{equation}
\normalsize
 We only need to bound $\sum_{ i \in \cup_{u \in R_3} C_u} \theta(i)^2$ to finish the proof. Note that every $C_u$ with $u \in R_3$ contains at least two nodes assigned to different clusters by $\hat z$. Then, we have $|R_2| + 2 |R_3| \leq r$. Since $|R_1| + |R_2| + |R_3| = r$, we have $|R_3| \leq |R_1|$. Hence, we obtain
\begin{align}
    \sum_{ i \in \cup_{u \in R_3} C_u} \theta(i)^2 &\leq |R_3| \max_{a \in [r]} \onormSize{}{\mtheta_{z^{-1}(a)}}^2 \\
    & \leq |R_1| \max_{a \in [r]} \onormSize{}{\mtheta_{z^{-1}(a)}}^2  \\
    & \leq {\max_{a \in [r]} \onormSize{}{\mtheta_{z^{-1}(a)}}^2 \over \min_{a \in [r]} \onormSize{}{\mtheta_{z^{-1}(a)}}^2 } \sum_{i \in \cup_{u \in R_1} z^{-1}(u)} \theta(i)^2\\
    & \leq 2 \sum_{i \in S} \theta(i)^2,
\end{align}
where the last inequality holds by the balanced assumption on $\mtheta$ when $p$ is large enough, and the fact that $ \cup_{u \in R_1} z^{-1}(u) \subset S$.
\end{proof}

\begin{lem}[Low-rank matrix estimation] \label{lem:lowrank} Let $\mY = \mX + \mE \in \bbR^{m \times n}$, where {$n > m$ and} $\mE$ contains independent mean-zero sub-Gaussian entries with bounded variance $\sigma^2$. Suppose $\text{rank}(\mX) = r$. Consider the least square estimator 
        \begin{equation}
            \hat \mX = \argmin_{\mX' \in \bbR^{m \times n}, \text{rank}(\mX') \leq r} \onormSize{}{\mX' - \mY}_F^2.
        \end{equation}
        There exist positive constants $C_1, C_2$ such that
        \begin{equation}
            \onormSize{}{\hat \mX - \mX}_F^2 \leq C_1 \sigma^2 nr,
        \end{equation}
        with probability at least $1 - \exp(-C_2 nr)$.

        \end{lem}

\begin{proof}[Proof of Lemma~\ref{lem:lowrank}] Note that $\onormSize{}{\hat \mX - \mY}_F^2 \leq \onormSize{}{\mX - \mY}_F^2$ by the definition of least square estimator. 
{
We have 
        \begin{align}
            \onormSize{}{\hat \mX - \mX}_F^2
            & \leq 2 \ang{ \hat \mX - \mX, \mY - \mX} \\ & \leq 2 \onormSize{}{\hat \mX - \mX}_F \sup_{\mT \in \bbR^{m \times n}, \text{rank}(\mT) \leq 2r, \onormSize{}{\mT}_F = 1} \ang{\mT, \mY - \mX} \label{eq:lem9_1}
        \end{align}
        \normalsize
        with probability at least $1 - \exp(-C_2 nr)$, where the second inequality follows by re-arrangement.
        
        Consider the SVD for matrix $\mT = \mU \Sigma \mV^T$ with orthogonal matrices $\mU \in \bbR^{m \times 2r}, \mV \in \bbR^{n \times 2r}$ and diagonal matrix $\Sigma \in \bbR^{2r \times 2r}$. We have 
        \begin{align}
            \sup_{\mT \in \bbR^{m \times n}, \text{rank}(\mT) \leq 2r, \onormSize{}{\mT}_F = 1} \ang{\mT, \mY - \mX} 
            =& \sup_{\mT \in \bbR^{m \times n}, \text{rank}(\mT) \leq 2r, \onormSize{}{\mT}_F = 1} \ang{\mU \Sigma, \mE \mV} \\
            =&  \sup_{\mv \in \bbR^{2nr}} \mv^T \me \leq C \sigma \sqrt{nr}, \label{eq:lem9_2}
        \end{align}
        with probability $1 - \exp(-C_2 nr)$, where $C, C_2$ are two positive constants, the vectorization $\me = \text{Vec}(\mE \mV) \in \bbR^{2nr}$ has independent mean-zero sub-Gaussian entries with bounded variance $\sigma^2$ due to the orthogonality of $\mV$, and the last inequality follows from \citet[Theorem 1.19]{rigollet2015high}.

        Combining inequalities~\eqref{eq:lem9_1} and \eqref{eq:lem9_2}, we obtain the desired conclusion.
        }

\end{proof}

\subsection{Proofs of Theorem~\ref{thm:stats} (Achievability) and Theorem~\ref{thm:refinement}}\label{sec:statprove2}

\begin{proof}[Proof of Theorem~\ref{thm:stats} (Achievability) and Theorem~\ref{thm:refinement}]

The proofs of Theorem~\ref{thm:stats} (Achievability) and Theorem~\ref{thm:refinement} share the same idea. We prove the contraction step by step. In each step, we show the specific procedures for the algorithm loss and address the MLE loss by stating the difference. 

We consider dTBM~\eqref{eq:model_tensor} with symmetric mean tensor, parameters $(z, \tS, \mtheta)$, {fixed $r\geq 1, K \geq 2$}, and i.i.d.\ noise.  Let $(\hat z,\hat \tS,  \hat \mtheta)$ denote the MLE in \eqref{eq:mle}, and $(z^{(0)}_k, \tS^{(0)}, \mtheta^{(0)}_k)$ denote parameters related to the initialization.  {For the case $r = 1$, $\ell(z^{(t)}_k, z) = 0$ trivially for all $t \geq 0, k \in [k]$. Hence, we focus on the proof of the first mode clustering $z^{(t+1)}_1$ with $r \geq 2$}; 
the extension for other modes can be obtained similarly. We drop the subscript $k$ in the matricizations $\mTheta, \mM_k, \mS_k, \mX_k$ and in estimates $z^{(0)}_k, z^{(t+1)}_k, z^{(t)}_k$ for ease of the notation.  Without loss of generality, we assume that the variance $\sigma = 1$, and that the identity permutation minimizes the initial misclustering error; i.e., $\pi^{(0)} = \argmin_{\pi \in \Pi} \sum_{i \in [p]}\ind \offf{z^{(0)}(i) \neq \pi \circ z(i) }$ and $\pi^{(0)}(a) = a$ for all $ a \in [r]$, and so for $\hat z$.

{\bf Step 1 (Notation and conditions).} We first introduce additional notations and the necessary conditions used in the proof. We will verify that the conditions hold in our context under high probability in the last step of the proof. 

{

\textbf{Notation.}
\begin{enumerate}[wide]

\item Projection. We use $\mI_d$ to denote the identity matrix of dimension $d$. For a vector $\mv \in \bbR^d$, let $\text{Proj}(\mv) \in \bbR^{d \times d}$ denote the projection matrix to $\mv$. Then, $\mI_d - \text{Proj}(\mv)$ is the projection matrix to the orthogonal complement $\mv^{\perp}$. 

    \item We define normalized membership matrices
    \begin{equation}
        \mW = \mM \of{ \text{diag}(\mone_{p}^T \mM) }^{-1}, \mW^{(t)} = \mM^{(t)} \of{ \text{diag}(\mone_{p}^T \mM^{(t)}) }^{-1},
    \end{equation}
    \normalsize
    weighted normalized membership matrices
    \begin{align}
        &\mP = \mTheta \mM ( \text{diag}( \onormSize{}{\mtheta_{z^{-1}(1)}}^2 , \ldots, \onormSize{}{\mtheta_{z^{-1}(r)}}^2 ) )^{-1}, \  \hat \mP = \hat \mTheta \hat \mM ( \text{diag}( \onormSize{}{\hat \mtheta_{z^{-1}(1)}}^2 , \ldots, \onormSize{}{\hat \mtheta_{z^{-1}(r)}}^2 ) )^{-1},
    \end{align}
    and the dual normalized and dual weighted normalized membership matrices
    \begin{align}
         &\mV = \mW^{\otimes (K-1)}, \quad \mV^{(t)} =\of{ \mW^{(t)}}^{\otimes (K-1)}, \quad \mQ = \mP^{\otimes K-1}, \quad  \hat \mQ = \hat \mP^{\otimes K-1}.
    \end{align}
    Also, let $\mB = (\mTheta \mM)^{\otimes (K-1)}, \hat \mB = (\hat \mTheta \hat \mM)^{\otimes (K-1)}$. By the definition, we have $\mB^T \mQ =\hat \mB^T \hat \mQ =\mI_{r^{K-1}}$.
    \item We use $\tS^{(t)}$ to denote the estimator of $\tS$ in the $t$-th iteration, $\hat \tS$ for MLE, $\tilde \tS$ to denote the oracle estimator of $\tS$ given true assignment $z$, and $\bar \tS$ for weighted oracle estimator; i.e.,
    \begin{align}
        &\tS^{(t)} = \tY \times_1 \of{\mW^{(t)}}^T \times_2 \cdots \times_K \of{\mW^{(t)}}^T,& \quad &\tilde \tS = \tY \times_1 \mW^T \times_2 \cdots \times_K \mW^T, \\
        & \hat \tS =  \tY \times_1 \hat \mP^T \times_2 \cdots \times_K \hat \mP^T,& \quad &\bar \tS  = \tY \times_1 \mP^T \times_2 \cdots \times_K \mP^T.
    \end{align}
    \item We define the matricizations of tensors
    \begin{align}
        \mS = \mat(\tS), \  \mY = \mat(\tY), \  \mX = \mat(\tX), \  \mE = \mat(\tE),
    \end{align}
    \begin{equation}
        \mS^{(t)} = \mat(\tS^{(t)}), \   \hat \mS = \mat(\hat \tS), \  \tilde \mS = \mat(\tilde \tS), \  \bar \mS = \mat( \bar \tS ).
    \end{equation}
    \item We define the extended core tensor on $K-1$ modes
    \begin{equation}
        \mA = \mS \mB^T , \quad \bar \mA = \bar \mS \mB^T, \quad \hat \mA = \hat \mS \hat \mB^T.
    \end{equation}
    By the assumption in parameter space~\eqref{eq:family}, we have $\mA = \mP \mX = \mW \mX, \quad \hat \mA = \hat \mP \hat \mX = \hat \mW \hat \mX.$

    \item We define the angle-based misclustering loss in the $t$-th iteration and loss for MLE
    \begin{align}
        L^{(t)} &= \frac{1}{p}  \sum_{i \in [p]} \theta(i) \sum_{b \in [r]}  \ind \{ z^{(t)}(i) = b \} \onormSize{}{ [ \mS_{ z(i):}  ]^s - [ \mS_{b:}  ]^s  }^2, \\
         L(\hat z)& = \frac{1}{p} \sum_{i \in [p]} \theta(i)^2 \sum_{b \in [r]} \ind\{ \hat z(i) = b\} \onormSize{}{[\mA_{z(i):}]^s - [\mA_{b:}]^s}^2.
    \end{align}
    We also define the loss for oracle and weighted oracle estimators
    \small
    \begin{align}
         \xi & = \frac{1}{p} \sum_{i \in [p]} \theta(i) \sum_{b \in [r]} \ind \Big\{\ang{ \mE_{i:} \mV, [  \tilde \mS_{z(i):} ]^s - [  \tilde \mS_{b:} ]^s}  \leq - \frac{ \theta(i) m}{4} \onormSize{}{ [ \mS_{z(i):}  ]^s - [ \mS_{b:}  ]^s  }^2 \Big\} \cdot \onormSize{}{ [ \mS_{z(i):}  ]^s - [ \mS_{ b:}  ]^s  }^2,\\
         \xi' & = \frac{1}{p} \sum_{i \in [p]} \theta(i)^2 \sum_{b \in [r]} \ind \Big\{\ang{\mE_{i:}, [\bar \mA_{z(i):}]^s - [\bar \mA_{b:}]^s } \leq - \frac{m'}{4}  \sqrt{\frac{p^{K-1}}{r^{K-1}}} \onormSize{}{ [\mA_{z(i):}]^s -  [\mA_{b:}]^s  }_F^2 \Big\} \cdot \onormSize{}{[\mA_{z(i):}]^s - [\mA_{b:}]^s}^2.
    \end{align}
    \normalsize
    where $m$ and $m'$ are some positive universal constants.
\end{enumerate}
Then we introduce the necessary conditions in Condition~\ref{cond:origin}.

\begin{condition}(Intermediate results) \label{cond:origin} Let $\mathbb{O}_{p, r}$ denote the collection of all the $p$-by-$r$ matrices with orthonormal columns. We have 
\begin{equation}\label{eq:cond1}
    \onormSize{}{\mE \mV}_\sigma \lesssim \sqrt{\frac{r^{K-1}}{p^{K-1}}} \of{ p ^{1/2}+ r^{(K-1)/2}}, \  \onormSize{}{\mE \mV}_F \lesssim \sqrt{\frac{ r^{2(K-1)}}{p^{K-2}}},\ \onormSize{}{ \mW_{a:}^T \mE \mV} \lesssim \frac{r^K}{p^{K/2}}, \ a \in [r], 
\end{equation}
\begin{equation}\label{eq:cond2}
    \sup_{ \mU_k \in \mathbb{O}_{p, r}, k = 2,\ldots, K } \onormSize{}{ \mE (\mU_{2} \otimes \cdots \otimes \mU_K)}_\sigma \lesssim \of{ \sqrt{r^{K-1}} + K\sqrt{pr}},
\end{equation}
\begin{equation}\label{eq:cond3}
    \sup_{ \mU_k \in \mathbb{O}_{p, r}, k = 2,\ldots, K } \onormSize{}{ \mE (\mU_{2} \otimes \cdots \otimes \mU_K)}_F \lesssim \of{ \sqrt{pr^{K-1}} + K\sqrt{pr}},
\end{equation}
\begin{equation}\label{eq:cond_oracle}
    \xi \leq \exp\of{ - M \frac{\Delta_{\min}^2 p^{K-1}}{r^{K-1}}}, \quad  \xi' \lesssim \exp\of{ -  \frac{\Delta_{\min}^2 p^{K-1}}{r^{K-1}}},
\end{equation}

\begin{equation}\label{eq:cond_intial}
    L^{(t)} \leq \frac{\bar C}{\tilde C} \frac{\Delta_{\min}^2}{r \log p}, \quad \text{for} \quad t = 0, 1, \ldots, T, \quad  L(\hat z) \leq \frac{\bar C}{\tilde C} \frac{\Delta_{\min}^2}{r \log p},
\end{equation}
where $M$ is a positive universal constant in inequality~\eqref{eq:constant_M}, $\bar C, \tilde C$ are positive universal constants in {the proof of Theorem~\ref{thm:initial}} and assumption SNR $\geq \tilde C p^{-K/2} \log p$, respectively. Further, inequality~\eqref{eq:cond1} holds by replacing $\mV$ to $\mV^{(t)}, \mQ, \hat \mQ$ and $\mW_{:a}$ to $\mW_{:a}^{(t),T}, \mP_{:a}^T, \hat \mP_{:a}^T$ when initialization condition~\eqref{eq:cond_intial} holds.
\end{condition}

}

{\bf Step 2 (Misclustering loss decomposition).} Next, we derive the upper bound of $L^{(t+1)}$ for $t = 0 ,1, \ldots, T-1$. By Sub-algorithm~\hyperref[alg:main]{2}, we update the assignment in $t$-th iteration via
    \begin{equation}
        z^{(t+1)}(i) = \argmin_{a \in [r]} \onormSize{}{ [ \mY_{i:} \mV^{(t)}  ]^s - [\mS_{a:}^{(t)}]^s }^2,
    \end{equation}
    following the facts that $\onormSize{}{\ma^s - \mb^s}^2 = 1 - \cos(\ma, \mb)$ for vectors $\ma,\mb$ of same dimension and $\mat(\tY^\text{d}) = \mY \mV^{(t)}$ where $\tY^\text{d}$ is the reduced tensor defined in Step 8 of Sub-algorithm~\hyperref[alg:main]{2}. Then the event $z^{(t+1)}(i) = b$ implies
    \begin{equation}\label{eq:event}
        \onormSize{}{ [ \mY_{i:} \mV^{(t)}  ]^s - [\mS_{b:}^{(t)}]^s }^2 \leq \onormSize{}{ [ \mY_{i:} \mV^{(t)}  ]^s - [\mS_{z(i):}^{(t)}]^s }^2.
    \end{equation}
Note that the event~\eqref{eq:event} also holds for the degenerate entity $i$ with $\onormSize{}{\mY_{i:} \mV^{(t)}} = 0$ due to the convention that $\ma^s = {\bf 0}$ if $\ma = {\bf 0}$. 
    Arranging the terms in \eqref{eq:event} yields the decomposition
    \begin{align}
        &2 \ang{ \mE_{i:} \mV, [  \tilde \mS_{z(i):} ]^s - [  \tilde \mS_{b:} ]^s }  \leq \onormSize{}{\mX_{i:} \mV^{(t)}} \of{ - \onormSize{}{ [ \mS_{z(i):}  ]^s - [ \mS_{b:}  ]^s  }^2 + G_{ib}^{(t)} + H_{ib}^{(t)} } +  F_{ib}^{(t)},\label{eq:decomp}
    \end{align}
    \normalsize
    where
    \small
    \begin{align}
        F_{ib}^{(t)} &= 2 \ang{\mE_{i:} \mV^{(t)}, \of{ [  \tilde \mS_{z(i):} ]^s  -  [  \mS_{z(i):}^{(t)} ]^s  }  - \of{  [  \tilde \mS_{b:} ]^s  -  [  \mS_{b:}^{(t)} ]^s  }  }  + 2 \ang{ \mE_{i:} \of{ \mV - \mV^{(t)} }, [  \tilde \mS_{z(i):} ]^s - [  \tilde \mS_{b:} ]^s  }, \label{eq:f}\\
        G_{ib}^{(t)} &=   \of{ \onormSize{}{ [\mX_{i:} \mV^{(t)}]^s -  [  \mS_{z(i):}^{(t)}  ]^s}^2 -  \onormSize{}{ [\mX_{i:} \mV^{(t)}]^s -  [  \mW_{:z(i)}^T \mY \mV^{(t)} ]^s}^2}   \\
        & \quad \quad \quad \quad \quad -   \of{ \onormSize{}{ [\mX_{i:} \mV^{(t)}]^s -  [  \mS_{b:}^{(t)}  ]^s}^2 -  \onormSize{}{ [\mX_{i:} \mV^{(t)}]^s -  [  \mW_{:b}^T \mY \mV^{(t)} ]^s}^2 }, \label{eq:g} \\
        H_{ib}^{(t)} &=   \onormSize{}{ [\mX_{i:} \mV^{(t)}]^s -  [  \mW_{:z(i)}^T \mY \mV^{(t)} ]^s}^2 - \onormSize{}{ [\mX_{i:} \mV^{(t)}]^s -  [  \mW_{:b}^T \mY \mV^{(t)} ]^s}^2  + \onormSize{}{ [ \mS_{z(i):}  ]^s - [ \mS_{b:}  ]^s  }^2 . \label{eq:h}
    \end{align}
    \normalsize

   Therefore, the event $ \ind \offf{ z^{(t+1)}(i) = b }$ can be upper bounded as
   \begin{align}
         \ind \offf{ z^{(t+1)}(i) = b  }& \leq \ind \offf{ z^{(t+1)}(i) = b ,   \ang{ \mE_{j:} \mV, [  \tilde \mS_{z(i):} ]^s - [  \tilde \mS_{b:} ]^s }  \leq - \frac{1}{4} \onormSize{}{\mX_{i:} \mV^{(t)}} \onormSize{}{ [ \mS_{z(i):}  ]^s - [ \mS_{b:}  ]^s  }^2 } \\
         & + \ind \offf{z^{(t+1)}(i) = b, \frac{1}{2}\onormSize{}{ [ \mS_{z(i):}  ]^s - [ \mS_{b:}  ]^s  }^2 \leq \onormSize{}{\mX_{i:} \mV^{(t)}}^{-1} F_{ib}^{(t)} + G_{ib}^{(t)} + H_{ib}^{(t)} }. \label{eq:ind_b}
    \end{align}
    \normalsize
   Note that 
    \begin{align}
     \onormSize{}{\mX_{i:} \mV^{(t)}}  &= \theta(i) \onormSize{}{\mS_{i:} (\mTheta \mM)^{\otimes (K-1),T} \mW^{(t), \otimes^{K-1}}}\\
         &\geq \theta(i) \onormSize{}{\mS_{z(i):}} \lambda_r^{K-1}(\mTheta \mM) \lambda^{K-1}_r(\mW^{(t)})\\
         &\geq\theta(i) m , \label{eq:xv_lower}
    \end{align}
    where the first inequality follows from the property of eigenvalues; the last inequality follows from Lemma~\ref{lem:singular_thetam}, Lemma~\ref{lem:membership}, and assumption that $\min_{a \in [r]} \onormSize{}{\mS_{z(i):}} \geq c_3>0$; and $m >0$ is a positive constant related to $ c_3$. Plugging the lower bound of $ \onormSize{}{\mX_{i:} \mV^{(t)}}$ \eqref{eq:xv_lower} into the inequality \eqref{eq:ind_b} gives
    \begin{equation}\label{eq:decomp1}
        \ind \offf{ z^{(t+1)}(i) = b  } \leq A_{ib} + B_{ib},
    \end{equation}
    where 
    \begin{align}
        A_{ib} &= \ind \Bigg\{z^{(t+1)}(i) = b, \ang{ \mE_{i:} \mV, [  \tilde \mS_{z(i):} ]^s - [  \tilde \mS_{b:} ]^s }  \leq - \frac{\theta(i) m}{4} \onormSize{}{ [ \mS_{z(i):}  ]^s - [ \mS_{b:}  ]^s  }^2 \Bigg\},\\
        B_{ib} &= \ind \Bigg\{z^{(t+1)}(i) = b, \frac{1}{2}\onormSize{}{ [ \mS_{z(i):}  ]^s - [ \mS_{b:}  ]^s  }^2  \leq (\theta(i) m)^{-1} F_{ib}^{(t)} + G_{ib}^{(t)} + H_{ib}^{(t)} \Bigg\}.
    \end{align}
    \normalsize
 Taking the weighted summation of~\eqref{eq:decomp1} over $i \in [p]$ yields 
    \begin{equation}
        L^{(t+1)} \leq \xi + \frac{1}{p}\sum_{i \in [p] }  \sum_{b \in [r]/z(i)}  \zeta_{ib}^{(t)}, 
    \end{equation}
    where $\xi$ is the oracle loss such that 
    \begin{equation}\label{eq:decomp2}
        \xi =   \frac{1}{p} \sum_{i \in [p] } \theta(i) \sum_{b \in [r]/z(i)} A_{ib} \onormSize{}{ [ \mS_{ z(i):}  ]^s - [ \mS_{b:}  ]^s  }^2.
    \end{equation}
    Similarly to $\xi$ in~\eqref{eq:decomp2}, we define
    \begin{equation}
        \zeta_{ib}^{(t)} =   \theta(i) B_{ib} \onormSize{}{ [ \mS_{ z(i):}  ]^s - [ \mS_{b:}  ]^s  }^2.
    \end{equation}
    
    {

    \textbf{Now, we show the decomposition for MLE loss.}
     
    By the definition of Gaussian MLE, the estimator $\hat \mtheta$ satisfies $\hat \theta(i)  = \ang{ \mY_{i:} , \hat \mA_{\hat z(i):} }/{\onormSize{}{ \hat \mA_{\hat z(i):}}_F^2}$ for all $i \in [p]$.
Hence, we have 
\begin{equation}
    \hat z(i) = \argmin_{a \in [r_1]} \onormSize{}{ [\mY_{ i:}]^s - [\hat \mA_{a:}]^s}_F^2,
\end{equation}
and the decomposition 
 \begin{equation}
        L(\hat z) \leq \xi' + \frac{1}{p}\sum_{i \in [p] }  \sum_{b \in [r]/z(i)}  \zeta_{ib}', 
    \end{equation}
    where $ \zeta_{ib}' =   \theta(i)^2 B_{ib}' \onormSize{}{ [ \mA_{ z(i):}  ]^s - [ \mA_{b:}  ]^s  }^2$ and 
    \begin{align}
        A'_{ib} &= \ind \Bigg\{  \hat z(i)  = b,   \ang{\mE_{i:}, [\bar \mA_{z(i):}]^s - [\bar \mA_{b:}]^s }  \leq - \frac{m'}{4}  \sqrt{\frac{p^{K-1}}{r^{K-1}}} \onormSize{}{ [\mA_{z(i):}]^s -  [\mA_{b:}]^s  }_F^2 \Bigg\}, \\
    B'_{ib} &= \ind \Bigg\{ \hat z(i)  = b,     - \frac{1}{2}  \onormSize{}{ [\mA_{z(i):}]^s -  [\mA_{b:}]^s  }_F^2   \leq \sqrt{\frac{r^{K-1}}{(m')^2 p^{K-1}}}  \hat F_{ib} + \hat G_{ib} + \hat H_{ib} \Bigg\}
    \end{align}
    \normalsize
    with terms 
    \begin{align}
    \hat F_{ib} &= 2 \ang{ \mE_{i:}, ([\bar  \mA_{z(i):}]^s - [\hat \mA_{a:}]^s) - ([\bar  \mA_{b:}]^s - [\hat \mA_{b:}]^s) }, \\
    \hat G_{ib}& = \of{ \onormSize{}{\mX_{i:}^s - [\hat \mA_{z(i):}]^s}_F^2 - \onormSize{}{ \mX^s_{i:} - [\mP_{:z(i)}^T \mY \hat \mQ \hat \mB^T ]^s }_F^2  }  - \of{ \onormSize{}{\mX^s_{i:} - [\hat \mA_{b:}]^s}_F^2 - \onormSize{}{ \mX^s_{i:} - [\mP_{:b}^T \mY \hat \mQ \hat \mB^T]^s }_F^2   }, \\
   \hat  H_{ib}& = \onormSize{}{ \mX^s_{i:} - [\mP_{:z(i)}^T \mY \hat \mQ \hat \mB^T]^s }_F^2 - \onormSize{}{ \mX^s_{i:} - [\mP_{:b}^T \mY \hat \mQ \hat \mB^T]^s }_F^2  + \onormSize{}{ \mA^s_{z(i):} -  \mA^s_{b:}  }_F^2.
\end{align}
    
    }

    {\bf Step 3 (Derivation of contraction inequality).} In this step we derive the upper bound of $\zeta_{ib}$ and obtain the contraction inequality~\eqref{eq:proof_5_ineq}. 

    Choose the constant $\tilde C$ in the condition SNR $\geq \tilde C p^{-K/2} \log p$ that satisfies the condition of Lemma~\ref{lem:upper_fgh}, inequalities~\eqref{eq:tilde_c1}, and \eqref{eq:tilde_c2}. Note that  \begin{align}
        \zeta_{ib}^{(t)} &= \theta(i) \onormSize{}{ [ \mS_{ z(i):}  ]^s - [ \mS_{b:}  ]^s  }^2 \ind \offf{z^{(t+1)}(i) = b, \frac{1}{2}\onormSize{}{ [ \mS_{z(i):}  ]^s - [ \mS_{b:}  ]^s  }^2 \leq (\theta(i) m)^{-1} F_{ib}^{(t)} + G_{ib}^{(t)} + H_{ib}^{(t)} }\\
        & \leq \theta(i) \onormSize{}{ [ \mS_{ z(i):}  ]^s - [ \mS_{b:}  ]^s  }^2 \ind \offf{z^{(t+1)}(i) = b, \frac{1}{4}\onormSize{}{ [ \mS_{z(i):}  ]^s - [ \mS_{b:}  ]^s  }^2 \leq (\theta(i) m)^{-1} F_{ib}^{(t)} + G_{ib}^{(t)} } \\
        & \leq 64 \ind \offf{ z^{(t+1)}(i) = b} \of{  \frac{(F_{ib}^{(t)})^2}{ cm^2\onormSize{}{ [ \mS_{z(i):}  ]^s - [ \mS_{b:}  ]^s  }^2} + \frac{\theta(i)(G_{ib}^{(t)})^2}{ \onormSize{}{ [ \mS_{z(i):}  ]^s - [ \mS_{b:}  ]^s  }^2}    }
    \end{align}
    where the first inequality follows from the inequality~\eqref{eq:hib} in Lemma~\ref{lem:upper_fgh}, and the last inequality follows from the assumption that $\min_{i \in [p]} \theta(i) \geq c>0$. Following \citet[Step 4, Proof of Theorem 2]{han2022exact} and Lemma~\ref{lem:upper_fgh}, we have 
\begin{equation}\label{eq:f_sum}
          \frac{1}{p} \sum_{i \in [p]}\sum_{b \in [r]/z(i)} \ind \offf{ z^{(t+1)}(i) = b} \frac{(F_{ib}^{(t)})^2}{ cm^2 \onormSize{}{ [ \mS_{z(i):}  ]^s - [ \mS_{b:}  ]^s  }^2} \leq \frac{C_0 \bar C}{c m^2 \tilde C^2} L^{(t)},
    \end{equation}
    for a positive universal constant $C$ and 
    \small
    \begin{align}
         \frac{1}{p} \sum_{i \in [p]}\sum_{b \in [r]/z(i)} \ind \offf{ z^{(t+1)}(i) = b} \frac{\theta(i)(G_{ib}^{(t)})^2}{ \onormSize{}{ [ \mS_{z(i):}  ]^s - [ \mS_{b:}  ]^s  }^2}  &\leq \frac{1}{512} \frac{1}{p} \sum_{i \in [p]}\theta(i) \sum_{b \in [r]/z(i)}   \ind \offf{ z^{(t+1)}(i) = b}  (\Delta_{\min}^2 + L^{(t)}) \\
         & \leq \frac{1}{512} (L^{(t+1)} + L^{(t)}), \label{eq:g_sum}
    \end{align}
    \normalsize
 where the last inequality follows from the definition of $L^{(t)}$ and the constraint of $\mtheta$ in parameter space~\eqref{eq:family}. For $\tilde C$ also satisfies 
    \begin{equation}\label{eq:tilde_c3}
         \frac{C_0 \bar C}{c m^2 \tilde C^2} \leq \frac{1}{512},
    \end{equation}
    we have 
    \begin{equation}\label{eq:zeta_upper}
        \frac{1}{p}\sum_{i \in [p] }  \sum_{b \in [r]/z(i)}  \zeta_{ib}^{(t)} \leq  \frac{1}{8} L^{(t+1)}  + \frac{1}{4} L^{(t)}.
    \end{equation}
    Plugging the inequality~\eqref{eq:zeta_upper} into the decomposition~\eqref{eq:decomp2}, we obtain the contraction inequality 
    \begin{equation}\label{eq:decomp_t}
          L^{(t+1)} \leq \frac{3}{2} \xi  + \frac{1}{2} L^{(t)},
    \end{equation}
     where $\frac{1}{2}$ is the contraction parameter. 
     
     Therefore, with $\tilde C$ satisfying inequalities~\eqref{eq:tilde_c3}, \eqref{eq:tilde_c1} and \eqref{eq:tilde_c2}, we obtain the conclusion in Theorem~\ref{thm:refinement} via inequality~\eqref{eq:decomp_t} combining the inequality~\eqref{eq:cond_oracle} in Condition~\ref{cond:origin} and Lemma~\ref{lem:mis}. 
     
     \textbf{We also have the contraction inequality for MLE.}
     
     Following the same derivation of~\eqref{eq:decomp_t} with the upper bound of $\hat F_{ib}, \hat G_{ib}, \hat H_{ib}$ in Lemma~\ref{lem:intermediate}, we also have
    \begin{equation}\label{eq:decomp_t_mle}
          L(\hat z) \leq \frac{3}{2} \xi'  + \frac{1}{2} L(\hat z),
    \end{equation}
    which indicates the conclusion $\ell(\hat z,z) \lesssim \Delta_{\min}^2  \exp \of{ - \frac{p^{K-1}}{r^{K-1}} \Delta_{\min}^2  }$.

    {\bf Step 4 (Verification of Condition~\ref{cond:origin}).} Last, we verify the Condition~\ref{cond:origin} under high probability to finish the proof. Note that the inequalities~\eqref{eq:cond1}, \eqref{eq:cond2}, and \eqref{eq:cond3} describe the property of the sub-Gaussian noise tensor $\tE$, and the readers can find the proof directly in \citet[Step 5, Proof of Theorem 2]{han2022exact}. The initial condition~\eqref{eq:cond_intial} for MLE is satisfied by Lemma~\ref{lem:poly_mle_degree}. Here, we include only the verification of inequalities~\eqref{eq:cond_oracle} and \eqref{eq:cond_intial} for algorithm estimators. 
    
    Now, we verify the oracle loss condition~\eqref{eq:cond_oracle}. Recall the definition of $\xi$,
    \begin{align}
        \xi & = \frac{1}{p} \sum_{i \in [p]} \theta(i) \sum_{b \in [r]} \ind \Big\{\ang{ \mE_{i:} \mV, [  \tilde \mS_{z(i):} ]^s - [  \tilde \mS_{b:} ]^s}  \leq - \frac{ \theta(i) m}{4} \onormSize{}{ [ \mS_{z(i):}  ]^s - [ \mS_{b:}  ]^s  }^2 \Big\}  \cdot \onormSize{}{ [ \mS_{z(i):}  ]^s - [ \mS_{ b:}  ]^s  }^2.
    \end{align}
    Let $e_i = \mE_{i:} \mV$ denote the aggregated noise vector for all $i \in [p]$, and $e_i$'s are independent zero-mean sub-Gaussian vector in $\bbR^{r^{K-1}}$. The entries in $e_i$ are independent zero-mean sub-Gaussian variables with sub-Gaussian norm upper bounded by $m_1\sqrt{r^{K-1}/p^{K-1}}$ with some positive constant $m_1$. We have the probability inequality
    \begin{align}
        &\bbP \of{ \ang{ e_i, [  \tilde \mS_{z(i):} ]^s - [  \tilde \mS_{b:} ]^s } \leq - \frac{ \theta(i) m}{4} \onormSize{}{ [ \mS_{z(i):}  ]^s - [ \mS_{b:}  ]^s  }^2}\leq P_1 + P_2 + P_3,
    \end{align}
    where 
    \begin{align}
        P_1 &= \bbP \of{ \ang{e_i, [  \mS_{z(i):} ]^s - [ \mS_{b:} ]^s} \leq -\frac{\theta(i)m}{8}  \onormSize{}{ [ \mS_{z(i):}  ]^s - [ \mS_{b:}  ]^s  }^2 }, \\
        P_2 &= \bbP \of{ \ang{e_i, [ \tilde  \mS_{z(i):} ]^s - [ \mS_{z(i):} ]^s} \leq -\frac{\theta(i)m}{16}  \onormSize{}{ [ \mS_{z(i):}  ]^s - [ \mS_{b:}  ]^s  }^2 },\\
        P_3 &= \bbP \of{ \ang{e_i, [  \mS_{b:} ]^s - [\tilde  \mS_{b:} ]^s} \leq -\frac{\theta(i)m}{16}   \onormSize{}{ [ \mS_{z(i):}  ]^s - [ \mS_{b:}  ]^s  }^2 }.
    \end{align}
    \normalsize
    For $P_1$, notice that the inner product $\ang{e_j,\mS_{z(j):}^s - \mS_{b:}^s} $ is a sub-Gaussian variable with sub-Gaussian norm bounded by $m_2 \sqrt{ r^{K-1}/ p^{K-1} } \onormSize{}{ \mS_{z(i):}^s -  \mS_{b:}^s }$ with some positive constant $m_2$. Then, by Chernoff bound, we have  
    \begin{equation}\label{eq:p1_2}
         P_1 \lesssim \exp \of{  - \frac{p^{K-1}}{r^{K-1}}   \onormSize{}{ [ \mS_{z(j):}  ]^s - [ \mS_{b:}  ]^s  }^2 }.
    \end{equation}
    
    For $P_2$ and $P_3$, we only need to derive the upper bound of $P_2$ due to the symmetry. By the law of total probability, we have 
    \begin{equation}\label{eq:p2}
        P_2 \leq P_{21} + P_{22},
    \end{equation}
    where with some positive constant $t>0$,
    \begin{align}
        P_{21}& =  \bbP \of{ t \leq   \onormSize{}{ [ \tilde  \mS_{z(i):} ]^s - [ \mS_{z(i):} ]^s }},\\
        P_{22} &= \bbP \Bigg(  \ang{e_i, [ \tilde  \mS_{z(i):} ]^s - [ \mS_{z(i):} ]^s} \leq -\frac{\theta(i)m}{16} \cdot \onormSize{}{ [ \mS_{z(i):}  ]^s - [ \mS_{b:}  ]^s  }^2   \bigg| \onormSize{}{ [ \tilde  \mS_{z(i):} ]^s - [ \mS_{z(i):} ]^s } < t  \Bigg).
    \end{align}
    
    For $P_{21}$, note that the term $\mW^T_{: z(i)} \mE \mV = \frac{ \sum_{ j \neq i, j \in [p]} \ind \{ z(j) = z(i) \} e_j }{\sum_{j \in [p]} \ind \{ z(j) = z(i) \} }$  is a sub-Gaussian vector with sub-Gaussian norm bounded by $m_3 \sqrt{r^K/ p^K}$ with some positive constant $m_3$. This implies
        \begin{align}
        P_{21} & \leq \bbP\of{ t \onormSize{}{\mS_{z(i):} } \leq \onormSize{}{ \tilde  \mS_{z(i):}  -  \mS_{z(i):} }  } \bbP\of{c_3 t \leq  \onormSize{}{ \mW^T_{: z(i)} \mE \mV  } }\lesssim \exp\of{ - \frac{p^K t^2}{r^K} }, \label{eq:p21}
    \end{align}
    where the first inequality follows from the basic inequality in Lemma~\ref{lem:norm_diff}, the second inequality follows from the assumption that $ \min_{a \in [r]} \onormSize{}{\mS_{z(i):}} \geq c_3>0$ in~\eqref{eq:family}, and the last inequality follows from the Bernstein inequality.
    
    For $P_{22}$, the inner product $ \ang{e_i, [ \tilde  \mS_{z(i):} ]^s - [ \mS_{z(i):} ]^s}$ is also a sub-Gaussian variable with sub-Gaussian norm $ m_4 \sqrt{ r^{K-1}/ p^{K-1} } t$, conditioned on $\onormSize{}{ [ \tilde  \mS_{z(i):} ]^s - [ \mS_{z(i):} ]^s } < t$ with some positive constant $m_4$. Then, by Chernoff bound, we have 
    \begin{equation}\label{eq:p22}
        P_{22} \lesssim \exp \of{  - \frac{p^{K-1}}{r^{K-1} t^2}   \onormSize{}{ [ \mS_{z(j):}  ]^s - [ \mS_{b:}  ]^s  }^4 }. 
     \end{equation}
    We take $t =  \onormSize{}{ [ \mS_{z(i):}  ]^s - [ \mS_{b:}  ]^s  }$ in $P_{21}$ and $P_{22}$, and plug the inequalities \eqref{eq:p21} and \eqref{eq:p22} into to the upper bound for $P_2$ in~\eqref{eq:p2}. We obtain that 
    \begin{equation}\label{eq:p2_2}
        P_2 \lesssim \exp \of{  - \frac{p^{K-1}}{r^{K-1}}   \onormSize{}{ [ \mS_{z(i):}  ]^s - [ \mS_{b:}  ]^s  }^2 }.
    \end{equation}
    Combining the upper bounds~\eqref{eq:p1_2} and \eqref{eq:p2_2} gives 
    \begin{align}
        &\bbP \of{ \ang{ e_i, [  \tilde \mS_{z(i):} ]^s - [  \tilde \mS_{b:} ]^s }  \leq - \frac{\theta(i) m}{4} \onormSize{}{ [ \mS_{z(i):}  ]^s - [ \mS_{b:}  ]^s  }^2}  \lesssim  \exp \of{  - \frac{p^{K-1}}{r^{K-1}}   \onormSize{}{ [ \mS_{z(i):}  ]^s - [ \mS_{b:}  ]^s  }^2 }.\label{eq:xi_p}
    \end{align}

    Hence, we have 
    \begin{align}
        \bbE \xi & =  \frac{1}{p} \sum_{i \in [p]} \theta(i) \sum_{b \in [r]} \bbP \Bigg\{  \ang{ \mE_{i:} \mV, [  \tilde \mS_{z(i):} ]^s - [  \tilde \mS_{b:} ]^s }   \leq - \frac{\theta(i) m}{4} \onormSize{}{ [ \mS_{z(i):}  ]^s - [ \mS_{b:}  ]^s  }^2 \Bigg\} \onormSize{}{ [ \mS_{z(i):}  ]^s - [ \mS_{ b:}  ]^s  }^2 \\
        &\lesssim \frac{1}{p} \sum_{i \in [p]}  \theta(i) \max_{i \in [p], b \in [r]}  \onormSize{}{ [ \mS_{z(i):}  ]^s - [ \mS_{b:}  ]^s  }^2  \cdot \exp \of{  - \frac{p^{K-1}}{r^{K-1}}   \onormSize{}{ [ \mS_{z(i):}  ]^s - [ \mS_{b:}  ]^s  }^2 } \\
        & \leq \exp \of{  - M \frac{p^{K-1}}{r^{K-1}}   \Delta_{\min}^2 }, \label{eq:constant_M}
    \end{align}
    where $M$ is a positive constant, the first inequality follows from the constraint that $\sum_{i \in [p]} \theta(i) = p$, and the last inequality follows from~\eqref{eq:xi_p}.
    
    By Markov's inequality, we have 
    \begin{align}
        & \bbP\of{ \xi \lesssim \bbE \xi + \exp \of{  - \frac{M p^{K-1}}{2 r^{K-1}}   \Delta_{\min}^2 } }  \geq 1 -  C \exp \of{  - \frac{M p^{K-1}}{2 r^{K-1}}   \Delta_{\min}^2 },
    \end{align}
    and thus the condition~\eqref{eq:cond_oracle} holds with probability at least $1 -  C \exp \of{  - \frac{M p^{K-1}}{2 r^{K-1}}   \Delta_{\min}^2 }$ for some constant $C > 0$.
    
    {

    \textbf{The initialization condition for MLE also holds.}
    
    For $\xi'$, notice that $\ang{\mE_i, \mA^s_{a:} - \mA^s_{b:} }$ is a sub-Gaussian vector with variance bounded by $\onormSize{}{\mA^s_{a:} - \mA^s_{b:}}^2$ and 
    \begin{align}
        \bbP \of{ t \leq \onormSize{}{ [\bar \mA_{a:}]^s - \mA^s_{a:} }} & \leq \of{ t \leq \onormSize{}{ [\mP_{:a}^T \mY \mQ]^s - [\mP_{:a}^T \mX \mQ]^s }} \\
        &\leq \bbP( t \min_{a \in [r]}\onormSize{}{\mS_{a:}} \leq \onormSize{}{\mP_{:a}^T \mE \mQ} ) \\
        &\lesssim \exp \of{ - \frac{p^{K}t^2}{r^{K}}},
    \end{align}
   where the first inequality follows from the property in later inequality \eqref{eq:angle_aq}. We also have 
   \begin{equation}
        \xi' \lesssim  \of{  -  \frac{p^{K-1}}{r^{K-1}}   \Delta_{\min}^2 }.
   \end{equation}
    
    }
    
    Finally, we verify the bounded loss condition~\eqref{eq:cond_intial} for algorithm estimator by induction.  With output $z^{(0)}$ from Sub-algorithm~\hyperref[alg:main]{2} and the assumption SNR $\geq \tilde C p^{-K/2} \log p$, by Theorem~\ref{thm:initial}, we have 
    \begin{equation}
        L^{(0)} \leq \frac{\bar C \Delta_{\min}^2}{\tilde C r \log p},\quad \text{    when $p$ is large enough}.
    \end{equation}
Therefore, the condition~\eqref{eq:cond_intial} holds for $t = 0$. Assume that the condition~\eqref{eq:cond_intial} also holds for all $t \leq t_0$. Then, by the decomposition~\eqref{eq:decomp_t}, we have 
    \begin{align}
         L^{(t_0+1)} &\leq \frac{3}{2} \xi + \frac{1}{2} L^{(t_0)} \\
         & \leq \exp \of{  - M \frac{p^{K-1}}{r^{K-1}}   \Delta_{\min}^2 } + \frac{\Delta_{\min}^2}{r \log p} \\
         & \leq \frac{\bar C}{\tilde C}\frac{\Delta_{\min}^2}{r\log p} ,
    \end{align}
    where the second inequality follows from the condition~\eqref{eq:cond_oracle} and the last inequality follows from the assumption that $\Delta_{\min}^2 \gtrsim p^{-K/2} \log p$. Thus, the condition~\eqref{eq:cond_intial} holds for $t_0 + 1$, and the condition~\eqref{eq:cond_intial} is proved by induction.
\end{proof}

{\bf Useful Lemmas for the Proof of Theorem~\ref{thm:refinement}}

\begin{lem}[Singular-value property of membership matrices]\label{lem:membership} Under the setup of Theorem~\ref{thm:refinement}, suppose that the condition~\eqref{eq:cond_intial} holds. Then, for all $a \in [r]$, we have $|\of{z^{(t)}}^{-1}(a)| \asymp p/r$. Moreover, we have 
\begin{align}
     &\lambda_{r}(\mM) \asymp \onorm{\mM}_{\sigma} \asymp \sqrt{p/r}, \  \lambda_{r}(\mW) \asymp \onorm{\mW}_{\sigma} \asymp \sqrt{r/p},\\
     &\lambda_{r}(\mP) \asymp \onorm{\mP}_{\sigma}    \asymp  \min_{a \in [r] } \onormSize{}{\mtheta_{z^{-1}(a)}}^{-1}  \lesssim  \sqrt{r/p}. \label{eq:membership1}
\end{align}
The inequalities~\eqref{eq:membership1} also hold by replacing $\mM$ and $\mW$ to $\mM^{(t)}$ and $\mW^{(t)}$ respectively. 
Further, we have 
\begin{equation}\label{eq:WWsingular}
   \lambda_{r}(\mW\mW^T) \asymp \onorm{\mW\mW^T}_{\sigma} \asymp r/p,
\end{equation}
which is also true for $\mW^{(t)}\mW^{(t),T}$.
\end{lem}

\begin{proof}[Proof of Lemma~\ref{lem:membership}] The proof for the inequality~\eqref{eq:membership1} for $\mM, \mW$ can be found in \citet[Proof of Lemma 4]{han2022exact}. The inequalities for $\mP$ follows the same derivation with balance assumption on $\mtheta$ and $\min_{i \in [p]}\theta(i) \geq c$.

For inequality~\eqref{eq:WWsingular}, note that for all $k \in [r]$,
\begin{align}
    \lambda_k(\mW\mW^T) &= \sqrt{\text{eigen}_k (\mW \mW^T\mW\mW^T)} \asymp \sqrt{\frac{r}{p} \text{eigen}_k (\mW \mW^T)} = \sqrt{\frac{r}{p} \lambda^2_k(\mW)} \asymp \frac{r}{p},
\end{align}
where $\text{eigen}_k(\mA)$ denotes the $k$-th largest eigenvalue of the square matrix $\mA$, the first inequality follows the fact that $\mW^T\mW$ is a diagonal matrix with elements of order $r/p$, and the second equation follows from the definition of singular value.
\end{proof}

\begin{lem}[Upper bound for $F_{ib}^{(t)}, G_{ib}^{(t)}$ and $H_{ib}^{(t)}$]\label{lem:upper_fgh} Under the Condition~\ref{cond:origin} and the setup of Theorem~\ref{thm:refinement}  {with fixed $r \geq 2$},  assume the constant $\tilde C$ in the condition SNR $\geq \tilde C p^{-K/2} \log p$ is large enough to satisfy the inequalities~\eqref{eq:tilde_c1} and \eqref{eq:tilde_c2}. {As $p \rightarrow \infty$}, we have 
\begin{align}
    &\max_{i \in [p]} \max_{b \neq z(i)} \frac{\of{F_{ib}^{(t)}}^2}{ \onormSize{}{ [  \mS_{z(i):} ]^s  -  [ \mS_{b:} ]^s }^2 } \lesssim  \frac{r L^{(t)}}{ \Delta_{\min}^2} \onormSize{}{\mE_{i:} \mV}^2 + \of{1 +  \frac{ r L^{(t)}}{ \Delta_{\min}^2}} \onormSize{}{\mE_{i:} (\mV - \mV^{(t)})}^2 ,\label{eq:fib}
\end{align}
    \begin{equation}\label{eq:gib}
        \max_{i \in [p]} \max_{b \neq z(i)} \frac{ \of{  G_{ib}^{(t)} }^2  }{ \onormSize{}{ [\mS_{z(i):} ]^s - [\mS_{b:}]^s }^2}  \leq  \frac{1}{512} \of{ \Delta_{\min}^2 +L^{(t)}},
    \end{equation}
    \begin{equation}\label{eq:hib}
        \max_{i \in [p]} \max_{b \neq z(i)}  \frac{ \aabs{  H_{ib}^{(t)} }  }{ \onormSize{}{ [\mS_{z(i):} ]^s - [\mS_{b:}]^s }^2} \leq \frac{1}{4}.
    \end{equation}
    \normalsize
    
    Similarly, when the SNR $\geq \tilde C p^{-(K-1)}\log p$ with a large constant $\tilde C$, we have 
     \begin{equation}\label{eq:fib_cor}
        \max_{i \in [p]} \max_{b \neq z(i)} \frac{\of{\hat F_{ib}}^2}{ \onormSize{}{ [  \mA_{z(i):} ]^s  -  [ \mA_{b:} ]^s }^2 } \lesssim p^{K-1} \frac{r L(\hat z)}{ \Delta_{\min}^2}
    \end{equation}
    \begin{equation}\label{eq:gib_cor}
        \max_{i \in [p]} \max_{b \neq z(i)} \frac{ \of{  \hat G_{ib} }^2  }{ \onormSize{}{ [\mA_{z(i):} ]^s - [\mA_{b:}]^s }^2}  \leq  \frac{1}{512} \of{ \Delta_{\min}^2 +L(\hat z)},
    \end{equation}
    \begin{equation}\label{eq:hib_cor}
        \max_{i \in [p]} \max_{b \neq z(i)}  \frac{ \aabs{  \hat H_{ib} }  }{ \onormSize{}{ [\mA_{z(i):} ]^s - [\mA_{b:}]^s }^2} \leq \frac{1}{4}.
    \end{equation}
    \end{lem}
    \normalsize

\begin{proof}[Proof of Lemma~\ref{lem:upper_fgh}]
We prove the the first three inequalities in Lemma~\ref{lem:upper_fgh} separately.
\begin{enumerate}[wide]
    \item Upper bound for $F_{ib}^{(t)}$, i.e., inequality \eqref{eq:fib}. Recall the definition of $F_{ib}^{(t)}$,
    \begin{align}
         F_{ib}^{(t)} &= 2 \ang{\mE_{i:} \mV^{(t)}, \of{ [ \tilde \mS_{z(i):} ]^s  - [\mS_{z(i):}^{(t)} ]^s  }  - \of{  [  \tilde \mS_{b:} ]^s  - [ \mS_{b:}^{(t)} ]^s  }  }  + 2 \ang{ \mE_{i:}(\mV - \mV^{(t)} ), [ \tilde \mS_{z(i):} ]^s - [ \tilde \mS_{b:} ]^s  }.
    \end{align}
    \normalsize
    By Cauchy-Schwartz inequality, we have 
    \begin{align}
        \of{F_{ib}^{(t)}}^2 
        & \leq 8 \of{\ang{\mE_{i:} \mV^{(t)}, \of{ [  \tilde \mS_{z(i):} ]^s  -  [  \mS_{z(i):}^{(t)} ]^s  }  - \of{  [  \tilde \mS_{b:} ]^s  -  [  \mS_{b:}^{(t)} ]^s  }  }  }^2 \\
         & \quad \quad \quad \quad \quad + 8\of{ \ang{ \mE_{i:} (\mV - \mV^{(t)} ), [  \tilde \mS_{z(i):} ]^s - [ \tilde  \mS_{b:} ]^s  }}^2 \\
         & \leq 8 \of{ \onormSize{}{\mE_{i:} \mV}^2 + \onormSize{}{\mE_{i:} ( \mV - \mV^{(t)} )}^2 } \max_{a \in [r]^s} \onormSize{}{  [  \tilde \mS_{a:} ]^s - [   \mS_{a:}^{(t)} ]^s  } \\
          & \quad \quad \quad \quad \quad +\onormSize{}{\mE_{i:} ( \mV - \mV^{(t)} )}^2  \onormSize{}{  [  \tilde \mS_{z(i):} ]^s - [  \tilde \mS_{b:} ]^s  }. \label{eq:fib_decomp}
    \end{align}
    \normalsize
    Note that for all $a \in [r]$,
    \begin{align}
        \onormSize{}{  [  \tilde \mS_{a:} ]^s - [   \mS_{a:}^{(t)} ]^s  }^2 &= \onormSize{}{ [ \mW_{:a}^T \mY \mV ]^s - [ \mW_{:a}^{(t),T} \mY \mV^{(t)} ]^s }^2 \\
        &\leq 2 \onormSize{}{ [ \mW_{:a}^T \mY \mV ]^s  - [ \mW_{:a}^{(t),T} \mY \mV ]^s  }^2  + 2 \onormSize{}{[ \mW_{:a}^{(t),T} \mY \mV ]^s - [ \mW_{:a}^{(t),T} \mY \mV^{(t)} ]^s }^2\\
        & \lesssim \frac{ r^2 (L^{(t)})^2}{ \Delta_{\min}^2} + \frac{ r r^{2K} + p r^{K+2}}{p^K} \frac{L^{(t)}}{\Delta_{\min}^2}\\
        & \lesssim  r L^{(t)}+  \frac{ r r^{2K} + p r^{K+2}}{p^K} \frac{L^{(t)}}{\Delta_{\min}^2}\\
        & \lesssim  r L^{(t)},\label{eq:fib_1}
    \end{align}
    where the second inequality follows from the inequalities~\eqref{eq:inter3} and \eqref{eq:inter4} in Lemma~\ref{lem:intermediate}, the third inequality follows from the condition~\eqref{eq:cond_intial} in Condition~\ref{cond:origin}, and the last inequality follows from the assumption that $\Delta_{\min}^2 \geq \tilde C p^{-K/2}\log p$. 
    
    Note that 
    \begin{align}
        \onormSize{}{  [  \tilde \mS_{z(i):} ]^s - [  \tilde \mS_{b:} ]^s  }^2
        &=  \onormSize{}{  [  \tilde \mS_{z(i):} ]^s -[ \mS_{z(i):} ]^s +  [ \mS_{z(i):} ]^s - [ \mS_{b:} ]^s + [ \mS_{b:} ]^s -[  \tilde \mS_{b:} ]^s  }^2 \\
        & \lesssim \onormSize{}{ [ \mS_{z(i):} ]^s - [ \mS_{b:} ]^s }^2 + \max_{a \in [r]} \onormSize{}{ [ \mS_{a:} ]^s -[  \tilde \mS_{a:} ]^s }^2 \\
        & \lesssim \onormSize{}{ [ \mS_{z(i):} ]^s - [ \mS_{b:} ]^s }^2  + \max_{a \in [r]} \frac{1}{\onormSize{}{\mS_{a:} }^2} \onormSize{}{ \mW_{:a}^T \mE \mV}^2\\
        & \lesssim \onormSize{}{ [ \mS_{z(i):} ]^s - [ \mS_{b:} ]^s }^2, \label{eq:fib_2}
    \end{align}
    \normalsize
    where the second inequality follows from Lemma~\ref{lem:norm_diff}, and the last inequality follows from the assumptions on $\onormSize{}{\mS_{a:}}$ in the parameter space~\eqref{eq:family}, the inequality~\eqref{eq:cond1} in Condition~\ref{cond:origin} and the assumption $\Delta_{\min}^2 \gtrsim p^{-K/2}\log p$. 
    
    Therefore, we finish the proof of inequality~\eqref{eq:fib} by plugging the inequalities~\eqref{eq:fib_1} and \eqref{eq:fib_2} into the upper bound~\eqref{eq:fib_decomp}.
    
    \item Upper bound for $G_{ib}^{(t)}$, i.e., inequality~\eqref{eq:gib}. By definition of $G_{ib}^{(t)}$, we rearrange terms and obtain
    \begin{align}
        G_{ib}^{(t)} &=   \of{ \onormSize{}{ [\mX_{i:} \mV^{(t)}]^s -  [  \mS_{z(i):}^{(t)}  ]^s}^2 -  \onormSize{}{ [\mX_{i:} \mV^{(t)}]^s -  [  \mW_{:z(i)}^T \mY \mV^{(t)} ]^s}^2}   \\
        & \quad \quad \quad \quad \quad -   \of{ \onormSize{}{ [\mX_{i:} \mV^{(t)}]^s -  [  \mS_{b:}^{(t)}  ]^s}^2 -  \onormSize{}{ [\mX_{i:} \mV^{(t)}]^s -  [  \mW_{:b}^T \mY \mV^{(t)} ]^s}^2 } \\
        &= 2 \ang{  [\mX_{i:} \mV^{(t)}]^s , \of{[  \mW_{:z(i)}^T \mY \mV^{(t)} ]^s -  [  \mS_{z(i):}^{(t)}  ]^s } - \of{ [  \mW_{:b}^T \mY \mV^{(t)} ]^s -  [  \mS_{b:}^{(t)}  ]^s }}\\
      &= G_1 + G_2 - G_3, \label{eq:gib_decomp}
    \end{align}
    where 
    \begin{align}
        G_1 &= \onormSize{}{ [  \mW_{:z(i)}^T \mY \mV^{(t)} ]^s -  [  \mS_{z(i):}^{(t)}  ]^s}^2 - \onormSize{}{ [  \mW_{:b}^T \mY \mV^{(t)} ]^s -  [  \mS_{b:}^{(t)}  ]^s }^2,\\
        G_2 &= 2 \ang{ [\mX_{i:} \mV^{(t)}]^s  -  [  \mW_{:z(i)}^T \mY \mV^{(t)} ]^s,   [  \mW_{:z(i)}^T \mY \mV^{(t)} ]^s -  [  \mS_{z(i):}^{(t)}  ]^s},\\
        G_3 & = 2 \ang{ [\mX_{i:} \mV^{(t)}]^s  -  [  \mW_{:b}^T \mY \mV^{(t)} ]^s,  [  \mW_{:b}^T \mY \mV^{(t)} ]^s -  [  \mS_{b:}^{(t)}  ]^s}.
    \end{align}
    \normalsize
    
    For $G_1$, we have 
    \begin{align}
        |G_1|^2 &\leq \aabs{ \onormSize{}{ [ \mW_{:z(i)}^T \mY \mV^{(t)} ]^s - [ \mS_{z(i):}^{(t)}  ]^s}^2 - \onormSize{}{ [ \mW_{:b}^T \mY \mV^{(t)} ]^s -  [ \mS_{b:}^{(t)}  ]^s }^2 }^2 \\
        & \leq \max_{a \in [r]} \onormSize{}{ [ \mW_{:a}^T \mY \mV^{(t)} ]^s - [ \mW_{:a}^{(t),T} \mY \mV^{(t)}  ]^s}^4\\
        & \leq C^4 \frac{r^4}{\Delta_{\min}^4} (L^{(t)})^4 +  \frac{r^2 r^{4K} + p^2r^{2K+4}}{p^{2K}} \frac{(L^{(t)})^2}{\Delta_{\min}^4} \\
        & \leq C^4 \frac{\bar C}{\tilde C^3} \of{\Delta_{\min}^4 +  \Delta_{\min}^2 L^{(t)}}, \label{eq:g1}
    \end{align}
    where the third inequality follows from the inequality~\eqref{eq:inter5} in Lemma~\ref{lem:intermediate} and the last inequality follows from the assumption that $\Delta_{\min}^2 \geq \tilde C p^{-K/2}\log p$ and inequality~\eqref{eq:cond_intial} in Condition~\ref{cond:origin}.
    
    For $G_2$, noticing that $[\mX_{i:} \mV^{(t)}]^s = [ \mW_{z(i):}^T\mX \mV^{(t)}]^s$, we have 
    \begin{align}
        |G_2|^2 &\leq 2\onormSize{}{ [\mX_{i:} \mV^{(t)}]^s  -  [  \mW_{:z(i)}^T \mY \mV^{(t)} ]^s}^2 \onormSize{}{ [  \mW_{:z(i)}^T \mY \mV^{(t)} ]^s -  [  \mS_{z(i):}^{(t)}  ]^s }^2 \\
        & \leq \frac{2}{ \onormSize{}{ \mW_{z(i):}^T\mX \mV^{(t)}}^2} \max_{a \in [r]}\onormSize{}{ \mW_{:a}^T \mE \mV^{(t)} }^2 \max_{a \in [r]} \onormSize{}{ [  \mW_{:a}^T \mY \mV^{(t)} ]^s -  [   \mW_{:a}^{(t),T} \mY \mV^{(t)}  ]^s}^2\\
        & \leq C' \frac{ r^{2K-1} + K p r^{K+1} }{p^{K}} \of{ \frac{r^2}{\Delta_{\min}^2} (L^{(t)})^2 + \frac{r r^{2K} + pr^{K+2}}{p^{K}} \frac{L^{(t)}}{\Delta_{\min}^2} }\\
        & \leq \frac{C'}{\tilde C^2} \Delta_{\min}^2 L^{(t)},\label{eq:g2}
    \end{align}
    where $C'$ is a positive universal constant, the second inequality follows from Lemma~\ref{lem:norm_diff}, the third inequality follows from the inequality~\eqref{eq:cond2} in Condition~\ref{cond:origin}, the inequalities~\eqref{eq:inter5} and \eqref{eq:j11_wxvt} in the proof of Lemma~\ref{lem:intermediate},  and the last inequality follows from the assumption $\Delta_{\min}^2 \geq \tilde C p^{-K/2} \log p$ and inequality~\eqref{eq:cond_intial} in Condition~\ref{cond:origin}.
    
    For $G_3$, note that by triangle inequality
    \begin{align}
        \onormSize{}{[\mX_{i:} \mV^{(t)}]^s  -  [  \mW_{:b}^T \mX \mV^{(t)} ]^s}^2 & \leq \onormSize{}{\mS_{z(i):}^s - \mS_{b:}^s}^2 + 2\max_{a \in [r]}\onormSize{}{[\mW_{:a}^T\mX \mV^{(t)}]^s - [\mW_{:a}^T\mX \mV]^s }^2\\
        & \leq \onormSize{}{\mS_{z(i):}^s - \mS_{b:}^s}^2 + C \frac{r^2 (L^{(t)})^2}{\Delta_{\min}^2}, \label{eq:g3xv_wxv}
    \end{align}
    where the last inequality follows from the inequality~\eqref{eq:j1} in the proof of Lemma~\ref{lem:intermediate} and $C$ is a positive constant.
   Then we have 
   \begin{align}
        |G_3|^2 & \leq 2 \onormSize{}{[\mX_{i:} \mV^{(t)}]^s  -  [  \mW_{:b}^T \mY \mV^{(t)} ]^s}^2 \max_{a \in [r]} \onormSize{}{ [  \mW_{:a}^T \mY \mV^{(t)} ]^s -  [   \mW_{:a}^{(t),T} \mY \mV^{(t)}  ]^s}^2 \\
        & \leq 2 \of{ \onormSize{}{[\mX_{i:} \mV^{(t)}]^s  -  [  \mW_{:b}^T \mX \mV^{(t)} ]^s}^2  +  \onormSize{}{ [  \mW_{:b}^T \mY \mV^{(t)} ]^s  -  [  \mW_{:b}^T \mX \mV^{(t)} ]^s}^2  } \\
        & \quad \times \max_{a \in [r]} \onorm{ [  \mW_{:a}^T \mY \mV^{(t)} ]^s -  [   \mW_{:a}^{(t),T} \mY \mV^{(t)}  ]^s}^2\\
        & \leq C^2 \of{\onormSize{}{\mS_{z(i):}^s - \mS_{b:}^s}^2 + C \frac{r^2 (L^{(t)})^2}{\Delta_{\min}^2} } \of{ \frac{r^2(L^{(t)})^2}{ \Delta_{\min}^2} + \frac{r r^{2K} + pr^{K+2}}{p^{K}} \frac{L^{(t)}}{\Delta_{\min}^2}}+\frac{C'}{\tilde C^2} \Delta_{\min}^2 L^{(t)} \\
        & \leq \frac{C^2 \bar C^2}{\tilde C} \onormSize{}{\mS_{z(i):}^s - \mS_{b:}^s}^2 (\Delta_{\min}^2 + L^{(t)}) +  \frac{C^3 C' \bar C^2}{\tilde C^2 }\of{\Delta_{\min}^4 +  \Delta_{\min}^2 L^{(t)}},\label{eq:g3}
    \end{align}
    \normalsize
    where the third inequality follows from the same procedure to derive~\eqref{eq:g1} and \eqref{eq:g2}, and the last inequality follows from the assumption $\Delta_{\min}^2 \geq \tilde C p^{-K/2} \log p$ and inequality~\eqref{eq:cond_intial} in Condition~\ref{cond:origin}.
    
    Choose the $\tilde C$ such that
    \begin{equation}\label{eq:tilde_c1}
        3 \of{  C^4 \frac{\bar C}{\tilde C^3}  + \frac{C'}{\tilde C^2} + \frac{C^2 \bar C^2}{\tilde C} + \frac{C^3 C' \bar C^2}{\tilde C^2 }} \leq \frac{1}{512}.
    \end{equation}
    Then,  we finish the proof of inequality~\eqref{eq:gib} by plugging the inequalities~\eqref{eq:g1}, \eqref{eq:g2}, and \eqref{eq:g3} into the upper bound~\eqref{eq:gib_decomp}.  
    
    \item Upper bound for $H_{ib}^{(t)}$, i.e., the inequality~\eqref{eq:hib}. By definition of $H_{ib}$, we rearrange terms and obtain
    \begin{align}
        H_{ib} &= \onormSize{}{ [\mX_{i:} \mV^{(t)}]^s -  [  \mW_{:z(i)}^T \mY \mV^{(t)} ]^s}^2 - \onormSize{}{ [\mX_{i:} \mV^{(t)}]^s -  [  \mW_{:b}^T \mY \mV^{(t)} ]^s}^2 + \onormSize{}{ [ \mS_{z(i):}  ]^s - [ \mS_{b:}  ]^s  }^2 \\
        & =  \onormSize{}{ [ \mX_{i:} \mV^{(t)}]^s -  [  \mW_{:z(i)}^T \mY \mV^{(t)} ]^s}^2 \\
        & \quad \quad \quad \quad \quad + \of{ \onormSize{}{ [ \mS_{z(i):}  ]^s - [ \mS_{b:}  ]^s  }^2 - \onormSize{}{ [\mX_{i:} \mV^{(t)}]^s-  [  \mW_{:b}^T \mX \mV^{(t)} ]^s }  } \\
        & \quad \quad \quad \quad \quad  - \of{\onormSize{}{ [\mX_{i:} \mV^{(t)}]^s-  [  \mW_{:b}^T \mY \mV^{(t)} ]^s } -  \onormSize{}{ [\mX_{i:} \mV^{(t)}]^s-  [  \mW_{:b}^T \mX \mV^{(t)} ]^s }   } \\
        & = H_1 + H_2 + H_3,
    \end{align}
    \normalsize
    where 
    \begin{align}
        H_1 &= \onormSize{}{ [ \mX_{i:} \mV^{(t)}]^s -  [  \mW_{:z(i)}^T \mY \mV^{(t)} ]^s}^2 - \onormSize{}{ [ \mW_{:b}^T \mX \mV^{(t)}]^s -  [ \mW_{:b}^T \mY \mV^{(t)}]^s}^2, \\
        H_2 &=  \onormSize{}{ [ \mS_{z(i):}  ]^s - [ \mS_{b:}  ]^s  }^2 - \onormSize{}{ [\mX_{i:} \mV^{(t)}]^s-  [  \mW_{:b}^T \mX \mV^{(t)} ]^s }^2  , \\
        H_3 & = 2 \ang{  [\mX_{i:} \mV^{(t)}]^s  -  [  \mW_{:b}^T \mX \mV^{(t)} ]^s,  [  \mW_{:b}^T \mY \mV^{(t)} ]^s -  [  \mW_{:b}^T \mX \mV^{(t)} ]^s }.
    \end{align}
    \normalsize
    For $H_1$, we have 
     \begin{align}
        |H_1| \leq  \frac{4 \max_{a \in [r]}\onormSize{}{ \mW_{:a}^T \mE \mV^{(t)} }^2 }{ \onormSize{}{ \mW_{z(i):}^T\mX \mV^{(t)}}^2}
         \leq  \frac{ r^{2K-1} + K p r^{K+1} }{p^{K}} 
        \leq \tilde C^{-2} \onormSize{}{ [ \mS_{z(i):}  ]^s - [ \mS_{b:}  ]^s  }^2\label{eq:h1},
    \end{align}
    following the derivation of $G_2$ in inequality~\eqref{eq:g2} and the assumption that $\Delta_{\min}^2 \geq \tilde C p^{-K/2} \log p$.

    For $H_2$, by the inequality~\eqref{eq:g3xv_wxv}, we have 
    \begin{align}
        |H_2| &\lesssim 2\max_{a \in [r]}\onormSize{}{[\mW_{:a}^T\mX \mV^{(t)}]^s - [\mW_{:a}^T\mX \mV]^s }^2 \lesssim \frac{r^2 (L^{(t)})^2}{\Delta_{\min}^2} \leq C \frac{\bar C^2}{ \tilde C^2} \onormSize{}{ [ \mS_{z(i):} ]^s - [ \mS_{a:}]^s }^2, \label{eq:h2}
    \end{align}
    where the last inequality follows from the condition~\eqref{eq:cond_intial} in Condition~\ref{cond:origin}.
    
      For $H_3$,  by Cauchy-Schwartz inequality, we have 
      \begin{align}
           |H_3| & \lesssim \onormSize{}{[\mX_{i:} \mV^{(t)}]^s  -  [  \mW_{:b}^T \mX \mV^{(t)} ]^s} |H_1|^{1/2} \leq 2 \tilde C^{-1} \onormSize{}{ [ \mS_{z(i):} ]^s - [ \mS_{a:}]^s }^2,\label{eq:h3}
      \end{align}
    following the inequalities~\eqref{eq:g3xv_wxv} and \eqref{eq:h1}.
    
    Choose $\tilde C$ such that 
    \begin{equation}\label{eq:tilde_c2}
        \tilde C^{-2} + C \frac{\bar C^2}{ \tilde C^2} + \tilde C^{-1} \leq \frac{1}{4}.
    \end{equation}
     Therefore, we finish the proof of inequality~\eqref{eq:hib} combining inequalities~\eqref{eq:h1}, \eqref{eq:h2}, and \eqref{eq:h3}.

\end{enumerate}
{
  
\textbf{Next, we show the upper bounds for $\hat F_{ib}, \hat G_{ib}$ and $\hat H_{ib}$.}

By Lemma~\ref{lem:angle_gap_x}, we have 
\begin{equation}
    \onormSize{}{\mS_{a:}^s - \mS_{b:}^s}  = (1 + o(1)) \onormSize{}{\mA_{a:}^s - \mA_{b:}^s}.
\end{equation}
Also, notice that the matrix product of $ \mB^T$ corresponds to the padding operation in Lemma~\ref{lem:pad}, and the padding weights are balanced such that $\onormSize{}{\mv \mB}  = (1 + o(1)) \max_{a} \onormSize{}{\mtheta_{z^{-1}(a)}}^{(K-1)/2} \onormSize{}{\mv} $ for all $\mv \in \bbR^{r(K-1)}$. For two vectors $\mv_1, \mv_2 \in \bbR^{r^{K-1}}$, we have 
\begin{equation}\label{eq:b_angle}
    \onormSize{}{\mv_1^s - \mv_2^s} = (1 + o(1)) \onormSize{}{[\mv_1 \mB^T]^s -  [\mv_2 \mB^T]^s  }.
\end{equation}
The equation~\eqref{eq:b_angle} also holds for $\hat \mB^T$.

Note that for all $i \in [p]$ we have 
\begin{align}
    \onormSize{}{\mA_{i:} \hat \mQ} &= \onormSize{}{\mS_{z(i:)} \mB^T \hat \mQ }\\
    &=  \onormSize{}{\mS_{z(i:)} \hat \mD^{\otimes (K-1)} }\\
    &=  (1 + o(1)) \onormSize{}{\mS_{z(i:)}} \\
    &=  (1 + o(1)) \max_{a} \onormSize{}{\mtheta_{z^{-1}(a)}}^{-(K-1)/2} \onormSize{}{\mA_{i:}},  \label{eq:aq}
\end{align}
\normalsize
where the third inequality follows from the singular property of MLE confusion matrix~\eqref{eq:mle_confusion} and the last inequality follows from the fact that $\mA_i = \mS_{z(i:)} \mB^T$ and Lemma~\ref{lem:membership}. Above equation indicates that $\mA_{i:}$ is the span space of the singular values as $p \rightarrow \infty$. Also, notice that the row space of $\mP_{:a}^T \mY \hat \mQ  \hat \mB^T$ is equal to the column space of $\hat \mQ$, and $\mA_{i:} \neq \mP_{:a}^T \mY \hat \mQ  \hat \mB^T$ in noisy case. 

Hence, for all $a \in [r]$, we have 
\begin{align}
    \onormSize{}{ [\mX_i \hat \mQ]^s - [\mP_{:a}^T \mY \hat \mQ ]^s }
    &= \onorm{ \frac{\mA_{z(i:)} \hat \mQ  }{\onormSize{}{\mA_{z(i:)} \hat \mQ }} - \frac{\mP_{:a}^T \mY \hat \mQ  }{ \onormSize{}{\mP_{:a}^T \mY \hat \mQ} } } \\
    &= (1+o(1))  \onorm{ \frac{\mA_{z(i:)} }{\onormSize{}{\mA_{z(i:)} }} - \frac{\mP_{:a}^T \mY \hat \mQ \hat \mB^T  }{ \onormSize{}{\mP_{:a}^T \mY \hat \mQ \hat \mB^T} } } \\
    &= (1+o(1)) \onormSize{}{ [\mX_i]^s - [\mP_{:a}^T \mY \hat \mQ \hat \mB^T ]^s } \label{eq:angle_aq}
\end{align}
where the second equation follows from~\eqref{eq:aq}, $\onormSize{}{\mP_{:a}^T \mY \hat \mQ \hat \mB^T} = (1 + o(1)) \max_{a} \onormSize{}{\mtheta_{z^{-1}(a)}}^{(K-1)/2} \onormSize{}{\mP_{:a}^T \mY \hat \mQ } $, and singular property of $\hat \mB^T$. Similar result holds after replacing $\mP_{:a}^T \mY \hat \mQ$ by  $\mP_{:a}^T \mY\mQ$ or $\mP_{:a}^T \mY \hat \mQ$.

We are now ready to show the upper bounds for $\hat F_{ib},\hat G_{ib}$ and $\hat H_{ib}$.

 For $\hat F_{ib}$, we have 
 \begin{align}
    (\hat F_{ib})^2 &\leq \onormSize{}{\mE_{i:}}^2 \onormSize{}{ [\bar  \mA_{a:}]^s - [\hat \mA_{a:}]^s }^2\\
    &\leq \onormSize{}{\mE_{i:}}^2  \off{ \onormSize{}{ [ \bar \mS_{a:} \mB^T]^s - [\bar \mS_{a:} \hat \mB^T]^s } +  \onormSize{}{ [ \bar \mS_{a:} \hat \mB^T]^s -  [\hat \mS_{a:} \hat \mB^T]^s}  }^2 \\
    &\lesssim \onormSize{}{\mE_{i:}}^2  \off{ \onormSize{}{ [ \bar \mS_{a:} \mB^T \hat \mQ]^s - [\bar \mS_{a:}]^s } +  \onormSize{}{ [ \bar \mS_{a:} ]^s -  [ \hat \mS_{a:} ]^s}  }^2.
\end{align}
Following similar derivations in inequalities~\eqref{eq:fib_1}, \eqref{eq:fib_2}, and the upper bound for $J_{1}$ in the proof of Lemma~\ref{lem:intermediate}, respectively, we have 
\begin{equation}
    \onormSize{}{ [ \bar \mS_{a:} ]^s -  [ \hat \mS_{a:} ]^s}  \lesssim r L(\hat z), \quad  \onormSize{}{ [ \bar \mS_{a:} ]^s -  [ \bar \mS_{b:} ]^s} \lesssim \onormSize{}{\mS_{a:}^s - \mS_{b:}^s}^2, 
\end{equation}
and 
\begin{align}
    \onormSize{}{ [ \bar \mS_{a:} \mB^T \hat \mQ]^s - [\bar \mS_{a:}]^s } \lesssim L(\hat z).
\end{align}
We then obtain the upper bound for $\hat F_{ib}$ by noticing that $\onormSize{}{\mE_i}^2 \lesssim p^{K-1}$.

For $\hat G_{ib}$ and $\hat H_{ib}$, by the property~\eqref{eq:angle_aq}, we have 
\begin{align}
        (1+o(1))\hat G_{ib}& = \of{ \onormSize{}{[\mX_{i:}\hat \mQ]^s - [\hat \mS_{a:}]^s}_F^2 - \onormSize{}{ [\mX_{i:} \hat \mQ]^s - [\mP_{:a}^T \mY \hat \mQ ]^s }_F^2  } \\
        & \quad - \of{ \onormSize{}{[\mX_{i:} \hat \mQ]^s - [\hat \mS_{b:}]^s}_F^2 - \onormSize{}{ [\mX_{i:} \hat \mQ]^s - [\mP_{:b}^T \mY \hat \mQ ]^s }_F^2   }. \\
        (1+o(1))\hat H_{ib} & = \onormSize{}{ [\mX_{i:} \hat \mQ]^s - [\mP_{:a}^T \mY \hat \mQ ]^s }_F^2 - \onormSize{}{ [\mX_{i:} \hat \mQ]^s - [\mP_{:b}^T \mY \hat \mQ]^s }_F^2 + \onormSize{}{ \mA^s_{a:} -  \mA^s_{b:}  }_F^2.
    \end{align}
    We obtain the upper bounds following the proof for inequalities~\eqref{eq:gib} and \eqref{eq:hib}.
}

\end{proof}

     \begin{lem}[Relationship between misclustering loss and intermediate parameters]\label{lem:intermediate} Under the Condition~\ref{cond:origin} and the setup of Theorem~\ref{thm:refinement} {with fixed $r \geq 2$, as $p \rightarrow \infty$}, we have
    \begin{equation}\label{eq:inter1}
        \onormSize{}{\mV - \mV^{(t)}}_{\sigma} \lesssim \sqrt{\frac{r^{K-1}}{p^{K-1}}}  \frac{r}{\Delta_{\min}^2} L^{(t)},
    \end{equation}
    \begin{equation}\label{eq:inter2}
        \onormSize{}{\mE (\mV - \mV^{(t)})}_{\sigma} \lesssim \sqrt{ \frac{r^{K-1}(p r^{K-1} + p r)}{p^{K-1}}} 
 \frac{r}{\Delta_{\min}^2} L^{(t)},
 \end{equation}
 \begin{align}
     &\max_{b \in [r]} \onormSize{}{ [ \mW_{:b}^T \mY \mV ]^s   -  [ \mW_{:b}^{(t), T} \mY \mV  ]^s } \leq C \of{\frac{r L^{(t)}}{\Delta_{\min}} + \sqrt{ \frac{r^{2K} + pr^{K+1}}{p^{K}} } \frac{\sqrt{L^{(t)}}}{\Delta_{\min}}}, \label{eq:inter3}
 \end{align}
 \begin{align}
     &\max_{b \in [r]} \onormSize{}{  [  \mW_{:b}^{(t), T}  \mY \mV]^s - [  \mW_{:b}^{(t), T}  \mY \mV^{(t)}]^s }  \leq C  \of{ \sqrt{ \frac{rr^{2K} + p r^{K+2}}{p^{K}}}  \frac{\sqrt{L^{(t)}}}{\Delta_{\min}} +  \frac{r L^{(t)}}{\Delta_{\min}}}, \label{eq:inter4}
 \end{align}
 \begin{align}
     & \max_{b \in [r]} \onormSize{}{ [\mW_{:b}^T \mY \mV^{(t)} ]^s   -  [ \mW_{:b}^{(t), T} \mY \mV^{(t)} ]^s } \leq C \of{\frac{rL^{(t)}}{\Delta_{\min}} + \sqrt{ \frac{rr^{2K} + p r^{K+2}}{p^{K}}}  \frac{\sqrt{L^{(t)}}}{\Delta_{\min}}}, \label{eq:inter5}
 \end{align}
 \normalsize
    for some positive universal constant $C$. In addition, the inequality~\eqref{eq:inter4} also holds by replacing $\mW_{:b}^{(t)}$ to $\mW_{:b}$. Further, the above inequalities holds after replacing $\mW$ to $\mP$, $\mV$ to $\mQ$, and $L^{(t)}$ to $L(\hat z)$.
    \end{lem}

    \begin{proof}[Proof of Lemma~\ref{lem:intermediate}] We follow and use several intermediate conclusions in \citet[Proof of Lemma 5]{han2022exact}. We prove each inequality separately.
    \begin{enumerate}[wide]
    \item Inequality~\eqref{eq:inter1}. By \citet[Proof of Lemma 5]{han2022exact}, we have 
    \begin{equation}
         \onormSize{}{\mV - \mV^{(t)}}_{\sigma} \lesssim \sqrt{\frac{r^{K-1}}{p^{K-1}}} r \ell^{(t)}.
    \end{equation}
    Then, we complete the proof of inequality~\eqref{eq:inter1} by applying Lemma~\ref{lem:mis} to the above inequality.
    \item Inequality~\eqref{eq:inter2}. By \citet[Proof of Lemma 5]{han2022exact}, we have 
    \begin{equation}
           \onormSize{}{\mE (\mV - \mV^{(t)})}_{\sigma} \lesssim \sqrt{ \frac{r^{K-1}(p r^{K-1} + p r)}{p^{K-1}}} r \ell^{(t)}.
    \end{equation}
    Also, we complete the proof of inequality~\eqref{eq:inter1} by applying Lemma~\ref{lem:mis} to the above inequality.
    
    \item Inequality~\eqref{eq:inter3}. We upper bound the desired quantity by triangle inequality,
    \begin{equation}
        \onormSize{}{ [ \mW_{:b}^T \mY \mV ]^s   -  [ \mW_{:b}^{(t), T} \mY \mV  ]^s } \leq I_1 + I_2 + I_3,
    \end{equation}
    where 
    \begin{align}
        I_1 &=  \onorm{ \frac{ \mW_{:b}^T \mY \mV  }{\onormSize{}{  \mW_{:b}^T \mX \mV }} - \frac{ \mW_{:b}^{(t), T} \mY \mV  }{ \onormSize{}{\mW_{:b}^{(t), T} \mX \mV }}   },\\
        I_2 & = \onorm{\of{ \frac{1}{\onormSize{}{\mW_{:b}^T \mY \mV}} -  \frac{1}{\onormSize{}{\mW_{:b}^T \mX \mV}}  } \mW_{:b}^T \mY \mV },\\
        I_3 &= \onorm{\of{ \frac{1}{\onormSize{}{\mW_{:b}^{(t), T} \mY \mV  }} -  \frac{1}{\onormSize{}{\mW_{:b}^{(t), T} \mX \mV  }}  } \mW_{:b}^{(t), T} \mY \mV }.
    \end{align}
    Next, we upper bound the quantities $I_1, I_2, I_3$ separately. 
    
    For $I_1$, we further bound $I_1$ by triangle inequality,
    \begin{equation}
        I_1 \leq I_{11} + I_{12},
    \end{equation}
    where
 \begin{equation}
        I_{11} = \onorm{ \frac{ \mW_{:b}^T \mX \mV  }{\onormSize{}{  \mW_{:b}^T \mX \mV }} - \frac{ \mW_{:b}^{(t), T} \mX \mV  }{ \onormSize{}{\mW_{:b}^{(t), T} \mX \mV }}   },\quad 
        I_{12} = \onorm{ \frac{\mW_{:b}^T \mE \mV}{ \onormSize{}{  \mW_{:b}^T \mX \mV }} - \frac{\mW_{:b}^{(t),T} \mE \mV}{ \onormSize{}{  \mW_{:b}^{(t),T} \mX \mV }} }.
    \end{equation}
    We first consider $I_{11}$. Define the confusion matrix $\mD = \mM^T \mTheta^T \mW^{(t)} = \entry{D_{ab}} \in \bbR^{r \times r}$ where 
    \begin{equation}
        D_{ab} = \frac{\sum_{i \in [p]}  \theta(i) \ind \offf{ z(i) = a, z^{(t)}(i) = b } }{\sum_{i \in [p]}  \ind\offf{  z^{(t)}(i) = b }}, \text{  for all } a, b \in [r].
    \end{equation}
    By Lemma~\ref{lem:membership},  we have $\sum_{i \in [p]}  \ind\offf{  z^{(t)}(i) = b } \gtrsim p/r$. Then, we have
    \begin{equation}\label{eq:dab_bound}
        \sum_{a \neq b, a,b \in [r]} D_{ab} \lesssim \frac{r}{p} \sum_{i \colon z^{(t)}(i) \neq z(i)} \theta(i) \lesssim\frac{L^{(t)}}{\Delta_{\min}^2}  \lesssim \frac{1}{\log p}, 
    \end{equation}
    \normalsize
    and for all $b \in [r]$,
    \begin{align}
        D_{bb} &= \frac{\sum_{i \in [p]}  \theta(i) \ind \offf{ z(i) = z^{(t)}(i) = b } }{\sum_{i \in [p]}  \ind\offf{  z^{(t)}(i) = b }} \geq \frac{c(\sum_{i \in [p]}  \ind\offf{  z^{(t)}(i) = b } - p\ell^{(t)})}{\sum_{i \in [p]}  \ind\offf{  z^{(t)}(i) = b }} \gtrsim 1 - \frac{1}{\log p},\label{eq:dbb_bound}
    \end{align}
    under the inequality \eqref{eq:cond_intial} in Condition~\ref{cond:origin}. By the definition of $\mW,\mW^{(t)},\mV$, we have 
    \begin{align}
        \frac{ \mW_{:b}^T \mX \mV  }{\onormSize{}{  \mW_{:b}^T \mX \mV }} = \off{\mS_{b:}}^s, \quad
        \frac{ \mW_{:b}^{(t), T} \mX \mV  }{ \onormSize{}{\mW_{:b}^{(t), T} \mX \mV }} = [D_{bb} \mS_{b:} + \sum_{a \neq b, a \in [r]} D_{ab} \mS_{a:}]^s.
    \end{align}

    Let $\alpha$ denote the angle between $\mS_{b:}$ and $D_{bb} \mS_{b:} + \sum_{a \neq b, a \in [r]} D_{ab} \mS_{a:}$. To roughly estimate the range of $\alpha$, we consider the inner product 
    \begin{align}
        \ang{\mS_{b:},  D_{bb} \mS_{b:} + \sum_{a \neq b, a \in [r]} D_{ab} \mS_{a:}}  &= D_{bb} \onorm{\mS_{b:}}^2 + \sum_{a \neq b} D_{ab}\ang{\mS_{b:}, \mS_{a:}}\\
        &\geq D_{bb} \onorm{\mS_{b:}}^2 -  \sum_{a \neq b, a \in [r]} D_{ab}   \onorm{\mS_{b:}} \max_{a \in [r]}  \onorm{\mS_{a:}} \\
        & \geq C,
    \end{align}
    where $C$ is a positive constant, and the last inequality holds when $p$ is large enough following the constraint of $\onorm{\mS_{b:}}$ in parameter space~\eqref{eq:family} and the bounds of $\mD$ in \eqref{eq:dab_bound} and \eqref{eq:dbb_bound}.
    
    The positive inner product between $\mS_{b:}$ and $D_{bb} \mS_{b:} + \sum_{a \neq b, a \in [r]} D_{ab} \mS_{a:}$ indicates $\alpha \in [0,  \pi/2)$, and thus $2\sin \frac{\alpha}{2} \leq \sqrt{2} \sin \alpha$. Then, by the geometry property of trigonometric function, we have 
    \begin{align}
       \onormSize{}{ [D_{bb} \mS_{b:} + \sum_{a \neq b, a \in [r]} D_{ab} \mS_{a:}] \sin \alpha} 
       &= \onormSize{}{(\mI_d - \text{Proj}(\mS_{b:})) \sum_{a \neq b, a \in [r]} D_{ab} \mS_{a:}}\\
       &\leq \sum_{a \neq b, a \in [r]} D_{ab} \onorm{(\mI_d - \text{Proj}(\mS_{b:}))\mS_{a:}  }\\
       &= \sum_{a \neq b, a \in [r]} D_{ab} \onorm{ \mS_{a:} \sin (\mS_{b:}, \mS_{a:}) }\\
       & \leq \sum_{a \neq b, a \in [r]} D_{ab} \onorm{\mS_{a:}} \onorm{ \mS^s_{b:} - \mS^s_{a:}}, \label{eq:i11_sin_num}
    \end{align}
    where the first inequality follows from the triangle inequality, and the last inequality follows from Lemma~\ref{lem:norm_diff}. Note that with bounds \eqref{eq:dab_bound} and \eqref{eq:dbb_bound}, when $p$ is large enough, we have
    \begin{align}
          \onormSize{}{\mW_{:b}^{(t), T} \mX \mV } &= \onormSize{}{ D_{bb} \mS_{b:} + \sum_{a \neq b, a \in [r]} D_{ab} \mS_{a:}} \geq D_{bb} \onorm{\mS_{b:}} - \sum_{a \neq b, a \in [r]} D_{ab} \onorm{ \mS_{a:}} \geq  C_1,\label{eq:i11_sin_dom}
    \end{align}
    for some positive constant $C_1$. Notice that $I_{11} = \sqrt{1 - \cos \alpha} = 2 \sin \frac{\alpha}{2}$. Therefore, we obtain
    \begin{align}
        I_{11} &\leq \sqrt{2} \sin \alpha \\
        & = \frac{  \onormSize{}{ [D_{bb} \mS_{b:} + \sum_{a \neq b, a \in [r]} D_{ab} \mS_{a:}] \sin \alpha}}{  \onormSize{}{ D_{bb} \mS_{b:} + \sum_{a \neq b, a \in [r]} D_{ab} \mS_{a:}} } \\
        & \leq \frac{ 1}{ C_1}  \sum_{a \neq b, a \in [r]} D_{ab} \onorm{\mS_{a:}} \onorm{ \mS^s_{b:} - \mS^s_{a:}} \\
        & \lesssim \frac{r}{p} \sum_{i \in [p]} \theta(i) \sum_{b \in [r]} \ind\offf{z^{(t)}(i) = b} \onorm{ \mS^s_{b:} - \mS^s_{a:}}\\
        & \leq \frac{r L^{(t)}}{\Delta_{\min}},\label{eq:i11}
    \end{align}
    where the second inequality follows from~\eqref{eq:i11_sin_num} and \eqref{eq:i11_sin_dom}, and the last two inequalities follow by the definition of $D_a$ and $L^{(t)}$, and the constraint of $\onorm{\mS_{b:}}$ in parameter space~\eqref{eq:family}.

    We now consider $I_{12}$. By triangle inequality, we have 
    \begin{align}
        I_{12} &\leq \frac{1}{\onormSize{}{  \mW_{:b}^T \mX \mV }}\onormSize{}{ (\mW_{:b}^T - \mW_{:b}^{(t),T}) \mE \mV }  + \frac{ \onormSize{}{(\mW_{:b}^T - \mW_{:b}^{(t),T}) \mX \mV} }{\onormSize{}{  \mW_{:b}^T \mX \mV } \onormSize{}{  \mW_{:b}^{(t),T} \mX \mV } }\onormSize{}{ \mW_{:b}^{(t),T} \mE \mV}.
    \end{align}
    By \citet[Proof of Lemma 5]{han2022exact}, we have 
    \begin{equation}\label{eq:i12_wwe}
        \onormSize{}{ (\mW_{:b}^T - \mW_{:b}^{(t),T}) \mE \mV }  \lesssim \sqrt{\frac{ r^{2K} + p r^{K+1}}{p^K}} \frac{\sqrt{L^{(t)}}}{\Delta_{\min}}.
    \end{equation}
    Notice that 
    \begin{align}
        \onormSize{}{(\mW_{:b}^T - \mW_{:b}^{(t),T}) \mX \mV} &\leq \onormSize{}{\mW_{:b}^T - \mW_{:b}^{(t),T}} \onorm{\mX\mV}_F \lesssim \frac{r^{3/2}L^{(t)}}{\sqrt{p}\Delta_{\min}^2} \onormSize{}{\mS} \onormSize{}{\mTheta \mM}_{\sigma} \lesssim \frac{ \sqrt{r L^{(t)}} }{\Delta_{\min}}, \label{eq:i12_wwx}
    \end{align}
    where the second inequality follows from~\citet[Inequality (121), Proof of Lemma 5]{han2022exact} and the last inequality follows from Lemma~\ref{lem:singular_thetam} and \eqref{eq:cond_intial} in Condition~\ref{cond:origin}.
     Note that $\onorm{\mW_{:b}^T \mX \mV} = \onorm{\mS_{b:}} \geq c_3$ and $\onormSize{}{  \mW_{:b}^{(t),T} \mX \mV } \geq C_1$ by inequality~\eqref{eq:i11_sin_dom}. Therefore, we have 
     \begin{align}
         I_{12} &\lesssim  \onormSize{}{ (\mW_{:b}^T - \mW_{:b}^{(t),T}) \mE \mV } +   \onormSize{}{(\mW_{:b}^T - \mW_{:b}^{(t),T}) \mX \mV}   \onormSize{}{ \mW_{:b}^{(t),T} \mE \mV} \\
         &\lesssim \sqrt{\frac{ r^{2K} + p r^{K+1}}{p^K}} \frac{\sqrt{L^{(t)}}}{\Delta_{\min}} + \frac{ \sqrt{r L^{(t)}} }{\Delta_{\min}}\sqrt{\frac{r^{2K}}{p^K}} \\
         & \lesssim \sqrt{\frac{ r^{2K} + p r^{K+1}}{p^K}} \frac{\sqrt{L^{(t)}}}{\Delta_{\min}},\label{eq:i12}
     \end{align}
     where second inequality follows from the inequalities~\eqref{eq:i12_wwe}, \eqref{eq:i12_wwx}, and \eqref{eq:cond1} in Condition~\ref{cond:origin}.

    Hence, combining inequalities~\eqref{eq:i11} and~\eqref{eq:i12} yields  
    \begin{equation}\label{eq:i1}
        I_1 \lesssim \frac{r L^{(t)}}{\Delta_{\min}} + \sqrt{\frac{ r^{2K} + p r^{K+1}}{p^K}} \frac{\sqrt{L^{(t)}}}{\Delta_{\min}}.
    \end{equation}
    
        For $I_2$ and $I_3$, recall that $\onorm{\mW_{:b}^T \mX \mV} = \onorm{\mS_{b:}} \geq c_3$ and $\onormSize{}{  \mW_{:b}^{(t),T} \mX \mV } \geq C_1$ by inequality~\eqref{eq:i11_sin_dom}. By triangle inequality and \eqref{eq:cond1} in Condition~\ref{cond:origin}, we have 
    \begin{equation}\label{eq:i2}
         I_2 \leq \frac{ \onormSize{}{ \mW_{:b}^T \mE \mV } }{ \onormSize{}{\mW_{:b}^T \mX \mV}  } \lesssim  \onormSize{}{ \mW_{:b}^T \mE \mV }  \lesssim \frac{r^K}{p^{K/2}},
    \end{equation}
    and 
    \begin{equation}\label{eq:i3}
        I_3 \leq \frac{ \onormSize{}{ \mW_{:b}^{(t),T} \mE \mV } }{ \onormSize{}{\mW_{:b}^{(t),T} \mX \mV}  }  \lesssim \onormSize{}{ \mW_{:b}^{(t),T} \mE \mV } \lesssim \frac{r^K}{p^{K/2}}.
     \end{equation}
       Therefore, combining the inequalities~\eqref{eq:i1}, \eqref{eq:i2}, and \eqref{eq:i3}, we finish the proof of inequality~\eqref{eq:inter3}.
    
    \item Inequality~\eqref{eq:inter4}. Here we only show the proof of inequality~\eqref{eq:inter4} with $\mW_{:b}^{(t)}$. The proof also holds by replacing $\mW_{:b}^{(t)}$ to $\mW_{:b}$, and we omit the repeated procedures.
    
    We upper bound the desired quantity by triangle inequality
    \begin{equation}
         \onormSize{}{  [   \mW_{:b}^{(t), T}  \mY \mV]^s -  [  \mW_{:b}^{(t), T}  \mY \mV^{(t)}]^s }  \leq J_1 + J_2 + J_3,
    \end{equation}
    where 
    \begin{align}
        J_1 &= \onorm{ \frac{\mW_{:b}^{(t),T} \mY \mV }{\onormSize{}{ \mW_{:b}^{(t),T} \mX \mV }} - \frac{\mW_{:b}^{(t),T} \mY \mV^{(t)} }{\onormSize{}{ \mW_{:b}^{(t),T} \mX \mV^{(t)} }}  },\\
        J_2 & = \onorm{\of{ \frac{1}{\onormSize{}{\mW_{:b}^{(t),T} \mY \mV}} -  \frac{1}{\onormSize{}{\mW_{:b}^{(t),T} \mX \mV}}  } \mW_{:b}^{(t),T} \mY \mV },\\
        J_3 &= \onorm{\of{ \frac{1}{\onormSize{}{\mW_{:b}^{(t),T} \mY \mV^{(t)}}} -  \frac{1}{\onormSize{}{\mW_{:b}^{(t),T} \mX \mV^{(t)}} }  } \mW_{:b}^{(t), T} \mY \mV^{(t)} }.
    \end{align}
    \normalsize
    Next, we upper bound the quantities $J_1, J_2, J_3$ separately. 

    For $J_1$, by triangle inequality,
we have 
    \begin{equation}
        J_1 \leq J_{11} + J_{12},
    \end{equation}
    where 
    \begin{equation}
        J_{11} =  \onorm{ \frac{\mW_{:b}^{(t),T} \mX \mV }{\onormSize{}{ \mW_{:b}^{(t),T} \mX \mV }} - \frac{\mW_{:b}^{(t),T} \mX \mV^{(t)} }{\onormSize{}{ \mW_{:b}^{(t),T} \mX \mV^{(t)} }}  }, \quad 
         J_{12} =  \onorm{ \frac{\mW_{:b}^{(t),T}  \mE \mV}{ \onormSize{}{\mW_{:b}^{(t),T}  \mX \mV} } - \frac{\mW_{:b}^{(t),T}  \mE \mV^{(t)}}{ \onormSize{}{\mW_{:b}^{(t),T}  \mX \mV^{(t)}} }   }.
    \end{equation}
    
    We first consider $J_{11}$.  Define the matrix $\mV^{k} \coloneqq \mW^{\otimes (k-1)} \otimes \mW^{(t), \otimes (K-k)}$ for $k = 2,\ldots, K-1$, and denote $\mV^1 = \mV^{(t)}, \mV^K = \mV$. Also, define the quantity
    \begin{equation}
        J_{11}^k = \onormSize{}{ [\mW_{:b}^{(t),T} \mX \mV^k ]^s - [\mW_{:b}^{(t),T} \mX \mV^{k+1} ]^s },
     \end{equation}
    for $k = 1,\ldots, K-1$. Let $\beta_k$ denote the angle between $\mW_{:b}^{(t),T} \mX \mV^k$ and $\mW_{:b}^{(t),T} \mX \mV^{k+1}$. With the same idea to prove $I_{11}$ in inequality~\eqref{eq:i11}, we bound $J_{11}^k$ by the trigonometric function of $\beta_k$. 
    
    To roughly estimate the range of $\beta_k$, we consider the inner product between $\mW_{:b}^{(t),T} \mX \mV^k$ and $\mW_{:b}^{(t),T} \mX \mV^{k+1}$. Before the specific derivation of the inner product, note that 
    \begin{equation}
        \mW_{:b}^{(t),T} \mX \mV^k = \mat_1(\tT_k), \quad  \mW_{:b}^{(t),T} \mX \mV^{k+1} = \mat_1(\tT_{k+1}),
    \end{equation}
    where 
    \begin{align}
        \tT_k &= \tX \times_1  \mW_{:b}^{(t),T} \times_2 \mW^T \times_3 \cdots \times_k \mW^T  \times_{k+1} \mW^{(t),T} \times_{k+2} \cdots \times_K \mW^{(t),T}\\
        \tT_{k+1} &= \tX \times_1  \mW_{:b}^{(t),T} \times_2 \mW^T \times_3 \cdots \times_k \mW^T \times_{k+1} \mW^{T} \times_{k+2} \cdots \times_K \mW^{(t),T}. 
    \end{align}
    Recall the definition of confusion matrix $\mD = \mM^T \mTheta^T \mW^{(t)} = \entry{D_{ab}} \in \bbR^{r \times r}$.
 We have 
    \begin{align}
         \ang{\mW_{:b}^{(t),T} \mX \mV^k , \mW_{:b}^{(t),T} \mX \mV^{k+1}} 
         &  =\ang{ \mat_{k+1}(\tT_k), \mat_{k+1}(\tT_{k+1}) } \label{eq:j11_tensormat} \\
  & =  \ang{ \mD^T \mS \mZ^k ,  \mS \mZ^k   }\\
         & = \sum_{b \in [r]} \of{ D_{bb} \onormSize{}{\mS_{b:}\mZ^k}^2 + \sum_{a \neq b, a \in [r]} D_{ab} \ang{\mS_{a:} \mZ^k, \mS_{b:} \mZ^k}} \\
         & \gtrsim (1 - \log p^{-1}) \min_{a \in [r]}\onormSize{}{\mS_{a:}\mZ^k}^2 - \log p^{-1} \max_{a \in [r]} \onormSize{}{\mS_{a:}\mZ^k}^2, \label{eq:j11inner}
    \end{align}
    \normalsize
    where $\mZ^k = \mD_{:b} \otimes \mI_r^{\otimes (k-1)} \otimes  \mD^{\otimes (K-k-1)}$, the equations follow by the tensor algebra and definitions, and the last inequality follows from the bounds of $\mD$ in \eqref{eq:dab_bound} and \eqref{eq:dbb_bound}. 
    
    Note that 
    \begin{align}
        \onorm{\mD}_{\sigma} \leq \onorm{\mD}_{F} \leq  \sqrt{\sum_{b \in [r]} D_{bb}^2 + (\sum_{a \neq b, a ,b \in [r]} D_{ab} )^2} \lesssim \sqrt{r  + \log^2 p^{-1}} \lesssim 1,\label{eq:singular_Du}
    \end{align}
    where the second inequality follows from inequality~\eqref{eq:dab_bound}, and the fact that for all $b \in [r]$,
    \begin{equation}
        D_{bb} \lesssim \frac{r}{p} \sum_{i\colon z(i) = b} \theta(i) \lesssim 1.
    \end{equation}
    Also, we have 
    \begin{equation}\label{eq:singular_Dl}
        \lambda_r(\mD) \geq \lambda_r(\mW^{(t)}) \lambda_r(\mTheta \mM) \gtrsim 1,
    \end{equation}
    following the Lemma~\ref{lem:singular_thetam} and Lemma~\ref{lem:membership}. Then, for all $k \in [K]$, we have 
    \begin{align}
        1 &\lesssim \onorm{\mD_{:b}} \lambda_r(\mD)^{K-k-1}  \leq \lambda_{r^{K-2}}(\mZ^{k}) \leq \onormSize{}{\mZ^k}_{\sigma} \leq \onorm{\mD_{:b}} \onorm{\mD}_{\sigma}^{K-k-1} \lesssim 1.\label{eq:singular_z}
    \end{align}
    Thus, we have bounds 
    \begin{equation}
         \max_{a \in [r]} \onormSize{}{\mS_{a:} \mZ^k} \leq \max_{a \in [r]} \onorm{\mS_{a:}} \onormSize{}{\mZ^k}_{\sigma} \lesssim 1, \quad 
        \min_{a \in [r]} \onormSize{}{\mS_{a:}\mZ^k}  \geq \min_{a \in [r]} \onorm{\mS_{a:}} \lambda_{r^{K-2}}(\mZ^{k}) \gtrsim 1.
    \end{equation}
    Hence, when $p$ is large enough, the inner product~\eqref{eq:j11inner} is positive, which implies $\beta_k \in [0, \pi/2)$ and thus $2 \sin \frac{\beta_k}{2} \leq \sqrt{2} \sin \beta_k$.
    
    Next, we upper bound the trigonometric function $\sin \beta_k$. Note that
    \begin{align}
         \sin \beta_k &= \sin ({\mD_{:b}^{T} \mS \mI_r^{\otimes k-1} \otimes \mD^{\otimes K-k} }, {\mD_{:b}^{T} \mS \mI_r^{\otimes k} \otimes \mD^{\otimes K-k-1} } )\\
         &\leq \sin \beta_{k1} + \sin \beta_{k2},
    \end{align}
    where
    \begin{align}
         \sin \beta_{k1} &= \sin ({\mD_{:b}^{T} \mS \mI_r^{\otimes k-1} \otimes \mD^{\otimes K-k} }, {\mD_{:b}^{T} \mS \mI_r^{\otimes k-1} \otimes \tilde \mD \otimes \mD^{\otimes K-k-1} } ),\\
         \sin \beta_{k2} & = \sin ({\mD_{:b}^{T} \mS \mI_r^{\otimes k-1} \otimes \tilde \mD \otimes \mD^{\otimes K-k-1} } ,  {\mD_{:b}^{T} \mS \mI_r^{\otimes k} \otimes \mD^{\otimes K-k-1} }   ),
    \end{align}
    and $\tilde \mD$ is the normalized confusion matrix with entries $\tilde \mD_{ab} = \frac{\sum_{i \in [p]} \theta(i) \ind\{z^{(t)} = b, z(i) = a\}}{\sum_{i \in [p] }\theta(i) \ind\{z^{(t)} = b\}}$.

    To bound $\sin \beta_{k1}$, recall Definition~\ref{def:stable} that for any cluster assignment $\bar z$ in the $\varepsilon$-neighborhood of true $z$,
    \begin{align}
        &\mp(\bar z)=(|\bar z^{-1}(1)|, \ldots,|\bar z^{-1}(r)|)^T,\quad \mp_{\mtheta}(\bar z)=(\onormSize{}{\mtheta_{\bar z^{-1}(1)}}_1,\ldots,\onormSize{}{\mtheta_{\bar z^{-1}(r)}}_1)^T.
    \end{align}
    
    {Note that we have $\ell^{(t)} \leq \frac{L^{(t)}}{\Delta_{\min}^2} \leq \frac{\bar C}{\tilde C} r \log^{-1}(p) $ by Condition~\ref{cond:origin} and Lemma~\ref{lem:mis}. Then,  with the locally linear stability assumption, the $\mtheta$ is $\ell^{(t)}$-locally linearly stable; i.e.,}
    \begin{equation}
         \sin (\mp(z^{(t)}), \mp_{\mtheta}(z^{(t)})) \lesssim \frac{L^{(t)}}{\Delta_{\min}}. 
    \end{equation}
    Note that $\text{diag}(\mp(z^{(t)})) \mD = \text{diag}(\mp_{\mtheta}(z^{(t)})) \tilde \mD$, and $\sin(\ma,\mb) =\min_{c \in \bbR} \frac{ \onormSize{}{\ma - c\mb} }{\onormSize{}{\ma}}$ for vectors $\ma,\mb$ of same dimension. Let $c_0 = \argmin_{c \in \bbR} \frac{\onormSize{}{\mp(z^{(t)}) - c\mp_{\mtheta}(z^{(t)})}}{\onormSize{}{\mp(z^{(t)})}}$. Then, we have
    \begin{align}
      \min_{c \in \bbR} \onormSize{}{ \mD - c \tilde \mD}_F
      &\leq  \onormSize{}{ \mI_r - c_0 \text{diag}(\mp(z^{(t)})) \text{diag}^{-1}(\mp_{\mtheta}(z^{(t)}))}_F \onormSize{}{\mD}_F\\
      &\lesssim \frac{ \onormSize{}{ \mp(z^{(t)}) - c_0\mp_{\mtheta}(z^{(t)})} }{\min_{a \in [r]} \onormSize{}{\mtheta_{z^{(t),-1}(a)} }_1  } \\
      & = \frac{\onormSize{}{ \mp(z^{(t)})} }{\min_{a \in [r]} \onormSize{}{\mtheta_{z^{(t),-1}(a)} }_1 } \sin (\mp(z^{(t)}), \mp_{\mtheta}(z^{(t)}))\\
      & \lesssim \frac{L^{(t)}}{\Delta_{\min}},
    \end{align}
    where the last inequality follows from Lemma~\ref{lem:membership}, the constraint $\min_{i \in [p]}\theta(i) \geq c>0$, $\onormSize{}{ \mp(z^{(t)})} \lesssim p$ and $\min_{a \in [r]} \onormSize{}{\mtheta_{z^{(t),-1}(a)} }_1  \gtrsim p$. 
    
    By the geometry property of trigonometric function, we have
    \begin{align}
        \sin \beta_{k1} &= \min_{c \in \bbR} \frac{\onormSize{}{\mD_{:b}^{T} \mS \mI_r^{\otimes k-1} \otimes (\mD - c \tilde \mD) \otimes \mD^{\otimes K-k - 1} } }{ \onormSize{}{\mD_{:b}^{T} \mS \mI_r^{\otimes k-1} \otimes \mD^{\otimes K-k}} }\\
        &\leq \frac{ \onormSize{}{\mD_{:b}^{T} \mS} \onormSize{}{\mD - c_0 \tilde \mD}_{\sigma} \onormSize{}{\mD}_{\sigma}^{K-k-1} }{\onormSize{}{\mD_{:b}^{T} \mS} \lambda_r^{K-k}(\mD) }\\
        & \lesssim  \onormSize{}{\mD - c_0 \tilde \mD}_{F}\\
        & \lesssim \frac{L^{(t)}}{\Delta_{\min}}, \label{eq:b1}
    \end{align}
    where the second inequality follows from the singular property of $\mD$ in \eqref{eq:singular_Du}, \eqref{eq:singular_Dl} and the constraint of $\mS$ in \eqref{eq:family}.
    
    To bound $\sin \beta_{k2}$, let $\mC = \text{diag}(\{\onormSize{}{\mS_{a:}}\}_{a \in [r]})$. We have 
    \begin{align}
        \sin \beta_{k2} &\lesssim \frac{ \onorm{\mD_{:b}^{T} \mS \mI_r^{\otimes k-1} \otimes (\mI_r - \tilde \mD) \otimes \mD^{\otimes K-k-1} } }{\onormSize{}{\mD_{:b}^{T} \mS \mI_r^{\otimes k} \otimes \mD^{\otimes K-k-1} }  }\\
        & \lesssim \frac{ \onormSize{}{(\mI_r - \tilde \mD^T) \mS \mZ^k}_F }{\onormSize{}{\mD_{:b}^{T} \mS} \lambda_r^{K-k-1}(\mD)} \\
        & \lesssim \onormSize{ }{(\mI_r - \tilde \mD^T) \mS \mC^{-1}}_F \onormSize{}{\mC \mZ^k}_{\sigma}\\
        & \lesssim \frac{r}{p} \sum_{i \in [p]} \theta(i) \sum_{b \in [r]} \ind\{z^{(t)}(i) = b\} \onormSize{}{\mS_{b:}^s - \mS_{z(i):}^s} \\
        & \lesssim \frac{L^{(t)}}{\Delta_{\min}}, \label{eq:b2} 
    \end{align}
    where the third inequality follows from the singular property of $\mD$ and the boundedness of $\mS$, and the fourth inequality follows from the definition of $\tilde \mD$, boundedness of $\mS$, the lower bound of $\mtheta$, and the singular property of $\mZ^k$ in inequality~\eqref{eq:singular_z}, and the last line follows from the definition of $L^{(t)}$.
    
  Combining~\eqref{eq:b1} and~\eqref{eq:b2} yields
    \begin{equation}
        \sin \beta_k \leq \sin \beta_{k1} + \sin \beta_{k2} \lesssim \frac{L^{(t)}}{\Delta_{\min}}.
    \end{equation}

       Finally, by triangle inequality, we obtain
    \begin{equation}\label{eq:j1}
        J_{11} \leq \sum_{k = 1}^{K-1} J_{11}^k \lesssim  \sum_{k = 1}^{K-1} \sin \beta_k  \lesssim (K-1)\frac{r L^{(t)}}{\Delta_{\min}}.
    \end{equation}

    We now consider $J_{12}$. By triangle inequality, we have 
    \begin{align}
        J_{12} &\leq \frac{ 1}{\onormSize{}{\mW_{:b}^{(t),T} \mX \mV }} \onormSize{}{\mW_{:b}^{(t),T} \mE (\mV-\mV^{(t)})}   +
        \frac{ \onormSize{}{\mW_{:b}^{(t),T} \mX (\mV-\mV^{(t)})} }{\onormSize{}{\mW_{:b}^{(t),T} \mX \mV }\onormSize{}{\mW_{:b}^{(t),T} \mX \mV^{(t)} }} \onormSize{}{\mW_{:b}^{(t),T} \mE \mV^{(t)} }.
    \end{align}
     Note that 
     \begin{align}
          \onormSize{}{\mW_{:b}^{(t),T} \mX \mV^{(t)}} &= \onormSize{}{\mD^T \mS \mZ^1} 
      \geq \lambda_r(\mD) \onorm{\mS} \lambda_{r^{K-2}}(\mZ^1) \gtrsim 1, \label{eq:j11_wxvt}
     \end{align}
    where the inequality follows from the bounds~\eqref{eq:singular_Dl} and \eqref{eq:singular_z}.
    
    By \citet[Proof of Lemma 5]{han2022exact}, we have 
    \begin{align}
        \onormSize{}{\mW_{:b}^{(t),T} \mE (\mV-\mV^{(t)})} \lesssim \sqrt{ \frac{r^{2K+1} + p r^{2 + K}}{p^K} } \frac{(K-1)\sqrt{L^{(t)}}}{\Delta_{\min}}. \label{eq:j12_wev}
    \end{align}
    Notice that 
    \begin{align}
         \onormSize{}{\mX (\mV^k-\mV^{k+1})}_F &\leq \onormSize{}{(\mI - \mD^T) \mS (\mI_r^{\otimes (k-1)} \otimes  \mD^{\otimes (K-k-1)})}_F\notag \\
        &\leq \onormSize{}{(\mW^T - \mW^{(t),T}) \mTheta \mM}_F \onorm{\mS}_F \onorm{ \mD}_{\sigma}^{K-k-1}\notag \\
        & \lesssim \onormSize{}{ \mW^T - \mW^{(t),T} } \onorm{\mTheta \mM}_{\sigma} \notag\\
        & \lesssim \frac{ \sqrt{r L^{(t)}}}{\Delta_{\min}}, \label{eq:step}
    \end{align}
    where the first inequality follows from the tensor algebra in 
    inequality~\eqref{eq:j11inner}, the second inequality follows from the fact  that $\mI = \mW^T \mTheta \mM$, and the last inequality follows from \citet[Proof of Lemma 5]{han2022exact}. It follows from~\eqref{eq:step} and Lemma~\ref{lem:membership} that
    \begin{align}
        \onormSize{}{\mW_{:b}^{(t),T} \mX (\mV-\mV^{(t)})}&\leq \onormSize{}{\mW_{:b}^{(t),T} }  \sum_{k = 1}^{K-1} \onormSize{}{\mX (\mV^k-\mV^{k+1})}_F\lesssim \frac{ \sqrt{r L^{(t)}}}{\sqrt{p} \Delta_{\min}}.\label{eq:j12_wxv}
    \end{align}
    \normalsize
    Note that $\onormSize{}{\mW_{:b}^{(t),T} \mX \mV }$ and $\onormSize{}{\mW_{:b}^{(t),T} \mX \mV^{(t)} }$ are lower bounded by inequalities~\eqref{eq:i11_sin_dom} and \eqref{eq:j11_wxvt}, respectively. We have 
    \begin{align}
        J_{12} &\lesssim \onormSize{}{\mW_{:b}^{(t),T} \mE (\mV-\mV^{(t)})}  + \onormSize{}{\mW_{:b}^{(t),T} \mX (\mV-\mV^{(t)})}\onormSize{}{\mW_{:b}^{(t),T} \mE \mV^{(t)} }\\
        &\lesssim \sqrt{ \frac{r^{2K+1} + p r^{2 + K}}{p^K} } \frac{\sqrt{L^{(t)}}}{\Delta_{\min}} + \frac{ \sqrt{r L^{(t)}}}{\sqrt{p} \Delta_{\min}} \sqrt{\frac{r^{2K}}{p^K}} \\
        & \lesssim  \sqrt{ \frac{r^{2K+1} + p r^{2 + K}}{p^K} } \frac{\sqrt{L^{(t)}}}{\Delta_{\min}},
    \end{align}
    where the second inequality follows from inequalities~\eqref{eq:j12_wev}, \eqref{eq:j12_wxv}, and the inequality \eqref{eq:cond1} in Condition~\ref{cond:origin}.

    For $J_2$ and $J_3$, recall that $\onormSize{}{\mW_{:b}^{(t),T} \mX \mV }$ and $\onormSize{}{\mW_{:b}^{(t),T} \mX \mV^{(t)} }$ are lower bounded by inequalities~\eqref{eq:i11_sin_dom} and \eqref{eq:j11_wxvt}, respectively. By triangle inequality and inequality \eqref{eq:cond1} in Condition~\ref{cond:origin}, we have 
    \begin{equation}\label{eq:j2}
        J_2 \leq \frac{ \onormSize{}{\mW_{:b}^{(t),T} \mE \mV}  }{  \onormSize{}{\mW_{:b}^{(t),T} \mX \mV} } \lesssim \onormSize{}{\mW_{:b}^{(t),T} \mE \mV} \lesssim \frac{r^K}{p^{K/2}},
    \end{equation}
    and 
    \begin{equation}\label{eq:j3}
        J_3 \leq \frac{ \onormSize{}{\mW_{:b}^{(t),T} \mE \mV^{(t)}}  }{  \onormSize{}{\mW_{:b}^{(t),T} \mX \mV^{(t)}} } \lesssim \onormSize{}{\mW_{:b}^{(t),T} \mE \mV} \lesssim \frac{r^K}{p^{K/2}}.
    \end{equation}
    
    Therefore, combining the inequalities~\eqref{eq:j1}, \eqref{eq:j2}, and \eqref{eq:j3}, we finish the proof of inequality~\eqref{eq:inter4}.

    \item Inequality~\eqref{eq:inter5}. By triangle inequality, we upper bound the desired quantity 
     \begin{align}
        &\onormSize{}{ [ \mW_{:b}^T \mY \mV^{(t)} ]^s   -  [ \mW_{:b}^{(t), T} \mY \mV^{(t)} ]^s } \\
        &\leq \onormSize{}{[ \mW_{:b}^T \mY \mV^{(t)} ]^s - [ \mW_{:b}^T \mY \mV ]^s }  + \onormSize{}{[\mW_{:b}^T \mY \mV ]^s  - [ \mW_{:b}^{(t),T} \mY \mV ]^s } + \onormSize{}{[\mW_{:b}^{(t),T} \mY \mV ]^s  -  [\mW_{:b}^{(t), T} \mY \mV^{(t)} ]^s }\\
        &\lesssim \frac{r L^{(t)}}{\Delta_{\min}} + \sqrt{ \frac{rr^{2K} + p r^{K+2}}{p^{K}}}  \frac{\sqrt{L^{(t)}}}{\Delta_{\min}},
    \end{align}
    \normalsize
    following the inequalities~\eqref{eq:inter3} and \eqref{eq:inter4}. Therefore, we finish the proof of inequality~\eqref{eq:inter5}.
        \end{enumerate}
  {
    
  \textbf{Next, we show the intermediate inequalities holds with $\mP, \mQ$ and $L(\hat z)$.}
        
        Consider the MLE confusion matrix $\hat \mD = \mM^T \mTheta^T \hat \mP = \entry{\hat D_{ab}} \in \bbR^{r \times r}$ with entries
\begin{align}
    &\hat D_{ab}
    = \frac{ \sum_{i \in [p]} \theta(i) \hat \theta(i) \ind\{ z(i) = a, \hat z(i) = b  \} }{\onormSize{}{\hat \mtheta_{\hat z^{-1}(b)}}^2}\\
    &\ = \frac{ \sum_{i \in [p]} (1 + o(p^{K-2}))(\hat \theta(i))^2 \ind\{ z(i) = a, \hat z(i) = b  \} }{\onormSize{}{\hat \mtheta_{\hat z^{-1}(b)}}^2}, \label{eq:mle_confusion}
\end{align}
\normalsize
where the second equation follows from Lemma~\ref{lem:poly_mle_degree}, and thus $\sum_{a \in [r]}\hat \mD_{ab} = 1 + o(1).$ By the derivation of~\eqref{eq:dab_bound}, \eqref{eq:dbb_bound}, \eqref{eq:singular_Dl}, and \eqref{eq:singular_Du}, we have 
\begin{equation}
     \sum_{a \neq b \in [r]} \hat D_{ab} \lesssim \frac{1}{p} \sum_{i \in [p]} \ind\{ \hat z(i) \neq z(i) \} (\hat \theta(i))^2 \lesssim \frac{1}{\log p}, \quad
     \hat D_{bb} \gtrsim 1- \frac{1}{\log p}, \quad \lambda_{\min}(\hat \mD) \asymp \onormSize{}{\hat \mD}_{\sigma} = (1 + o(1)).
\end{equation}
for all $a \neq b \in [r]$.

Now, we are ready to show the intermediate inequalities. First, by Lemma~\ref{lem:angle_gap_x} and $\min_{ i \in [p]} \theta(i) \geq c$, we have 
\begin{equation}
    \onormSize{}{\mS_{a:}^s - \mS_{b:}^s} \asymp \onormSize{}{\mA_{a:}^s - \mA_{b:}^s}.
\end{equation}
Then we can replace the $L^{(t)}$ by $L(\hat z)$ in the proof of Lemma~\ref{lem:intermediate}. The analogies of inequalities~\eqref{eq:inter1}, \eqref{eq:inter2}, \eqref{eq:inter3},  \eqref{eq:inter4}, and \eqref{eq:inter5}  hold by using the MLE confusion matrix and the definition of $L(\hat z)$. 

Particularly, for the analogy of \eqref{eq:inter4}, the usage of MLE confusion matrix avoids the stability condition on $\mtheta$. Let $\bar \mD$ be the normalized version of $\hat \mD$. The angle in inequality~\eqref{eq:b1} decays to 0 at speed $p^{-(K-2)} \lesssim \Delta_{\min}$ when $K \geq 3$, and the inequality~\eqref{eq:b2} holds by the fact that 
\begin{align}
    \onormSize{}{(\mI_r - \bar \mD) \mS \mC^{-1}}_F 
    &\lesssim \frac{r}{p} \sum_{i \in [p]} (\theta(i))^2 \sum_{b \in [r]} \onormSize{}{\mS_{b:}^s - \mS_{z(i):}^s} \lesssim \frac{r}{p} \sum_{i \in [p]} (\theta(i))^2 \sum_{b \in [r]} \onormSize{}{\mA_{b:}^s - \mA_{z(i):}^s}.
\end{align} 
  }     
       
    \end{proof}

\begin{lem}[Polynomial estimation error of MLE]\label{lem:poly_mle_degree}  Let $(\hat z, \hat \tS, \hat \mtheta)$ denote the MLE in \eqref{eq:mle} {with fixed $K \geq 2$ and symmetric mean tensor}, and $\hat \tX$ denote the mean tensor consisting of parameter $(\hat z, \hat \tS, \hat \mtheta)$. {With high probability going to 1 as $p \rightarrow \infty$, we} have 
\begin{equation}
    \onormSize{}{\tX - \hat \tX}_F^2 \lesssim \sigma^2 \of{ r^{K} + K p r},
\end{equation}
with probability going to 1. 
When {$\text{SNR} \gtrsim p^{-(K-1)} \log p$,  $\mtheta$ is balanced, and $\min_{i \in [p]} \theta(i) \geq c$ for some positive constant $c$}, the MLE satisfies
\begin{align}
    &\frac{1}{p} \sum_{i \in [p]} \ind\{ \hat z(i) \neq z(i) \} (\theta(i))^2 \lesssim \frac{1}{r \log p}, \quad \frac{1}{p} \sum_{i \in [p]} \ind\{ \hat z(i) \neq z(i) \} (\hat \theta(i))^2 \lesssim \frac{1}{r \log p} ,\  \text{and } L(\hat z) \lesssim \frac{\Delta_{\min}^2}{r \log p},
\end{align}
Further, we have 
\begin{equation}
     \theta(i)^2 = (1 + o(p^{-(K-2)})) \hat \theta(i)^2.
 \end{equation}
\end{lem}

\begin{proof}[Proof of Lemma~\ref{lem:poly_mle_degree}] Without loss of generality, we assume $\sigma^2 = 1$ and identity mapping minimizes the misclustering error for MLE.
For arbitrary two sets of parameters $(z, \tS, \mtheta), (z', \tS', \mtheta') \in \tP(\gamma)$ and corresponding mean tensors $\tX, \tX'$, we have 
\begin{align}
    &\text{rank}(\mat_k(\tX) - \mat_k(\tX')) \leq \text{rank}(\mat_k(\tX)) + \text{rank}(\mat_k(\tX)) \leq 2 r, \quad k \in [K].
\end{align}

Hence, we have
\begin{equation}\label{eq:dtbm_tucker}
    \tX - \tX' \in \tQ(2r, \ldots, 2 r),
\end{equation}
where $\tQ(r, \ldots, r) \coloneqq \{ \text{Tucker tensor with rank }(r, \ldots,  r) \} $.

Then, we obtain that 
\begin{align}
      \bbP(\FnormSize{}{\tX - \hat \tX_{ML}} \geq t) 
     &\leq 2  \bbP\of{ \sup_{ \tX, \tX' \in \tP(r, 
\ldots, r)} \ang{ \frac{\tX - \tX'}{\FnormSize{}{\tX - \tX'}}, \tE  } \geq t } \\
     &\leq 2  \bbP\of{ \sup_{ \tT \in \tQ(2r, \ldots, 2r) \cap \{ \onormSize{}{\tT}_F = 1\} } \ang{ \tT, \tE  } \geq t } \\
     &\lesssim  \exp(-Kpr),
\end{align}
with the choice $t\asymp \sigma \sqrt{(Kpr + r^K)}$. Here the first inequality follows from \citet[Lemma 1]{wang2019multiway}, the second inequality follows from~\eqref{eq:dtbm_tucker}, and the last inequality follows from \citet[Lemma E5]{han2022optimal}.

When $\Delta_{\min}^2 \gtrsim p^{-(K-1)} \log p$, we replace the vector $\hat x_{\hat z(i)}$ and $\hat \mX$ by our MLE estimator in the proof of Theorem~\ref{thm:initial}. With estimation error $\onormSize{}{\tX - \hat \tX}_F^2 \lesssim \of{ r^{K} + K p r}$ and $\Delta_{\min}^2 \gtrsim p^{-(K-1)} \log p$, we have 
\begin{align}
      \frac{1}{p} \sum_{i \in [p]} \ind\{ \hat z(i) \neq z(i) \} (\theta(i))^2 &\lesssim \frac{r^{K-1}}{\Delta_{\min}^2 p^{K}} \onormSize{}{\tX - \hat \tX}_F^2 \lesssim \frac{r^{K-2}}{p^{K-1} \Delta_{\min}^2}  \lesssim \frac{1}{r \log p},
\end{align}
and 
\begin{equation}
     L(\hat z) \lesssim \frac{\Delta_{\min}^2}{r \log p}.
\end{equation}
Above result holds for $\hat \theta(i)$ after switching the parameters $\mX$ with $\hat \mX$ and switch $\mtheta$ with $\hat \mtheta$ in the proof.

Last, notice that for all $a \in [r]$
\begin{align}
    (1 - O(1)) \frac{p^2}{r^2}\onormSize{}{\mW_{:a}^T \mX - \hat \mW_{:a}^T \hat \mX}_F^2 
    &\leq  \onormSize{}{ \sum_{\hat z(i) = z(i) = a} (\theta(i) \mW_{:a}^T \mX  - \hat \theta(i) \hat \mW_{:a}^T \hat \mX)}_F^2 \leq \onormSize{}{\tX - \hat \tX}_F^2 \leq p r,
\end{align}
where the first inequality follows from the facts that  $\ell(\hat z, z) \lesssim \frac{1}{\log p}, |z^{-1}(a)| \asymp p /r$, 
\begin{align}
     | z^{-1}(a)| - C \frac{p}{r} \ell(\hat z, z)  \leq |\hat z^{-1}(a)| \leq | z^{-1}(a)| + C \frac{p}{r} \ell(\hat z, z), 
\end{align}
\begin{equation}
  | z^{-1}(a)|  - C \frac{p}{r} \ell( \hat z,z)  \leq \sum_{ z(i) = z(i) = a} \theta(i) \leq | z^{-1}(a)|,  \quad 
    |\hat z^{-1}(a)|  - C \frac{p}{r} \ell(\hat z, z)  \leq \sum_{\hat z(i) = z(i) = a} \hat \theta(i) \leq |\hat z^{-1}(a)|.
\end{equation}

Hence, for all $i \in [p]$
\begin{align}
  (\theta(i) - \hat \theta(i))^2 \onormSize{}{ \mW_{:a}^T \mX }_F^2 - O(p) 
  &\leq  \onormSize{}{ (\theta(i) - \hat \theta(i))\mW_{:a}^T \mX }_F^2 - \onormSize{}{ \hat \theta(i) (\mW_{:a}^T \mX - \hat \mW_{:a}^T \hat \mX) }_F^2 \\
  &\leq \onormSize{}{\tX - \hat \tX}_F^2 \leq p r,\label{eq:theta_hat}
\end{align}
where the first inequality follows from $\onormSize{}{\mW_{:a}^T \mX - \hat \mW_{:a}^T \hat \mX}_F^2 \lesssim 1/p$ and $\hat \theta(i) \lesssim \frac{p}{r}$. Notice that for all $a \in [r]$
\begin{equation}
    \onormSize{}{\mW_{:a}^T \mX}_F^2 \geq \onormSize{}{\mS_{a:}}_F^2 \lambda_{\min}^{2(K-1)}(\mTheta \mM) \gtrsim p^{K-1}.
\end{equation}
The inequality indicates that $\theta(i)^2 = (1 + o(p^{-(K-2)})) \hat \theta(i)^2$. 

\end{proof}

\end{document}